\documentclass{cmslatex}
\usepackage[paperwidth=7in, paperheight=10in, margin=.875in]{geometry}
\usepackage[backref,colorlinks,linkcolor=red,anchorcolor=green,citecolor=blue]{hyperref}
\usepackage{amsfonts,amssymb,color}
\usepackage{amsmath}
\usepackage{graphicx}
\usepackage{cite}
\usepackage{enumerate}

\usepackage{multirow}
\usepackage{mathtools}

\allowdisplaybreaks

\newcommand{\x}{\mathbf{x}}

\renewcommand{\theequation}{\arabic{section}.\arabic{equation}}

   \allowdisplaybreaks
\begin{document}
 \title{Stabilized exponential time differencing schemes for the convective Allen-Cahn equation
 \thanks{Y. Cai's work was partially supported by National Natural Science Foundation of China grants 12171041 and 11771036.
 L. Ju's work was partially supported by US National Science Foundation grants DMS-1818438 and DMS-2109633. J. Li's work was partially supported by China Postdoctoral Science Foundation grant 2021M700476.}
 }
\author{
Yongyong Cai \thanks{Laboratory of Mathematics and Complex Systems and School of Mathematical Sciences, Beijing Normal University, Beijing 100875, China ({\tt yongyong.cai@bnu.edu.cn)}}
\and  Lili Ju \thanks{Department of Mathematics, University of South Carolina, Columbia, SC 29208, USA ({\tt ju@math.sc.edu})}
\and Rihui Lan \thanks{Department of Mathematics, University of South Carolina, Columbia, SC 29208, USA ({\tt rlan@mailbox.sc.edu})}
\and 
Jingwei Li \thanks{Laboratory of Mathematics and Complex Systems and School of Mathematical Sciences, Beijing Normal University, Beijing 100875, China  ({\tt jingwei@bnu.edu.cn)}}
}
 		\pagestyle{myheadings} \markboth{STABILIZED ETD SCHEMES FOR CONVECTIVE ALLEN-CAHN EQUATION}{Y. CAI, L. JU, R. LAN, AND J. LI}\maketitle
		
          \begin{abstract}
               The convective Allen-Cahn equation has been widely used to simulate multi-phase flows in many phase-field models. As a generalized form of the classic Allen-Cahn equation, the convective Allen-Cahn equation still preserves the maximum bound principle (MBP) in the sense that the time-dependent solution of the equation with appropriate initial and boundary conditions preserves for all time a uniform pointwise bound in absolute value. In this paper, we develop efficient first- and second-order exponential time differencing (ETD) schemes combined with the linear stabilizing technique to preserve the MBP unconditionally in the discrete setting. The space discretization is done using the upwind difference scheme for the convective term and the central difference scheme for the diffusion term, and both the mobility and nonlinear terms are approximated through the linear convex interpolation. The unconditional preservation of the MBP of the proposed schemes is proven, and their convergence analysis is presented. Various numerical experiments in two and three dimensions are also carried out to verify the theoretical results.
          \end{abstract}

\section{Introduction}\label{intro}
In this paper, we consider the convective Allen-Cahn equation taking the following form
\begin{align}\label{gac}
\partial_tu(\x,t)+\mathbf{v}(\x,t)\cdot \nabla u(\x,t)=M(u(\x,t))\left(\epsilon^2\Delta u(\x,t)+f(u(\x,t))\right),\quad \x \in \Omega, t >0,
\end{align}
subjected to suitable boundary conditions, e.g., periodic, homogeneous Neumann, or Dirichlet boundary conditions. Here, $\Omega\subset\mathbb{R}^{d} (d=1,2,3)$ is an open, connected, and bounded domain, $\mathbf{x}$ is the spatial variable, $t\ge 0$ is time,  $u(\x,t)$ is the unknown scalar-valued phase variable, $\epsilon>0$ is the parameter related to the interface width between two phases, $M(u)\ge 0$ is  the mobility function, $f(u)$ is the nonlinear reaction term, and $\mathbf{v}(\x,t)$ is a given velocity field satisfying the divergence-free condition $\nabla\cdot \mathbf{v}=0$. In addition, we assume that  $\mathbf{v}$ is differentiable, and $f(\cdot)$ and $M(\cdot)$ are locally $C^1(\mathbb{R})$ functions.

Compared with the classical Allen-Cahn equation \cite{AlCa79}, the convective Allen-Cahn equation \eqref{gac} is more complicated due to the presence of the velocity field. Nevertheless, the convective Allen-Cahn equation still preserves the maximum bound principle (MBP) in the sense that if the initial value and/or the boundary data are pointwisely bounded by a specific constant in absolute value, the solution is also bounded by the same constant everywhere for all time. It is well-known that the classic Allen-Cahn equation (the equation \eqref{gac}   without the convective term) with the periodic or homogeneous Neumann boundary condition can be regarded as the $L^2$ gradient flow with respect to the following energy functional
\begin{align}\label{energyorg}
E[u]=\frac{\epsilon^2}{2}(\nabla u(\x,t),\nabla u(\x,t))+(F(u(\x,t)),1),
\end{align}
where $F(\cdot)$ is the primitive fucntion of $-f$ (i.e., $f=-F'$) and $(\cdot,\;\cdot)$ represents the  inner product on $L^2(\Omega)$ with the corresponding norm denoted as $\|\cdot\|_0$. However,  the solution of the convective Allen-Cahn equation \eqref{gac}  is usually not guaranteed to decrease the energy $E[u]$ along the time.

The MBP is an important property of the convective Allen-Cahn equation \eqref{gac}, especially for the logarithmic type nonlinearity $F(\cdot)$ and  the degenerate mobility $M(\cdot)$. In such cases, if the numerical scheme does not satisfy MBP,  complex values may occur in numerical simulations due to the logarithm arithmetic, and/or the nonlinear mobility could become negative, which would lead to unphysical numerical solutions. For the classical Allen-Cahn equation,   the MBP-preserving numerical methods have been thoroughly studied. For the spatial discretizations, a partial list includes the works on the lumped-mass finite element method \cite{XFH19,XDF20}, finite difference method \cite{YDZ18}, finite volume method \cite{PGYF19,PGF19} and so on. For the temporal integration, the stabilized linear semi-implicit schemes were shown to preserve the MBP unconditionally for the first-order scheme \cite{TaYa16,XFY17,XHF20} and conditionally for the second-order scheme. Some nonlinear second-order schemes were presented to preserve the MBP for the Allen-Cahn type equations in \cite{LTZ20,ShXu20}. However, the works on the construction of numerical schemes for unconditional MBP preservation are rare. Shen et al.  analyzed the stabilizing MBP-preserving schemes with the finite difference discretization in space for the convective 
Allen-Cahn equation and designed an adaptive time-stepping scheme in \cite{STY16}, but unconditionally MBP-preserving schemes with order higher than one are still lacking. Thus, it is highly desirable to find alternative ways to develop higher-order time integration schemes which preserve the MBP unconditionally.

Combined with the linear stabilizing technique  \cite{XuTa06}, the exponential time differencing (ETD) method has been recently proposed to preserve the MBP.  
The ETD schemes are based on the variation-of-constants formula,   where the nonlinear terms are approximated by polynomial interpolations in time, followed by the exact integration of the resulting temporal integrals \cite{BKV98,CoMa02,HJW19,HoOs05}. Therefore, the ETD method is applicable to a large class of semilinear parabolic equations, especially for those with a stiff linear part \cite{HJW19,LJM19,JLQZ18,JZD15,JZZD15,wangc1,wangc2}. The first and second-order stabilized ETD schemes were also applied to the nonlocal Allen-Cahn equation and were proven to be unconditionally MBP-preserving in \cite{DJLQ19}. 
Then, an abstract framework on the MBP-preserving stabilized ETD schemes was established in \cite{DJLQ21} for a class of  semilinear parabolic equations and later was also applied to the conservative Allen-Cahn equations in \cite{LJCF21,JJLL21}. Recently, some third- and fourth-order MBP-preserving schemes were proposed for the Allen-Cahn equation by considering the integrating factor Runge-Kutta  schemes (IFRK)
\cite{IGG18,LLJF21,JLQY20,zhangh1,JLQ22}, and an arbitrarily high-order ETD multistep method was presented in \cite{zhouz1} by enforcing the maximum bound via extra cutoff postprocessing.

The main contribution of the current paper is to develop the unconditional MBP-preserving ETD schemes with first- and second-order temporal accuracies for the convective Allen-Cahn equation. Following the framework in \cite{DJLQ21}, we firstly reformulate the equation by introducing a stabilizing constant and approximate the nonlinear mobility function and nonlinear term with the convex linear interpolation in time.  With the upwind scheme for the convective term, the fully-discrete ETD scheme of the convective Allen-Cahn equation satisfies the MBP unconditionally. 

The rest of the paper is organized as follows. Section \ref{sec:pre} reviews the basic assumptions on the mobility and nonlinear potential function so that the model equation has a unique solution and satisfies the MBP. In Section \ref{sec:etd}, the first- and second-order ETD schemes for time integration are constructed and shown to preserve the discrete MBP unconditionally. Next we discuss the fully-discrete schemes with the upwind spatial discretization for the convective term. Then the error analysis is carried out in Section \ref{sec-conv}. A variety of numerical tests are performed to verify the theoretical results in Section \ref{sec:simul}, and finally we end this paper with some concluding remarks in  Section \ref{con}. 

\section{Preliminaries and maximum bound principle}\label{sec:pre}
In this section, we first review the convective Allen-Cahn equation \eqref{gac} and the assumptions on the mobility function and the nonlinear function under which the MBP holds. 
Suppose that  the initial value   of the equation \eqref{gac}  is given as
\begin{align*}
u(\x,0)=u_0(\x),\quad\x\in\overline{\Omega},
\end{align*}
where $u_0\in C(\overline\Omega)$.
We also impose the periodic boundary condition (only for a regular rectangular domain $\Omega=\prod_{i=1}^d(a_i,b_i)$), the homogeneous Neumann boundary condition or the following Dirichlet boundary condition
\begin{align}\label{eq:diri}
u(\x,t)=g(\x,t),\quad \x\in\partial\Omega,\, t>0,
\end{align} 
where $g\in C(\partial\Omega\times[0,\infty))$.
In practical applications, the nonnegative mobility function $M(\cdot)\ge0$ could lead to a degenerate parabolic equation. To avoid the technical difficulties arising from the possible degeneracy of $M(\cdot)$, if unspecified in the subsequent discussion, we always assume there exists a positive constant $\gamma>0$ such that
\begin{equation}\label{assump:mb}
M(s)\ge \gamma,\quad \forall\, s\in\mathbb{R}.
\end{equation}
Using classical theory for quasi-linear parabolic equations \cite{ladyzenskaya1967linear,taylor117partial}, 
under the condition \eqref{assump:mb}, the convective Allen-Cahn equation \eqref{gac}  admits a unique smooth solution. 

Following  the analysis for the MBP of semilinear parabolic equations \cite{DJLQ21}, we make the following assumption on the nonlinear function $f(u)$.
\begin{assumption}\label{assump-nonlinear}
	There exists a constant $\beta>0$ such that for the nonlinear function $f(u)$, it holds
	\begin{align*}
	f(\beta)\leq 0\leq f(-\beta).
	\end{align*}
\end{assumption}

Let $\|\cdot\|$ denote the maximum norm on $C(\overline{\Omega})$ and $\|\cdot\|_{\partial\Omega}$ on $C(\partial\Omega)$.
Under the condition \eqref{assump:mb} and Assumption \ref{assump-nonlinear}, for a smooth mobility function $M(\cdot)$, it holds that the convective Allen-Cahn equation \eqref{gac} satisfies the MBP \cite{Evans97,STY16}: if the initial value and/or the boundary data are  pointwisely bounded by the constant $\beta$ in absolute value, then its solution is also bounded by $\beta$ for all time, i.e., for the periodic or homogeneous Neumann boundary condition case,
\begin{equation}\label{eq:mbp-a}
\|u_0\|\leq \beta\quad \Longrightarrow \quad |u(\x, t)| \leq \beta, \;\;\; \forall\,\x\in\overline{\Omega},\, t>0.
\end{equation}
and for the Dirichlet boundary condition case,
\begin{equation}\label{eq:mbp-b}
\|u_0\|\leq \beta\;\; \& \;\; \|g(\cdot,t)\|_{\partial\Omega}\leq \beta,\; \;\forall\, t>0 \quad \Longrightarrow \quad |u(\x, t)| \leq \beta, \;\;\; \forall\,\x\in\overline{\Omega},\, t>0.
\end{equation}

\begin{remark}
We would like to note that by a standard approximation argument \cite{EG96,DD16}, it can be proved that for the case of general mobility $M(\cdot)\ge0$,  weak solutions to the convective Allen-Cahn equation \eqref{gac} also exist and satisfy the MBP in the a.e. sense.
\end{remark}
\begin{remark}\label{rmk} 
In view of the MBP by \eqref{eq:mbp-a} or \eqref{eq:mbp-b}, the condition  \eqref{assump:mb} on the mobility $M(u)$ can be relaxed to
	\begin{equation}
	M(s)\ge\gamma,\,\quad \forall\, |s|\leq\beta.
	\end{equation}
	In addition, it is  required that $M(\cdot)$ is  differentiable  when restricted in the interval $[-\beta,\beta]$.
\end{remark}

For a given scalar function $\phi(\x)$ and a divergence-free vector field $\mathbf{w}(\x)$, let us define a linear elliptic differential operator $\mathcal{L}_{\mathbf{w}}[\phi]$ as
\begin{equation}\label{eq:ldff}
\mathcal{L}_{\mathbf{w}}[\phi] u =\epsilon^2M(\phi(\x))\Delta u-\mathbf{w}(\x)\cdot\nabla u,
\end{equation}
and a modified nonlinear function $\widetilde{f}$ as
\begin{equation}\label{eq:tlf}
\widetilde{f}(u)=M(u)f(u).
\end{equation}
According to the fact that the mobility is nonnegative, Assumption \ref{assump-nonlinear} implies that $\widetilde{f}(\beta)\leq 0\leq\widetilde{f}(-\beta)$. For the above  operator $\mathcal{L}_{\mathbf{w}}[\phi]$, the following lemma holds \cite{DJLQ21}.
\begin{lemma}\label{lem-linear}
	Given $\phi(\x)\in C(\overline{\Omega})$ and $\mathbf{w}(\x)\in C(\overline{\Omega},\mathbb{R}^d)$, 
	then the second-order  linear elliptic differential operator $L_{\mathbf{w}}[\phi]$  under appropriate boundary condition (periodic, homogeneous Neumann or homogeneous Dirichlet) generates a contraction semigroup (denoted as $e^{t\mathcal{L}_{\mathbf{w}}[\phi]}$) with respect to the supremum norm on the subspace of $C(\overline{\Omega})$ that satisfies the corresponding boundary condition, i.e.,
	\begin{align*}
	\|e^{t\mathcal{L}_{\mathbf{w}}[\phi]}\|\leq 1,\quad \forall\, \phi\in C(\overline{\Omega}), \,t\geq 0.
	\end{align*}
	Furthermore, for any $\alpha>0$, we have
	\begin{align*}
	\|e^{t(\mathcal{L}_{\mathbf{w}}[\phi]-\alpha)}\|\leq e^{-\alpha t},\quad \forall\, \phi\in C(\overline{\Omega}), \,t\geq 0.
	\end{align*}
\end{lemma}

The convective Allen-Cahn equation  \eqref{gac} can be re-written as
\begin{align}\label{prob}
\partial_tu(\x,t)=\mathcal{L}_{\mathbf{v}(t)}[u(t)]u(\x,t)+\widetilde{f}(u(\x,t)), \quad \x\in\Omega,\, t>0.
\end{align}
Here we have omitted the dependence on $\x$ in $\mathcal{L}_{\mathbf{v}(t)}[u(t)]$ when there is no confusion, i.e., ${\mathbf{v}(t)} = {\mathbf{v}(\x,t)}$
and ${u(t)} = {{u}(\x,t)}$ are functions of $\x$.
In this paper, we mainly focus on establishing the MBP for finite difference approximations of the convective Allen-Cahn equation \eqref{gac} under Assumption \ref{assump-nonlinear} and the condition \eqref{assump:mb} for mobility. Two typical nonlinear potentials $F(u)$ with $f(u)=-F'(u)$ will be numerically studied, i.e. 
the double-well potential
\begin{align}\label{dw}
F(u)=\frac{1}{4}\left(u^{2}-1\right)^{2},
\end{align}
where {the corresponding maximum bound} $\beta\in[1,\infty)$ (c.f. Assumption \ref{assump-nonlinear}), and the Flory-Huggins potential
\begin{align}\label{fh}
F(u)=\frac{\theta}{2}[(1+u) \ln (1+u)+(1-u) \ln (1-u)]-\frac{\theta_{c}}{2}u^{2},
\end{align}
with  $ \theta_c>\theta>0$, where {the corresponding maximum  bound} $\beta\in[\rho,1)$ (c.f. Assumption \ref{assump-nonlinear})  with $\rho$ being the positive root of $f(\rho)=0$. 

For the numerical solution purpose, by introducing a stabilizing constant $\kappa>0$, the equation (\ref{prob}) can be further formulated into an equivalent form 
\begin{align}\label{stab}
\partial_tu(\x,t)=\mathcal{L}^\kappa_{\mathbf{v}(t)}[u(t)]u(\x,t)+N(u(\x,t)),\,\quad\forall\, \x\in\Omega,\, t>0,
\end{align}
where
\begin{align}\label{stab-linear}
\mathcal{L}^\kappa_{\mathbf{v}(t)}[u(t)]=\mathcal{L}_{\mathbf{v}(t)}[u(t)]-\kappa\mathcal{I},\quad N(u)=\kappa u+\widetilde{f}(u).
\end{align}
The stabilizing constant $\kappa$ is required that 
\begin{align}\label{coef}
\kappa\geq \max_{|\xi|\leq \beta}\left|\widetilde{f}'(\xi)\right|=\max_{|\xi|\leq \beta}\left|M'(\xi)f(\xi)+M(\xi)f'(\xi)\right|.
\end{align}
Note that \eqref{coef} is always  well-defined as long as $M(\xi)$ and $f(\xi)$ are continuously differentiable.

\begin{lemma}\label{lemma-stab-nonlinear}
	Under Assumption \ref{assump-nonlinear} and the choice of the stabilizing constant \eqref{coef}, we have
	
	\rm{(I)} $|N(\xi)|\leq\kappa\beta$ for any $\xi\in[-\beta,\beta]$;
	
	\rm{(II)} $|N(\xi_1)-N(\xi_2)|\leq 2\kappa|\xi_1-\xi_2|$ for any $\xi_1,\xi_2\in[-\beta,\beta]$.
\end{lemma}
\begin{proof}
	From \eqref{stab-linear}, we have $N'(\xi)=\kappa+\widetilde{f}'(\xi)$ and 
	$0\leq N'(\xi)\leq2\kappa$ 
	which implies (II). In addition, $N(\xi)$ is an increasing function.  Recalling Assumption \ref{assump-nonlinear}, for any $\xi\in[-\beta,\beta]$, we get 
	\[
	-\kappa\beta\leq-\kappa\beta+\widetilde{f}(-\beta)=N(-\beta)\leq N(\xi)\leq N(\beta)=\kappa\beta+\widetilde{f}(\beta)\leq\kappa\beta,
	\]
	and  the conclusion (I) holds.\hfill
\end{proof}

\section{Exponential time differencing for temporal approximation}\label{sec:etd}
In this section, we construct the temporal approximations of the convective Allen-Cahn equation \eqref{prob} by using the ETD approach. First of all, let us define the $\phi$-functions as follows: 
\begin{align*}
\phi_0(s)=e^{s},\quad\phi_1(s)=\frac{e^{s}-1}{s},\quad\phi_2(s)=\frac{e^{s}-s-1}{s^2}.
\end{align*}
Assuming that Assumption \ref{assump-nonlinear} and the condition \eqref{assump:mb} for mobility hold, we will begin with the equivalent form \eqref{stab} to develop the MBP-preserving ETD schemes of first- and second-order in this section.

\subsection{ETD schemes and discrete MBPs}
Choosing a time step size $\tau>0$,  we have the time steps as  $\{t^n=n\tau\}_{n\geq0}$. To construct the ETD schemes for time stepping of the model problems \eqref{gac}, below we take  the Dirichlet boundary condition case for illustration and focus on the stabilized form \eqref{stab} on the interval $[t^n,t^{n+1}]$, or equivalently, the equation for $w(s):=w(\x,s)=u(\x,t^n+s)$:
\begin{align}\label{etd}
\left\{
\begin{array}{ll}
\partial_sw(\x,s)=\mathcal{L}^\kappa_{\mathbf{v}(t^n+s)}[w(s)] w(\x,s)+N(w(\x,s)),&\quad \x\in \Omega, s\in (0,\tau],\\[4pt]
w(\x,s)=g(\x,t^n+s),&\quad \x\in \partial\Omega, s\in (0,\tau],\\[4pt]
w(\x,0)=u(\x,t^n),&\quad\x\in \Omega.
\end{array}
\right.
\end{align}

Let us start with $u^0 = u_0$ and denote 
$u^n$ as an approximation of $u(t_n)$. 

For the first-order-in-time approximation of (\ref{etd}), we use $\mathcal{L}^\kappa_{\mathbf{v}(t^n+s)}[u(t^n+s)]=\mathcal{L}^\kappa_{\mathbf{v}(t^n)}[u(t^n)]+O(\tau)$ and $N(u(t^n+s))=N(u(t^n))+O(\tau)$ to construct the first order ETD (ETD1) scheme: for $n\geq0$ and given $u^n$, find $u^{n+1}=w^n(\tau)$ by solving 
\begin{align}\label{etd1}
\left\{
\begin{array}{ll}
\partial_sw^n(\x,s)=\mathcal{L}^\kappa_{\mathbf{v}^n}[u^n(\x)] w^n(\x,s)+N(u^n(\x)),
&\quad \x\in \Omega, s\in (0,\tau],\\[4pt]
w^n(\x,s)=g(\x,t^n+s),&\quad \x\in \partial\Omega, s\in (0,\tau],\\[4pt]
w^n(\x,0)=u^n(\x),&\quad\x\in \Omega,
\end{array}
\right.
\end{align}
where  $\mathbf{v}^n=\mathbf{v}(t^n)$. Equivalently,
\begin{align*}
u^{n+1}(\x)=\phi_0(\tau\mathcal{L}_{\mathbf{v}^n}^\kappa[u^n(\x)])u^n(\x)+\tau\phi_1(\tau\mathcal{L}_{\mathbf{v}^n}^\kappa[u^n(\x)])N(u^n(\x)).
\end{align*}
We refer to \cite{DJLQ21} for efficient calculations of the products of matrix exponentials and vectors needed in numerical implementation of the proposed ETD1 scheme.

\begin{theorem}\label{thm-etd1-mbp}
	{\rm(Discrete MBP of the ETD1 scheme)} The ETD1 scheme  \eqref{etd1} preserves the discrete MBP unconditionally, i.e., for any time step size $\tau>0$, the ETD1 solution   satisfies $\|u^n\|\leq\beta$ for any $n>0$ if $\|u_0\|\leq \beta$ and  $ \|g(\cdot,t)\|_{\partial\Omega}\leq \beta$ for any $t>0$.
\end{theorem}
\begin{proof} 
	Since $\|u^0\|\leq \beta$, by using mathematical induction we  only need to show that $\|u^n\|\leq\beta$ implies
	$|u^{n+1}|\leq\beta$ for  $n\ge0$. First,  we know $u^n,w^n(s)\in C(\overline{\Omega})$ ($n\ge0$, $s\in[0,\tau]$) by standard PDE analysis for linear parabolic equations. Suppose that there exists $(\mathbf{x}^*,s^*)\in \overline{\Omega}\times[0,\tau]$ such that
	\begin{align*}
	w^n(\mathbf{x}^*,s^*)=\max_{0\leq s\leq\tau, \x\in\overline{\Omega}} w^n(\x,s).
	\end{align*}
	If $\x^*\in\partial\Omega$ or $s^*=0$, the boundary data $|g(\x,t)|\leq \beta$ for any $\x\in\partial{\Omega},\, t>0$ or the initial condition $\|u^n\|\leq\beta$ at $s^*=0$ implies $w^n(\mathbf{x}^*,s^*)\leq\beta$. Otherwise, $\x^*\in\Omega$ and $s^*\in(0,\tau]$, we have
	\begin{align*}
	\partial_sw^n(\mathbf{x}^*,s^*)\geq 0,\quad \nabla w^n(\mathbf{x}^*,s^*)=\mathbf{0},\quad \Delta w^n(\mathbf{x}^*,s^*) \leq 0,
	\end{align*}
	which implies
	\begin{align*}
	\mathcal{L}_{\mathbf{v}^n}[u^n]w^n(\mathbf{x}^*,s^*)\leq0,
	\end{align*}
	then
	\begin{equation*}
	\kappa w^n(\mathbf{x}^*,s^*)\leq N(u^n(\mathbf{x}^*,s^*)).
	\end{equation*}

	Since $\|u^n\|\leq\beta$,  Lemma \ref {lemma-stab-nonlinear} leads to  $N(u^n(\mathbf{x}^*,s^*))\leq \kappa \beta$, furthermore $w^n(\mathbf{x}^*,s^*)\leq\beta$. Similarly, if there exists $(\mathbf{x^{**}},s^{**})\in \Omega\times[0,T]$ such that 
	\begin{align*}
	w^n(\mathbf{x}^{**},s^{**})=\min_{0\leq s\leq\tau,\x\in\overline{\Omega}}w(\x,s),
	\end{align*}
	one can show $w^n(\mathbf{x}^{**},s^{**})\geq-\beta$. Thus, we conclude $\|u^{n+1}\|\leq\beta$.\hfill
\end{proof}

Next, we consider the second-order-in-time approximation of \eqref{etd} by using
\begin{align}\label{linear_interp}
\begin{aligned}
\mathcal{L}^\kappa_{\mathbf{v}(t^n+s)}[u(t^n+s)]&=\left(1-\frac{s}{\tau}\right)\mathcal{L}^\kappa_{\mathbf{v}^n}[u(t^n)]+\frac{s}{\tau}\mathcal{L}^\kappa_{\mathbf{v}^{n+1}}[u(t^{n+1})]+O(\tau^2),\\
N(u(t^n+s))&=\left(1-\frac{s}{\tau}\right)N(u(t^n))+\frac{s}{\tau}N(u(t^{n+1}))+O(\tau^2).
\end{aligned}
\end{align}
To construct a second-order-in-time approximation, we need to make the approximation of the differential operator part $\mathcal{L}^\kappa_{\mathbf{v}(t_n+s)}[u(t^n+s)]$ to be $s$-independent, so that ETD approach could be applied. {In order to do this,  we learn from \eqref{linear_interp} and reformulate \eqref{etd} as}
\begin{align*}
\partial_s w(\x,s)=\;&\frac12\left(\mathcal{L}^\kappa_{\mathbf{v}^n}[u(\x,t^n)]+\mathcal{L}^\kappa_{\mathbf{v}^{n+1}}[u(\x,t^{n+1})]\right)w(\x,s)\\
&{+\left(1-\frac{s}{\tau}\right)N(u(\x,t^n))+\frac{s}{\tau}N(u(\x,t^{n+1}))}\\
&+\left(\frac{s}{\tau}-\frac12\right)\left(\mathcal{L}^\kappa_{\mathbf{v}^{n+1}}[u(\x,t^{n+1})]-\mathcal{L}^\kappa_{\mathbf{v}^n}[u(\x,t^n)]\right)w(\x,s)+O(\tau^2).
\end{align*}
Integrating both sides of the above equation from $0$ to $\tau$ and applying  the midpoint rule to the fourth term on the RHS, we derive the second order ETD Runge-Kutta (ETDRK2) scheme: for $n\geq0$ and given $u^n$, find $u^{n+1}=w^n(\tau)$ by solving 
\begin{align}\label{etd2}
\left\{
\begin{array}{ll}
\partial_sw^n(\x,s)=\frac{1}{2}\left(\mathcal{L}_{\mathbf{v}^n}^\kappa[u^n(\x)]+\mathcal{L}_{\mathbf{v}^{n+1}}^\kappa[\hat{u}^{n+1}(\x)]\right)w^n(\x,s)\\[4pt]
\qquad\qquad\qquad\qquad+(1-\frac{s}{\tau})N(u^n(\x))+\frac{s}{\tau}N(\hat{u}^{n+1}(\x)),
&\quad\x\in \Omega,s\in (0,\tau],\\[4pt]
w(\x,s)=g(\x,t^n+s),&\quad\x\in \partial\Omega,s\in (0,\tau],\\[4pt]
w^n(\x,0)=u^n,&\quad\x\in \Omega,
\end{array}
\right.
\end{align}
where $\hat{u}^{n+1}$ is generated from the ETD1 scheme \eqref{etd1}. 
Equivalently,
\begin{eqnarray}
\hat{u}^{n+1}(\x)&=&\phi_0(\tau\mathcal{L}_{\mathbf{v}^n}^\kappa[u^n])u^n(\x)+\tau\phi_1(\tau\mathcal{L}_{\mathbf{v}^n}^\kappa[u^n(\x)])N(u^n(\x)),\notag\\
u^{n+1}(\x)&=&\phi_0\left(\frac{\tau}{2}(\mathcal{L}_{\mathbf{v}^n}^\kappa[u^n(\x)]+\mathcal{L}_{\mathbf{v}^{n+1}}^\kappa[\hat{u}^{n+1}(\x)])\right)u^n(\x)\notag\\
&&+\tau\phi_1\left(\frac{\tau}{2}(\mathcal{L}_{\mathbf{v}^n}^\kappa[u^n(\x)]+\mathcal{L}_{\mathbf{v}^{n+1}}^\kappa[\hat{u}^{n+1}(\x)])\right)u^n(\x)\notag\\
&&+\tau\phi_2\left(\frac{\tau}{2}(\mathcal{L}_{\mathbf{v}^n}^\kappa[u^n(\x)]+\mathcal{L}_{\mathbf{v}^{n+1}}^\kappa[\hat{u}^{n+1}(\x)])\right)\left(N(\hat{u}^{n+1}(\x))-N(u^n(\x))\right).\notag
\end{eqnarray}
It is worth noting that both ETD1 and ETDRK2 are linear schemes.
\begin{remark}
In view of the linear interpolation \eqref{linear_interp} for the differential operator part $\mathcal{L}^\kappa_{\mathbf{v}(t_n+s)}[u(t^n+s)]$, it is easy to verify that the midpoint rule is the exact evaluation for the ETDRK2 scheme \eqref{etd2}.
\end{remark}
\begin{theorem}\label{thm-etd2-mbp}
	{\rm(Discrete MBP of the ETDRK2 scheme)} The ETDRK2 scheme  \eqref{etd2} preserves the discrete MBP unconditionally, i.e., for any time step size $\tau>0$, the ETDRK2 solution  satisfies $\|u^n\|\leq\beta$ for any $n>0$ if $\|u_0\|\leq \beta$ and  $ \|g(\cdot,t)\|_{\partial\Omega}\leq \beta$ for any $ t>0$.
\end{theorem}
\begin{proof}
	The proof is quite similar to the ETD1 case and thus is omitted here for brevity. 
	\hfill
\end{proof}

\begin{remark}
For the case of periodic or homogeneous boundary condition, the ETD1 and ETDRK2 schemes can be similarly formulated as \eqref{etd1} and \eqref{etd2} by removing $w(\x,s)=g(\x,t^n+s)$ for $\x\in \partial\Omega,s\in (0,\tau]$ and imposing the respective boundary condition. Consequently, Theorems \ref{thm-etd1-mbp} and  \ref{thm-etd2-mbp} for the MBP still hold  by removing the extra boundary value requirement
$ \|g(\cdot,t)\|_{\partial\Omega}\leq \beta$ as done in \cite{DJLQ21}.
\end{remark}

\subsection{Fully discrete ETD schemes}\label{sec:fully-disc}
In this section, we focus on the spatial discretization. To this end, we recall the continuity of a function defined on a set $D\subset \mathbb{R}^d$ can be described as follows \cite{Rudin64}:
\[
w \text{~is continuous at~} \x^*\in D\Longleftrightarrow \forall\, \x_k\rightarrow \x^* \text{~in~} D \text{~implies~} w(\x_k)\rightarrow w(\x^*).
\]
Let $C(X)$ be the continuous function space over $X$, where $X$ is the set of all spatial grid points (boundary and interior points). 
Denote $U$ as the values of $u$ on $X$, i.e. $U(\x,t)=u(\x,t),\ \forall\, \x\in X, \;t>0$. 
The corresponding space-discrete equation of \eqref{prob} becomes an ordinary differential system taking the same form:
\begin{align}\label{prob-disc}
\partial_t U(\x,t)=\mathcal{L}_{\mathbf{v}(t),h}[U(t)]U(\x,t)+f(U(\x,t)),\quad t>0,\,\x\in X^*,
\end{align}
with $U(\x,0)=u_0(\x)$ for all $\x\in X$ where $X^*=X\cap \overline{\Omega}^+$ with $\overline{\Omega}^+=\prod\limits_{i=1}^d(a_i,b_i]$ for the periodic boundary condition case,  $X^*=X$ for the homogeneous Neumann boundary condition case,  and $X^*$ is the set of all interior grid points for the Dirichlet boundary condition case.

In order to establish the MBP for the discrete-in-space convective Allen-Cahn system \eqref{prob-disc}, as well as its time discretizations (to produce the fully discrete schemes)  proposed later, we make the following specific assumptions on the discrete operator $\mathcal{L}_{\mathbf{v}(t),h}[U(t)]$.
\begin{assumption}\label{assump-linear}
	For any given  $U\in C(X)$ and $\mathbf{v}\in (C(X))^d$, the operator $\mathcal{L}_{\mathbf{v},h}[U]$ satisfies that for any $W\in C(X)$ and $\x_0\in X^*$, if
	\begin{eqnarray*}
		W(\x_0)=\max_{\x\in X\cap\overline{\Omega}}W(\x),
	\end{eqnarray*}
	then $\mathcal{L}_{\mathbf{v},h}[U]W(\x)\big|_{\x=\x_0}\leq0$.
\end{assumption}

Next, we shall describe the concrete discrete scheme concerning the Dirichlet boundary condition for 1D ($\x=x\in\Omega$), and the extensions to  2D/3D rectangular domains and other boundary conditions are straightforward.
Given $\Omega=(a,b)$, the interval is divided into $N$ subintervals with uniform mesh size $h=(b-a)/N$ and  the grids points are $X=\{x_i=a+ih\}_{i=0}^{N}$ with $X^*=\{x_i,\,i=1,\ldots,N-1\}$ for the Dirichlet boundary condition. Let $U_i:=U_i(t)$  be the approximation of $u(x,t)$ at $x_i\in X$. We only need describe how  $\{U_i(t)\}_{i=1}^{N-1}$ evolve since the boundary values are explicitly given by $U_0(t)=g(x_0,t)$ and $U_N(t)=g(x_N,t)$. For convenience, we can view $U$ as a vector function of time $U=(U_1,U_2,\ldots,U_{N-1})^T$ when needed. 

The Laplace opertaor is discretized by the central finite difference method 
\[
u_{x x}\left(x_{i}, \cdot\right) \approx \frac{u(x_{i-1},\cdot)-2 u(x_{i},\cdot)+u(x_{i+1},\cdot)}{h^{2}},\quad i=1,2,\cdots,N-1,
\]
and the matrix $D_{h}$  as the discrete approximation of the Laplace operator under the homogeneous Dirichlet  boundary condition can be expressed as
\[
D_{h}=\frac{1}{h^{2}}\left[\begin{array}{ccccc}
-2 & 1 & & &\\
1 & -2 & 1 & & \\
& \ddots & \ddots & \ddots & \\
& & 1 & -2 & 1 \\
& & & 1 & -2
\end{array}\right]_{(N-1) \times (N-1)}.
\]
In addition, we need define $G_D= \frac{1}{h^{2}}(g(x_0,t),0,\cdots,0,g(x_N,t))^T$ for contribution of the boundary values to the Laplace operator. 

The convective term is  approximated by the upwind difference scheme
\[
\mathbf{v}u_{x}|_{x=x_i}\approx \mathbf{v}_i^{+} U_{i}^{-}+\mathbf{v}_i^{-} U_{i}^{+},
\]
where $\mathbf{v}_i^{+}$ and $\mathbf{v}_i^{-}$ are defined as
\[
\mathbf{v}_i^{+}=\max \{0, \mathbf{v}(x_i,t)\}, \quad \mathbf{v}_i^{-}=\min \{0, \mathbf{v}(x_i,t)\},
\]
and $U_{i}^{-}$ and $U_{i}^{+}$ are defined as
\[
U_{i}^{-}=\frac{-U_{i-1}+U_{i}}{h}, \quad U_{i}^{+}=\frac{U_{i+1}-U_{i}}{h}.
\]
Denote  $\mathbf{v}_i:=\mathbf{v}(x_i,t)$ ($i=1,\ldots,N-1$) and define
\[
A_{\mathbf{v}(t)}=\frac{1}{2 h}\left[\begin{array}{cccccc}
-2 & 1-\operatorname{sign}\left(\mathbf{v}_{1}\right) & &&\\
1+\operatorname{sign}\left(\mathbf{v}_{2}\right) & -2 & 1-\operatorname{sign}\left(\mathbf{v}_{2}\right) & \\
& \ddots & \ddots & \ddots \\
& & 1+\operatorname{sign}\left(\mathbf{v}_{N-2}\right) & -2 & 1-\operatorname{sign}\left(\mathbf{v}_{N-2}\right) \\
& & & 1+\operatorname{sign}\left(\mathbf{v}_{N-1}\right) & -2
\end{array}\right]
\]
and 
\[
\Lambda_{\mathbf{v}(t)}=\operatorname{diag}\left(\left|\mathbf{v}_{1}\right|, \ldots,\left|\mathbf{v}_{N-1}\right|\right), \quad \Lambda_U=\operatorname{diag}\left(M\left(U_{1}\right), \ldots, M\left(U_{N-1}\right)\right),
\]
where  $\operatorname{diag}(V)$ is the diagonal matrix with diagonal elements being components of the vector $V$.
We also define $G_C= \frac{1}{2h}((1+\operatorname{sign}(\mathbf{v}_{1}))g(x_0,t),0,\cdots,0,(1-\operatorname{sign}(\mathbf{v}_{N-1}))g(x_N,t))^T$ for the contribution of the boundary values to the convective term. 
Then  $\mathcal{L}_{\mathbf{v}(t),h}[U(t)]$ can be written in the matrix form as
\[
\mathcal{L}_{\mathbf{v}(t),h}[U(t)]= \epsilon^2\Lambda_{U(t)} D_{h}+\Lambda_{\mathbf{v}(t)}A_{\mathbf{v}(t)}.
\]

It is  easy to verify that at any fixed time $t^*$, the combination of  the upwind finite difference approximation of the convective term and the central difference approximation of the diffusion term  in $\mathcal{L}_{\mathbf{v}(t^*),h}[u(t^*)]$ satisfies Assumption \ref{assump-linear}.  The generated contraction semigroup can be simply treatd  as a matrix exponential $e^{t\mathcal{L}_{\mathbf{v}(t^*),h}[u(t^*)]}$.
The  discrete-in-space system of the stabilized formulation  \eqref{stab} then reads
\begin{align}\label{eq:spatial}
\partial_t U=\mathcal{L}_{\mathbf{v}(t),h}^\kappa [U(t)]U+(N(U)+ \epsilon^2\Lambda_UG_D+\Lambda_{\mathbf{v}(t)}G_C),\qquad t>0,
\end{align}
where $\mathcal{L}_{\mathbf{v}(t),h}^\kappa[U(t)]=\mathcal{L}_{\mathbf{v}(t),h}[U(t)]-\kappa I$.
The fully discrete ETD1 and ETDRK2 schemes then can be produced by applying ETD1 \eqref{etd1} and ETDRK2 \eqref{etd2} to \eqref{eq:spatial}, respectively.

In the following, we present the equivalent formulas of the above fully ETD1 and ETDRK2 schemes that are usually more convenient  for practical implementation. With the above 1D Dirichlet set-up, adopting the temporal discretizations in Section \ref{sec:etd},
we denote $U_i^n$ to be the numerical approximation of $u(x_i,t^n)$ ($i=1,\ldots,N-1$, $n\ge0$) and $U^n$ to be a function defined on $X^*$ (also viewed as the solution vector as $U^n=(U_1^n,\ldots,U_{N-1}^n)^T$).
 
The  fully discrete ETD1 scheme can be written as 
\begin{align}\label{etd1-fully}
U^{n+1}=\phi_0(\tau\mathcal{L}_{\mathbf{v}^n,h}^\kappa[U^n])U^n+\tau\phi_1(\tau\mathcal{L}_{\mathbf{v}^n,h}^\kappa[U^n])\widetilde{N}(U^n),
\end{align}
where $\mathbf{v}^n$ is understood as $\mathbf{v}(x,t^n)$ taking values at $x_i\in X^*$ and $\widetilde{N}(U^n) = N(U^n)+\epsilon^2\Lambda_{U^n}G^n_D+\Lambda_{\mathbf{v}^n}G^n_C)$.
The fully discrete ETDRK2 scheme reads  
\begin{eqnarray}\label{etd2rk-fully}
\hat{U}^{n+1}&=&\phi_0(\tau\mathcal{L}_{\mathbf{v}^n,h}^\kappa[U^n])U^n+\tau\phi_1(\tau\mathcal{L}_{\mathbf{v}^n,h}^\kappa[U^n])\widetilde{N}(U^n),\label{etd2rk-fully:1}\\
U^{n+1}&=&\phi_0\left(\frac{\tau}{2}(\mathcal{L}_{\mathbf{v}^n,h}^\kappa[U^n]+\mathcal{L}_{\mathbf{v}^{n+1},h}^\kappa[\hat{U}^{n+1}])\right)U^n\notag\\
&&+\tau\phi_1\left(\frac{\tau}{2}(\mathcal{L}_{\mathbf{v}^n,h}^\kappa[U^n]+\mathcal{L}_{\mathbf{v}^{n+1},h}^\kappa[\hat{U}^{n+1}])\right)\widetilde{N}(U^n)\notag\\
&&+\tau\phi_2\left(\frac{\tau}{2}(\mathcal{L}_{\mathbf{v}^n,h}^\kappa[U^n]+\mathcal{L}_{\mathbf{v}^{n+1},h}^\kappa[\hat{U}^{n+1}])\right)\left(\widetilde{N}(\hat{U}^{n+1})-\widetilde{N}(U^n)\right).\label{etd2rk-fully:2}
\end{eqnarray}
\begin{remark}
For the periodic boundary condition case, we have 
$X^*=\{x_i,\,i=1,\ldots,N\}$, $U_0=U_N$, $G_D=G_C={\bf 0}$, and 
\[
D_{h}=\frac{1}{h^{2}}\left[\begin{array}{ccccc}
-2 & 1 & & &1\\
1 & -2 & 1 & & \\
& \ddots & \ddots & \ddots & \\
& & 1 & -2 & 1 \\
1& & & 1 & -2
\end{array}\right]_{N \times N},
\]
\[
A_{\mathbf{v}(t)}=\frac{1}{2 h}\left[\begin{array}{cccccc}
-2 & 1-\operatorname{sign}\left(\mathbf{v}_{1}\right) & &&1+\operatorname{sign}\left(\mathbf{v}_{1}\right)\\
1+\operatorname{sign}\left(\mathbf{v}_{2}\right) & -2 & 1-\operatorname{sign}\left(\mathbf{v}_{2}\right) & \\
& \ddots & \ddots & \ddots \\
& & 1+\operatorname{sign}\left(\mathbf{v}_{N-1}\right) & -2 & 1-\operatorname{sign}\left(\mathbf{v}_{N-1}\right) \\
1-\operatorname{sign}\left(\mathbf{v}_{N}\right)& & & 1+\operatorname{sign}\left(\mathbf{v}_{N}\right) & -2
\end{array}\right],
\]
\[
\Lambda_{\mathbf{v}(t)}=\operatorname{diag}\left(\left|\mathbf{v}_{1}\right|, \ldots,\left|\mathbf{v}_{N}\right|\right), \quad \Lambda_U=\operatorname{diag}\left(M\left(U_{1}\right), \ldots, M\left(U_{N}\right)\right).
\]
For the homogeneous Neumann boundary condition case, we have 
$X^*=X=\{x_i,\,i=0,1,\ldots,N\}$, $G_D=G_C={\bf 0}$, and 
\[
D_{h}=\frac{1}{h^{2}}\left[\begin{array}{ccccc}
-2 & 2 & & &\\
1 & -2 & 1 & & \\
& \ddots & \ddots & \ddots & \\
& & 1 & -2 & 1 \\
& & & 2 & -2
\end{array}\right]_{(N+1) \times (N+1)},
\]
\[
A_{\mathbf{v}(t)}=\frac{1}{2 h}\left[\begin{array}{cccccc}
-2 & 2 & &&\\
1+\operatorname{sign}\left(\mathbf{v}_{1}\right) & -2 & 1-\operatorname{sign}\left(\mathbf{v}_{1}\right) & \\
& \ddots & \ddots & \ddots \\
& & 1+\operatorname{sign}\left(\mathbf{v}_{N-1}\right) & -2 & 1-\operatorname{sign}\left(\mathbf{v}_{N-1}\right) \\
& & & 2 & -2
\end{array}\right],
\]
\[
\Lambda_{\mathbf{v}(t)}=\operatorname{diag}\left(\left|\mathbf{v}_{0}\right|, \ldots,\left|\mathbf{v}_{N}\right|\right), \quad \Lambda_U=\operatorname{diag}\left(M\left(U_{0}\right), \ldots, M\left(U_{N}\right)\right).
\]
\end{remark}

Similar to Theorems \ref{thm-etd1-mbp} and \ref{thm-etd2-mbp}, we then obtain the unconditional MBP preservation for the above fully discrete ETD1 and ETDRK2 schemes as: 
\begin{equation}\label{fullmbp}
\|U^n\|_{X}\leq \beta,\quad \forall\, n\ge0, 
\end{equation}
where $\|\cdot\|_{X}$ denotes the maximum norm on $C(X)$.

\section{Convergence analysis}\label{sec-conv}
As an important application of the MBP, we now consider the convergence of the ETD1 \eqref{etd1-fully} and ETDRK2 \eqref{etd2rk-fully:1}-\eqref{etd2rk-fully:2} schemes. For simplicity of the analysis, we again take  the 1D problem with $\Omega=(a,b)$ and set the homogeneous Dirichlet  boundary condition. In this case,  $u(x,t)=g(x,t) = 0$ for $x\in\partial\Omega,\;t>0$, thus $U_0(t) = U_N(t) = 0$, $G_D=G_C={\bf 0}$ and $\widetilde{N}(U^n) = {N}(U^n)$ in \eqref{etd1-fully} and \eqref{etd2rk-fully:1}-\eqref{etd2rk-fully:2}.  The error estimation results can be naturally extended to the problems with other boundary conditions and rectangular domains in higher dimensions.

Let the spatial interpolation $\mathcal{I}_h:C(X)\rightarrow C(\overline{\Omega})$ or  $C(\overline{\Omega})\rightarrow C(\overline{\Omega})$ be the piecewise linear interpolation operator with respect to the nodes associated with the mesh. More precisely, for any function $v(x)\in C(X)$ or $C(\overline{\Omega})$ , the interpolation $\mathcal{I}_hv\in C(\overline{\Omega})$ is piecewise linear  and
\begin{align*}
\mathcal{I}_hv(x)=\sum_{i=1}^Nv(x_i)\psi_i(x),\quad \forall\,x\in\Omega,
\end{align*}
where $\psi_i(x)$  is the piecewise linear basis function satisfying $\psi_i(x_i)=1$ and $\psi_i(x_j)=0$ when $i\neq j$. Note that under the  homogeneous Dirichlet  boundary condition, $C(X^*)$ can be regarded 
as a linear subspace of $C(X)$  in the sense of isometric isomorphism through the zero  extension to the boundary nodes.
We have the following error estimates.

\begin{theorem}\label{thm-etd1-conv}
	For the fixed terminal time $T>0$, assume that  $\mathbf{v}\in C^1([0,T],C^1(\Omega))$, the exact solution $u$ to the model problem \eqref{stab} belongs to $C^1([0,T],C^3(\overline\Omega)\cap C_0(\overline\Omega))$  ($C_0(\overline{\Omega})=\{\phi(x)\in C(\overline{\Omega}),\,\phi|_{\partial\Omega}=0\}$) and $\{U^n\in C(X)\}_{n\geq0}$ is generated by the fully discrete ETD1 scheme \eqref{etd1-fully} with $U^0=u_0(X)$.  Then for any $\tau>0$, $h>0$, it holds that
	\begin{align*}
	\|u(\cdot,t^n)-\mathcal{I}_hU^n(\cdot)\|\leq C(\tau+h),\quad\forall\, t^n\leq T,
	\end{align*}
	where the constant $C>0$ is independent of $\tau$ and $h$.
\end{theorem}
\begin{proof}
	First of all, by Taylor expansion and triangle inequality, recalling the homogeneous boundary conditions, there holds
	\begin{align}
	\|u(\cdot,t^n)-\mathcal{I}_hU^n(\cdot)\|\leq& \|\mathcal{I}_hu(\cdot,t^n)-\mathcal{I}_hU^n(\cdot)\|+\|u(\cdot,t^n)-\mathcal{I}_hu(\cdot,t^n)\|\nonumber\\
	\leq &\|u(x,t^n)-U^n(x)\|_{X^*}+h\|\partial_xu(\cdot,t^n)\|.\label{eq:errordc}
	\end{align}
	Under the assumptions of Theorem \ref{thm-etd1-conv}, the MBP implies that
	$\|u(t)\|\leq\beta$ and $\|\mathcal{I}_hU\|\leq\beta$. 
	Hence, it suffices to consider the following ``error function" $e^n(x)\in C(X^*)$ with the homogeneous Dirichlet boundary condition  as
	\begin{equation}\label{eq:err}
	e^n(x)=U^n(x)-u(x,t^n),\quad x\in X^*.
	\end{equation}
	From $t^n$ to $t^{n+1}$, the ETD1 solution \eqref{etd1-fully} is given by $U^{n+1}(x)=W_1^n(x,\tau)$ for $x\in X^*$ with the function $W_1^n:=W_1^n(s)=W_1^n(x,s)$ solving 
	\begin{align}\label{etd1-disc}
	\left\{\begin{array}{ll}
	\partial_s {W_1^n}=\mathcal{L}_{\mathbf{v}^n,h}^\kappa [U^n]W_1^n+N\left(U^{n}\right), &\quad x\in X^*, s \in(0, \tau], \\[4pt]
	W_1^n(x,0)=U^{n}(x), & \quad x\in X^*,
	\end{array}\right.
	\end{align}
	where $\mathbf{v}^n=\mathbf{v}(t^n)$.
	Based on the spatial discretization \eqref{eq:spatial}, $u(x,t^n+s)$ satisfies the following equation
	for $x\in X^*$ and $s \in(0, \tau]$,
	\begin{equation}\label{eq:ceq}
	{\partial_s u(x,t^n+s)}=\mathcal{L}_{\mathbf{v}^n,h}^\kappa [U^n(x)]u(x,t^n+s)+N\left(u(x,t^n)\right)+R^n(x,s),
	\end{equation}
	where the  local truncation error function $R^n(s):=R^n(x,s)\in C(X^*)$  for any $s\in(0,\tau]$ is given by
	\begin{align}\label{eq:local}
	R^n(s)=&\mathcal{L}^\kappa_{\mathbf{v}(t_n+s)}[u(t^n+s)]u(t^n+s)-\mathcal{L}^\kappa_{\mathbf{v}^n,h}[u(t^n)]u(t^n+s)
	+N\left(u(t^n+s)\right)-N\left(u(t^{n})\right)\nonumber\\
	&+\mathcal{L}^\kappa_{\mathbf{v}^n,h}[u(t^n)]u(t^n+s)-\mathcal{L}^\kappa_{\mathbf{v}^n,h}[U^n]u(t^n+s).
	\end{align}
	Under the assumptions of Theorem \ref{thm-etd1-conv}, using Taylor expansion and the MBP properties, it is easy to check that
	\begin{align*}
	&\left\|\mathcal{L}^\kappa_{\mathbf{v}(t^n+s)}[u(t^n+s)]u(t^n+s)-\mathcal{L}^\kappa_{\mathbf{v}^n,h}[u(t^n)]u(t^n+s)\right\|_{X^*}\leq C_1(\tau+h),\\
	&\left\|\mathcal{L}^\kappa_{\mathbf{v}^n,h}[u(t^n)]u(t^n+s)-\mathcal{L}^\kappa_{\mathbf{v}^n,h}[U^n]u(t^n+s)\right\|_{X^*}\leq C_2\|u(t^n)-U^n\|_{X^*},\\
	&\left\|N\left(u(t^n+s)\right)-N\left(u(t^{n})\right)\right\|_{X^*}\leq C_3\tau,
	\end{align*}
	and 
	\begin{equation}\label{eq:local-bd}
	\|R^n(s)\|_{X^*}\leq C_2\|u(t^n)-U^n\|_{X^*}+(C_1+C_3)\tau+C_1h.
	\end{equation}
	Define the local error function $e_1^n(s):=e_1^n(x,s)=W_1^n(x,s)-u(x,t^n+s)\in C(X^*)$ for any $s\in[0,\tau]$ and subtracting \eqref{eq:ceq} from \eqref{etd1-disc}, we have
	\begin{equation}\label{eq:int1}\begin{cases}
	{\partial_s e_1^n(s)}=\mathcal{L}_{\mathbf{v},h}^\kappa[U^n]e_1^n(s)+N(U^n)-N(u(t^n))-R_1^n(s), & x\in X^*, s \in(0, \tau],\\
	e_1^n(x,0)=e^n(x), & x\in X^*.
	\end{cases}
	\end{equation}
	The MBP property and Lemma \ref{lemma-stab-nonlinear} imply
	\begin{equation}\label{eq:Nd}
	\|N(U^n)-N(u(t^n))\|_{X^*}\leq 2\kappa \|U^n-u(t^n)\|_{X^*}=2\kappa\|e^n\|_{X^*}.
	\end{equation}
	By Duhamel's principle, we obtain from \eqref{eq:int1} that
	\begin{align}\label{eq:intg}
	e_1^n(\tau)=e^{\tau \mathcal{L}_{\mathbf{v},h}^\kappa[U^n]}e^n+\int_0^\tau e^{(\tau-s)  \mathcal{L}_{\mathbf{v},h}^\kappa[U^n] }\left(N(U^n)-N(u(t^n))-R^n(s)\right)\,ds.
	\end{align}
	Since $\mathcal{L}_{\mathbf{v},h}[U^n]$ satisfies Assumption  \ref{assump-linear}, by Lemma \ref{lem-linear} the corresponding contraction semigroup enjoys the property of 
	$$\|e^{\tau \mathcal{L}_{\mathbf{v},h}^\kappa[U^n]}e^n\|_{X^*}\leq e^{-\tau\kappa}\|e^n\|_{X^*},$$ combining \eqref{eq:local-bd} and \eqref{eq:Nd}, we have from \eqref{eq:intg} that ($e^{n+1}=e^n_1(\tau)$)
	\begin{align}\label{etd1-ind}
	\left\|e^{n+1}\right\|_{X^*} \leq & \mathrm{e}^{-\kappa \tau}\left\|e^{n}\right\|_{X^*}+\int_{0}^{\tau}e^{-\kappa(\tau-s)}\left(\|R^n_1(s)\|_{X^*}+\|N\left(u^n\right)-N(u(t^n))\|_{X^*}\right)\,ds \notag\\
	\leq & \mathrm{e}^{-\kappa \tau}\left\|e^{n}\right\|_{X^*}+\left(2\kappa\|e^n\|_{X^*}+C_2\|e^n\|_{X^*}+C_4(\tau+h)\right)\int_{0}^{\tau} \mathrm{e}^{-\kappa(\tau-s)} \,ds \notag\\
	=& \mathrm{e}^{-\kappa \tau}\left\|e^{n}\right\|_{X^*}+\frac{1-\mathrm{e}^{-\kappa \tau}}{\kappa}\left(2\kappa\|e^n\|_{X^*}+C_2\|e^n\|_{X^*}+C_4(\tau+h)\right)\notag\\
	\leq &(1+C_5\tau)\left\|e^{n}\right\|_{X^*}+C_6 \tau(\tau+h),
	\end{align}
	where we have used the fact that $1-s\leq e^{-s} \leq 1+s$ for  $s>0$.  Therefore, recalling $e^0(x)=0$ for any $x\in X^*$), we obtain
	\[
	\begin{aligned}
	\left\|e^{n+1}\right\|_{X^*} & \leq(1+C_5 \tau)^{n+1}\left\|e^{0}\right\|_{X^*}+C_6 \tau(\tau+h) \sum_{k=0}^{n}(1+C_5\tau)^{k} \\
	&=\frac{C_6}{C_5}(\tau+h)\left[(1+C_5\tau)^{n+1}-1\right] \leq C_7e^{C_5 (n+1) \tau} (\tau+h),
	\end{aligned}
	\]
	and the desired estimates at $t^{n+1}$ hold in view of \eqref{eq:errordc}.  It is easy to check all the constants $C_j$ appearing in the proof are independent of $\tau$ and $h$. \hfill
\end{proof}

For the second order  ETDRK2 scheme \eqref{etd2rk-fully:1}-\eqref{etd2rk-fully:2}, the error estimates can be established as follows.
\begin{theorem}\label{thm-etd2rk-conv}
	For the fixed terminal time $T>0$, assume that $\mathbf{v}\in C^1([0,T],C^1(\Omega))$, the exact solution $u$ to the model problem \eqref{stab} belongs to $C^2([0,T],C^4(\overline\Omega)\cap C_0(\overline\Omega))$ and $\{U^n\in C(X)\}_{n\geq0}$ is generated by the fully discrete ETDRK2 scheme \eqref{etd2rk-fully:1}-\eqref{etd2rk-fully:2} with $U^0=u_0(X)$. Then for any $\tau,h\in(0,1]$, it holds that
	\begin{align*}
	\|u(\cdot,t^n)-\mathcal{I}_hU^n(\cdot)\|\leq C(\tau^2+h),\quad\forall\, t^n\leq T,
	\end{align*}
	where the constant $C>0$ is independent of $\tau$ and $h$.
\end{theorem}
\begin{proof} The proof is similar to that of the ETD1 case and we shall only sketch the proof below.
	First of all, the ETDRK2 solution is given by $U^{n+1}(x)=W_2^n(x,\tau)\in C(X^*)$ with the function $W_2^n(s)$ solving
	\begin{equation}\label{etd2-disc}
	\begin{cases}
	\partial_s W_2^n=\frac12\left(\mathcal{L}_{\mathbf{v}^n,h}^\kappa[U^n]+\mathcal{L}_{\mathbf{v}^{n+1},h}^\kappa[\hat{U}^{n+1}]\right)W_2^n\\
	\qquad\qquad\qquad\qquad+(1-\frac{s}{\tau})N(U^n)+\frac{s}{\tau}N(\hat{U}^{n+1}),&  x\in X^*, s \in(0, \tau],\\
	W_2^n(x,0)=U^{n}(x),& x\in X^*,
	\end{cases}
	\end{equation}
	where $\hat{U}^{n+1}\in C(X^*)$ is given by $\hat{U}^{n+1}=W_1^n(\tau)$ with $W_1^n(s)$ satisfying \eqref{etd1-disc}.  Based on the proof of Theorem \ref{thm-etd1-conv}, we know
	\begin{equation}\label{eq:est1}
	\left\|\hat{U}^{n+1}-u(t^{n+1})\right\|_{X^*}\leq (1+C_1\tau)\|e^n\|_{X^*}+C_2\tau(\tau+h).
	\end{equation}
	We define the local truncation error function $R^n(s):=R^n(x,s)\in C(X^*)$ ($x\in X^*$) for any $s\in(0,\tau]$  as
	\begin{align}\label{eq:locale2}
	&R^n(s)=\frac{\partial u(t^n+s)}{\partial s}-\frac12\left(\mathcal{L}_{\mathbf{v}^n,h}^\kappa[U^n]+\mathcal{L}_{\mathbf{v}^{n+1},h}^\kappa[\hat{U}^{n+1}]\right)u(t^n+s)\nonumber\\
	&\qquad\qquad-(1-\frac{s}{\tau})N(u(t^n))-\frac{s}{\tau}N(u(t^{n+1})),
	\end{align}
	which can be decomposed as follows in view of  \eqref{gac}, \eqref{eq:est1} and Taylor expansion (as in the ETD1 case and details are omitted):
	\begin{align}\label{eq:localdec}
	R^n(s)=R_1^n(s)+R_2^n(s)+R_3^n(s),
	\end{align} 
	with 
	\begin{align}
	&\|R_1^n(s)\|_{X^*}\leq C_4\left(\|e^n\|_{X^*}+\left\|\hat{U}^{n+1}-u(t^{n+1})\right\|_{X^*}\right)
	\leq C_5\left(\|e^n\|_{X^*}+\tau(\tau+h)\right),\label{eq:local:1}\\
	&R_2^n(s)=\left(\frac{s}{\tau}-\frac12\right)\left(\mathcal{L}^\kappa_{\mathbf{v}^{n+1}}[u(t^{n+1})]-\mathcal{L}^\kappa_{\mathbf{v}^n}[u(t^n)]\right)u(t^n),\label{eq:local:2}\\
	&
	\|R_3^n(s)\|_{X^*}\leq  C_6(\tau^2+h).\label{eq:local:3}
	\end{align}
	Define the local error function in the interval $[t^n,t^{n+1}]$ as $e_{2}^n(s)=W_2^n(x,s)-u\left(x,t^n+s\right)\in C(X^*)$. By subtracting \eqref{eq:locale2} from \eqref{etd2-disc}, it yields 
	\begin{equation}\label{eq:erreq2}\begin{cases}
	\partial_s e_2^n(x,s)=\frac12\left(\mathcal{L}_{\mathbf{v}^n,h}^\kappa[U^n]+\mathcal{L}_{\mathbf{v}^{n+1},h}^\kappa[\hat{U}^{n+1}]\right)e_2^n(x,s)+(1-\frac{s}{\tau})N(U^n)+\frac{s}{\tau}N(\hat{U}^{n+1})\\
	\qquad\qquad\quad-(1-\frac{s}{\tau})N(u(t^n))-\frac{s}{\tau}N(u(t^{n+1}))-R^n(x,s),  \quad s \in(0, \tau], \,x\in X^*,\\
	e_2^{n}(0,x)=e^{n}(x), \quad x\in X^*.
	\end{cases}
	\end{equation}
	Denoting
	\begin{equation*}
	\begin{split}
f^n(s):=f^n(x,s)=&(1-\frac{s}{\tau})N(U^n)+\frac{s}{\tau}N(\hat{U}^{n+1})
-(1-\frac{s}{\tau})N(u(t^n))-\frac{s}{\tau}N(u(t^{n+1}),
	\end{split}
	\end{equation*}
	by using the MBP, Lemma \ref{lemma-stab-nonlinear} and \eqref{eq:est1}, we have
	\begin{align}\label{eq:fns}
	\left\|f^n(s)\right\|_{X^*}
	\leq&2\kappa(1-\frac{s}{\tau})\left\|U^n-u(t^n)\right\|_{X^*}+2\kappa\frac{s}{\tau}\left\|\hat{U}^{n+1}-u(t^{n+1})\right\|_{X^*}\nonumber\\
	\leq& C_7\left(\|e^n\|_{X^*}+\tau(\tau+h)\right).
	\end{align}
	Applying Duhamel's principle to \eqref{eq:erreq2},  we get
	\begin{align}
	e^n_2(\tau)=&e^{\frac{\tau}{2}\left(\mathcal{L}_{\mathbf{v}^n,h}^\kappa[U^n]+\mathcal{L}_{\mathbf{v}^{n+1},h}^\kappa[\hat{U}^{n+1}]\right)}e^n\nonumber\\
	&
	+\int_0^\tau e^{\frac{\tau-s}{2}\left(\mathcal{L}_{\mathbf{v}^n,h}^\kappa[U^n]+\mathcal{L}_{\mathbf{v}^{n+1},h}^\kappa[\hat{U}^{n+1}]\right)}\left(f^n(s)-R_1^n(s)-R_2^n(s)-R_3^n(s)\right)\,ds.\label{eq:errorint2}
	\end{align}
	Since $\int_0^\tau R_2^n(s)\,ds=0$, by applying Taylor expansion to the matrix exponential, we have 
	\begin{align}
	&\left\|\int_0^\tau e^{\frac{\tau-s}{2}\left(\mathcal{L}_{\mathbf{v}^n,h}^\kappa[U^n]+\mathcal{L}_{\mathbf{v}^{n+1},h}^\kappa[\hat{U}^{n+1}]\right)}R_2^n(s)\,ds\right\|_{X^*}\nonumber\\
	&\qquad\leq C_8\tau^2\left\|\left(\mathcal{L}_{\mathbf{v}^{n+1},h}^\kappa[u(t^{n+1})]-\mathcal{L}_{\mathbf{v}^n,h}^\kappa[u(t^n)]\right)u(t^n)\right\|_{X^*}\leq C_9\tau^3,\label{eq:intest}
	\end{align}
	where $C_9$ depends on the $C^2([0,T],C^4(\overline\Omega))$ norm of $u(x,t)$.
	Combining \eqref{eq:local:1}-\eqref{eq:local:3}, \eqref{eq:fns} and \eqref{eq:intest} and noticing that $e^{n+1}=e_2^n(\tau)$  and the semigroup $e^{\frac{\tau}{2}\left(\mathcal{L}_{\mathbf{v}^n,h}^\kappa[U^n]+\mathcal{L}_{\mathbf{v}^{n+1},h}^\kappa[\hat{U}^{n+1}]\right)}$ has an upper bound $e^{-\kappa \tau}$ in the $C(X^*)$ norm,  we can get
	\begin{align}
	\|e^{n+1}\|_{X^*}\leq &e^{-\tau\kappa}\|e^n\|_{X^*}+C_{10}(\|e^n\|_{X^*}+\tau^2+h)\int_0^\tau e^{-\kappa(\tau-s)}\,ds\nonumber+C_9\tau^3\\
	\leq&(1+C_{11}\tau)\|e^n\|_{X^*}+C_{12}\tau(\tau^2+h).
	\end{align}
	Following the proof for the ETD1 case, we then conclude $\|e^n\|_{X^*}\leq C_{13}(\tau^2+h)$ for $0 \leq t^n\leq T$ and the error estimates in Theorem \ref{thm-etd2rk-conv} hold based on \eqref{eq:errordc}.  \hfill
\end{proof}

\begin{remark}
We note that the regularity assumptions in Theorems \ref{thm-etd1-conv} and \ref{thm-etd2rk-conv} are not sharp and can be weakened. The derived error estimates also hold for the cases of periodic and homogeneous Neumann boundary conditions. 
\end{remark}

\section{Numerical experiments}\label{sec:simul}
In this section, we present various 2D and 3D numerical experiments to demonstrate the accuracy and  the discrete MBP  of the proposed fully discrete ETD1 \eqref{etd1-fully} and ETDRK2 \eqref{etd2rk-fully:1}-\eqref{etd2rk-fully:2} schemes. ETDRK2 scheme is used for all examples, while the ETD1 scheme is only considered in the temporal convergence test due to the lack of  accuracy. 

\subsection{Convergence tests}
Let us take the domain $\Omega=(-0.5,0.5)^2$ and the terminal time $T=0.1$. We consider the 2D convective Allen-Cahn equation \eqref{gac} with the mobility function $M(u)\equiv1$ and $\epsilon=0.01$. In addition, the nonlinear reaction  $f=-F^\prime$ is given by the double-well potential case \eqref{dw}, thus the stabilizing coefficient is taken as $\kappa=2$ so that \eqref{coef} is satisfied \cite{DJLQ21}. The initial value is $u_0(x,y)=\cos(2\pi x)\cos(2\pi y)$ 
with the velocity field $\mathbf{v}=[1,1]^T$. The periodic boundary condition is imposed.

By fixing the space mesh size $h=1/1024$, we first test the convergence in time with various time step sizes. To compute the solution errors, we set the solution obtained by the ETDRK2 scheme with $\tau=1/1024$ as the referential value. The $L^\infty$  and $L^2$  errors of the numerical solutions at the terminal time $T=0.1$ with different time step sizes and their corresponding convergence rates in time are reported in Table \ref{tab-temporal-conv}, where the expected temporal convergence rates (order 1 for ETD1 and order 2 for ETDRK2) are clearly observed.
\begin{table}[!ht]
	\normalsize\small
	\begin{center}
		\renewcommand{\arraystretch}{1.1}
		\small\begin{tabular}{|c|cccc|cccc|}
			\hline
			\multirow{2}{2em}{$\tau$}& \multicolumn{4}{|c|}{ETD1} &\multicolumn{4}{c|}{ETDRK2}\\
			\cline{2-5}\cline{6-9}
			& $L^\infty$ Error & Rate & $L^2$ Error & Rate& $L^\infty$ Error & Rate & $L^2$ Error & Rate\\
			\hline
			$1/16$&6.3251e-3&-&3.6021e-3&-&9.9082e-5&-&5.9261e-5&-\\
			$1/32$&3.1267e-3&1.01&1.7792e-3&1.01&2.4905e-5&1.99&1.4890e-5&1.99\\
			$1/64$&1.5162e-3&1.04&8.6245e-4&1.04&6.2266e-6&1.99&3.7220e-6&2.00\\
			$1/128$&7.0833e-4&1.09&4.0283e-4&1.09&1.5406e-6&2.01&9.2082e-7&2.01\\
$1/256$&3.0373e-4&1.22&1.7272e-4 &1.22&3.6707e-7&2.06&2.1939e-7&2.06\\
			\hline
		\end{tabular}
	\caption{\normalsize  Results on $L^\infty$ and $L^2$ errors of the numerical solutions  at $T=0.1$ and their corresponding convergence rates in time for the fully discrete ETD schemes.}
		\label{tab-temporal-conv}
	\end{center}
\end{table}

Next, we test the convergence with respect to the spatial mesh size $h$ by fixing the temporal step size $\tau=1/2048$. The numerical solution of the convective Allen-Cahn equation obtained by the ETDRK2 scheme with $h=1/512$ is treated as the referential value.
 The $L^\infty$  and $L^2$  errors of the numerical solutions at the terminal time $T=0.1$ along the spatial mesh refinement and their corresponding convergence rates are presented in Table \ref{tab-spatial-conv}. It is observed that the spatial convergence rates gradually converge to 1, which is  consistent with the upwind finite difference approximation as expected. 

\begin{table}[!ht]
	\normalsize\small
	\begin{center}
		\renewcommand{\arraystretch}{1.1}
		\begin{tabular}{|c|cccc|}
			\hline
			$1/h$ & $L^\infty$ Error & Rate & $L^2$ Error & Rate\\
			\hline\hline
			8&1.4308e-1&-&7.6016e-2&-\\
			16&9.0813e-2&0.65&4.9845e-2&0.60\\
			32&4.8558e-2&0.90&2.7861e-2&0.83\\
			64&2.4858e-2&0.96&1.4636e-2&0.92\\
			128&1.2536e-2&0.98&7.4887e-3&0.96\\
			256&6.2909e-3& 0.99&3.7862e-3&0.98\\
			\hline
		\end{tabular}
	\caption{\normalsize  Results on $L^\infty$ and $L^2$ errors of the numerical solutions  at $T=0.1$ and their corresponding convergence rates  
		in space for the fully discrete ETDRK2 scheme.}
		\label{tab-spatial-conv}
	\end{center}
\end{table}

\subsection{MBP tests}
Next, we numerically simulate the 2D convective Allen-Cahn equation \eqref{gac} with $\epsilon=0.01$ to investigate the preservation of the discrete MBP in the long-time phase separation process. Let us take the computational domain $\Omega=(-0.5,0.5)^2$ and the terminal time $T=50$. 
We set the initial configuration 
$u_0(x,y)=0.9\sin(100\pi x)\sin(100\pi y)$, the velocity field $\mathbf{v}(x,y,t)=e^{-t}[\sin(2\pi x),-\cos(2\pi y)]^T$ and $M(u)=1-u^2$. 
The homogeneous Neumann boundary condition is imposed. The mesh size is set to be $h=1/64$. The ETDRK2 scheme is taken for all tests in this subsection.

\begin{figure}[!ht]
	\centering
	{\includegraphics[width=0.35\textwidth]{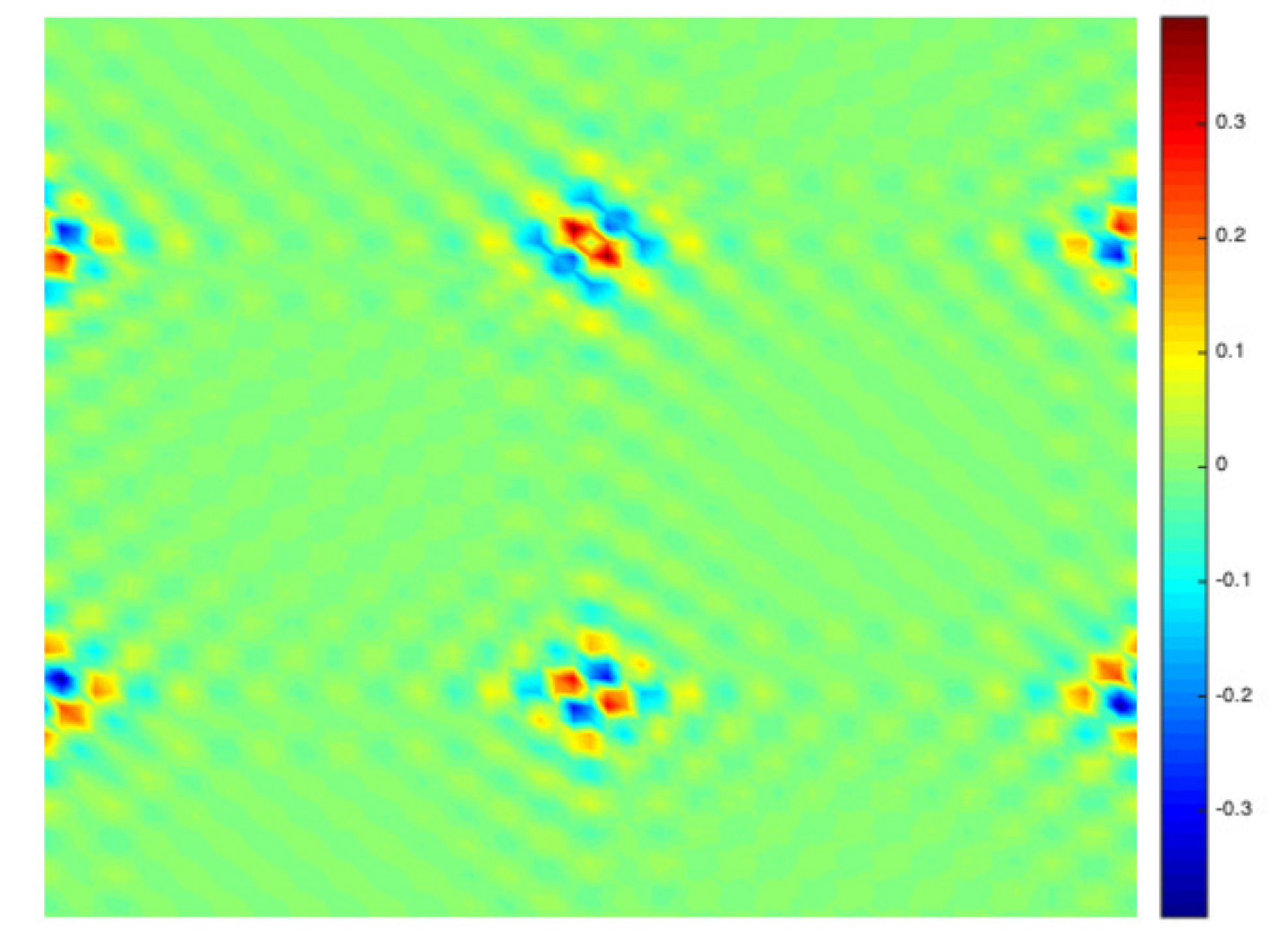}}
	{\includegraphics[width=0.35\textwidth]{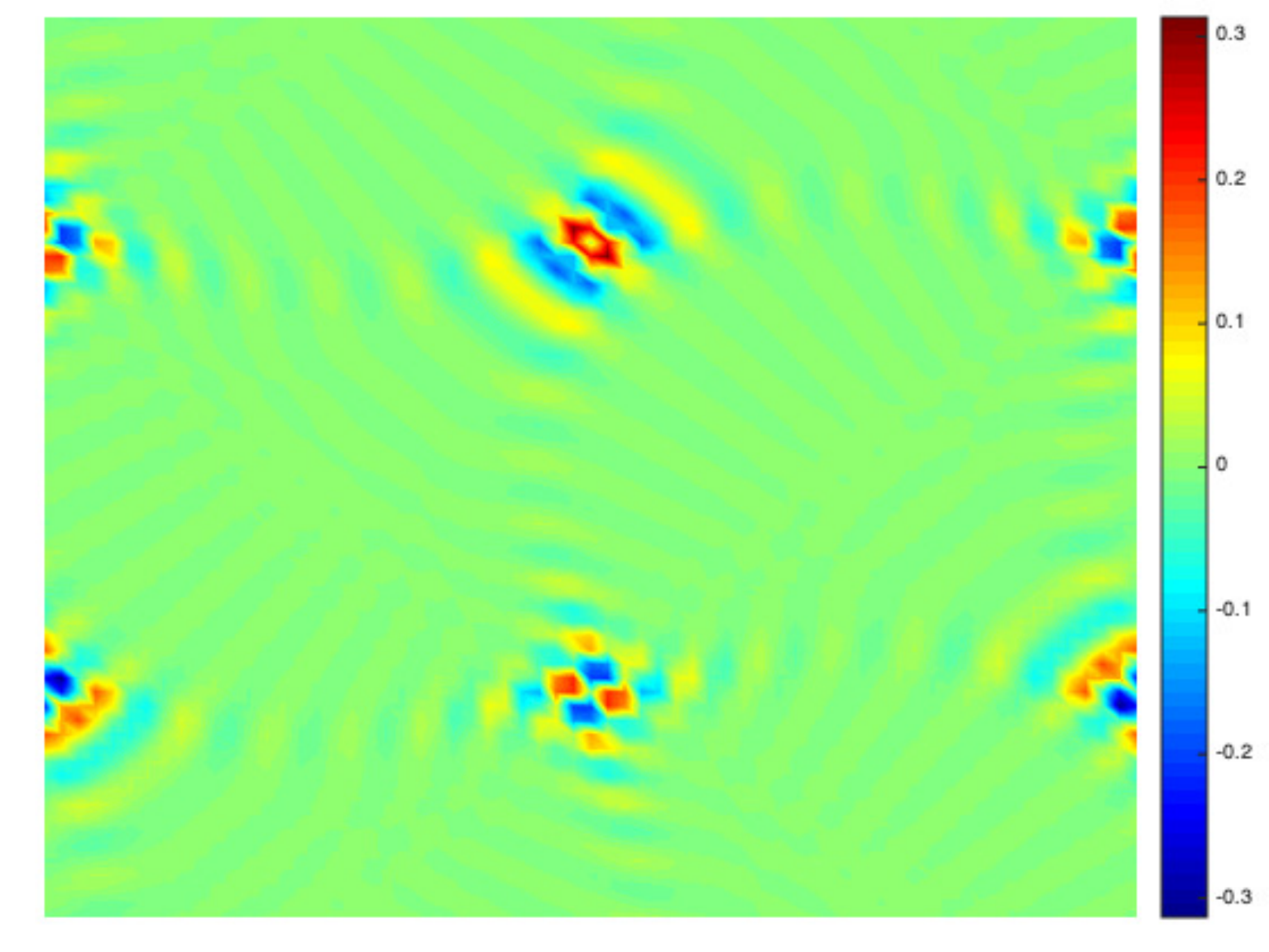}}
	{\includegraphics[width=0.35\textwidth]{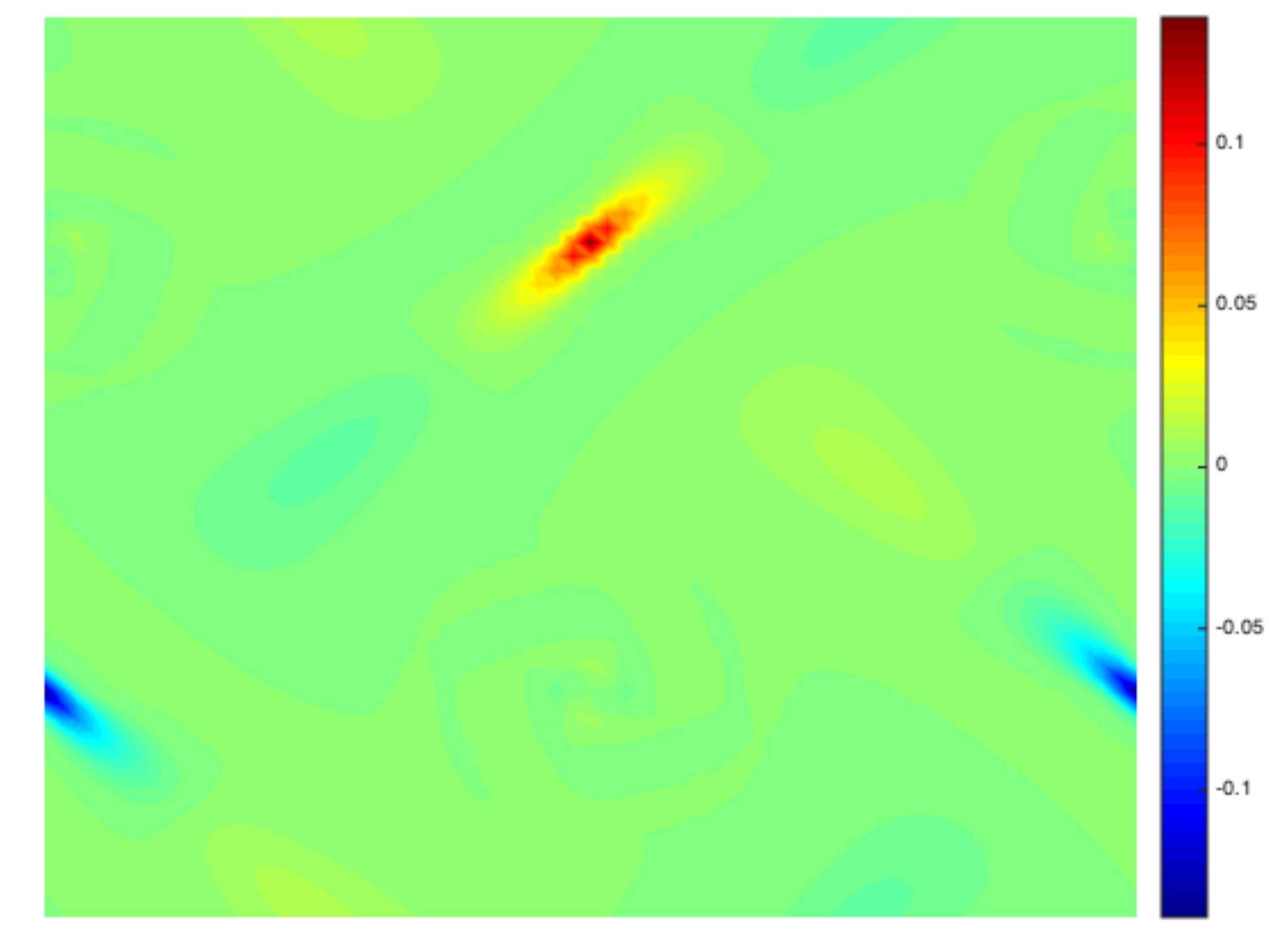}}
	{\includegraphics[width=0.35\textwidth]{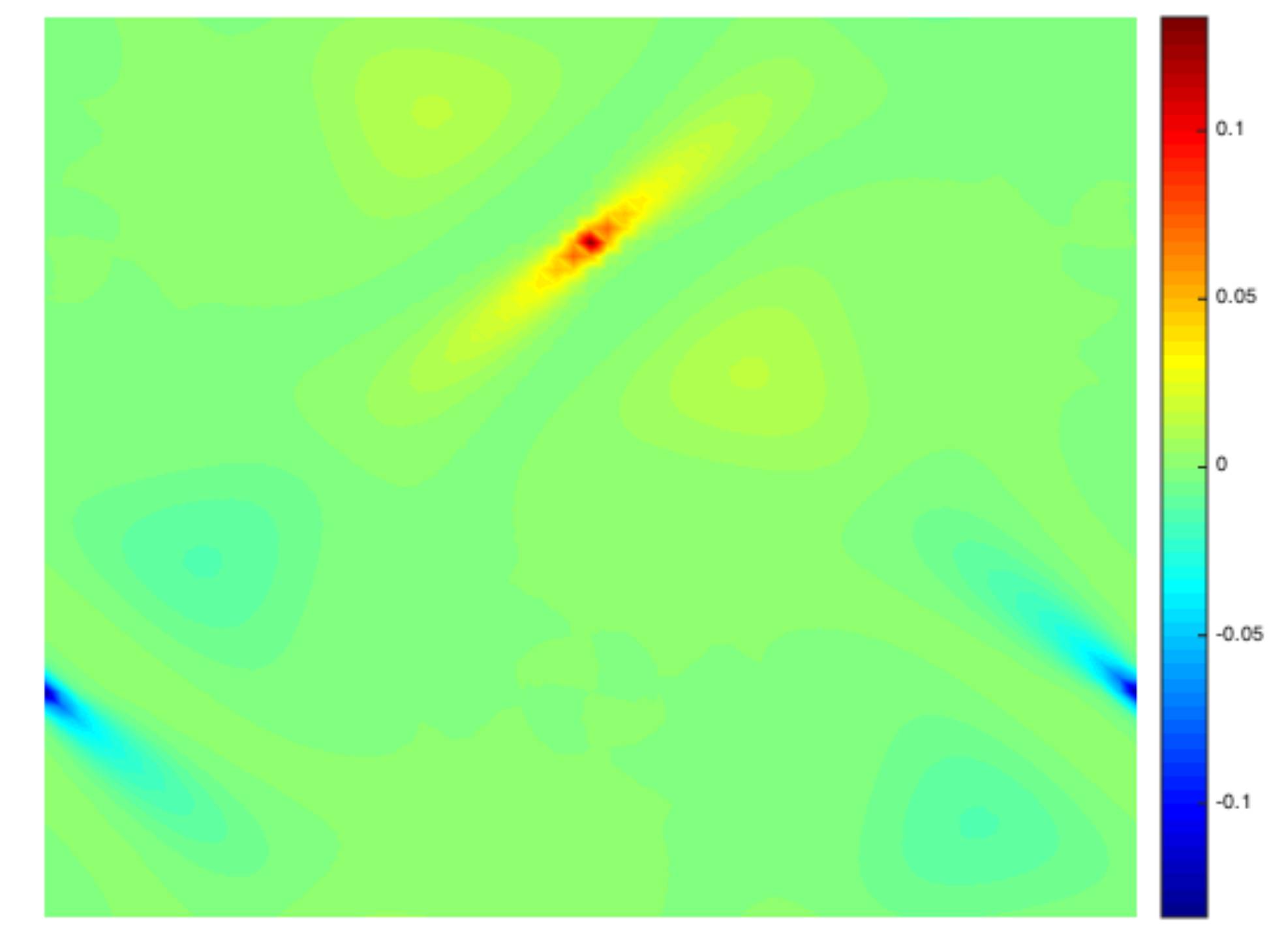}}
	{\includegraphics[width=0.35\textwidth]{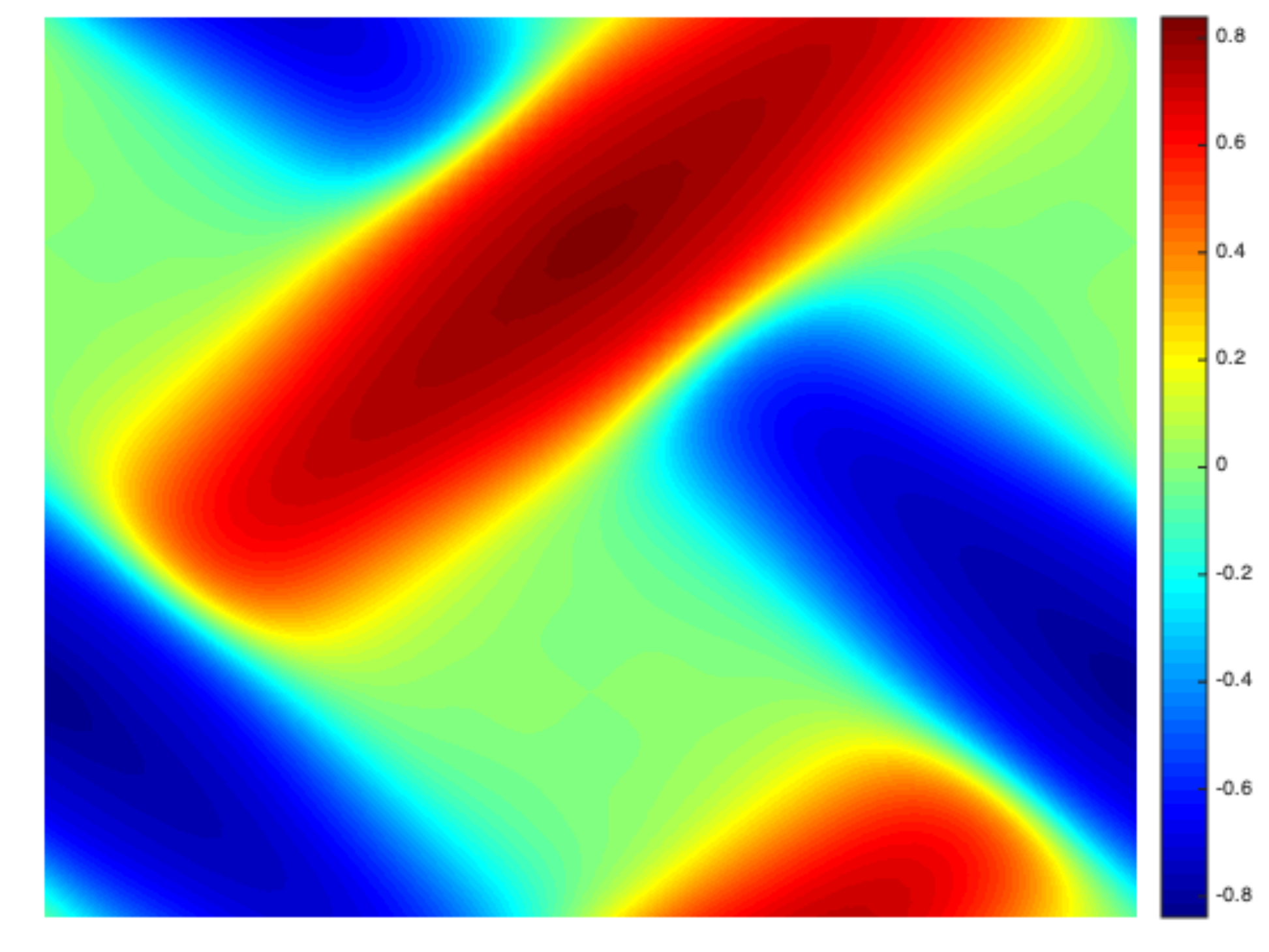}}
	{\includegraphics[width=0.35\textwidth]{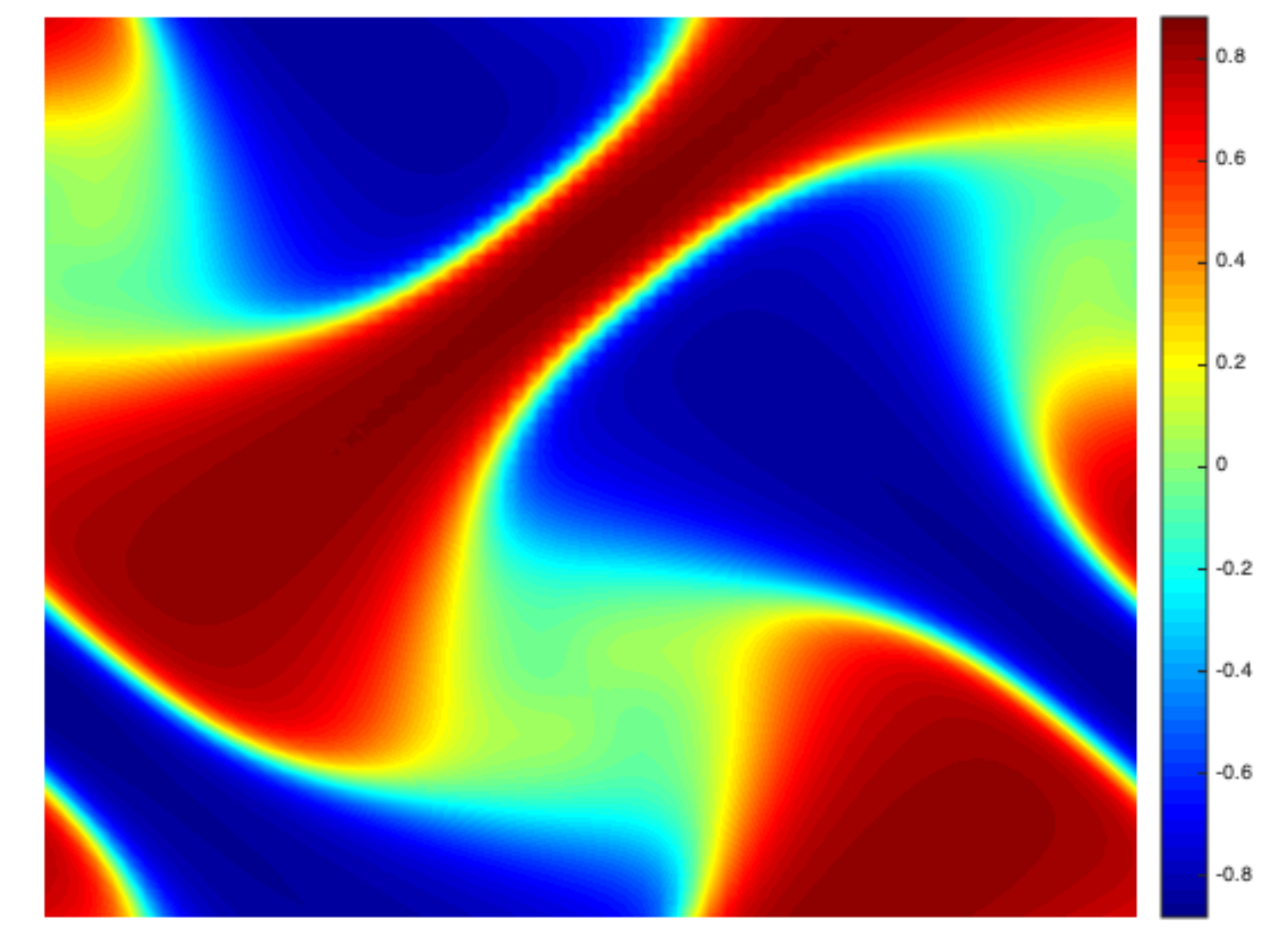}}
	{\includegraphics[width=0.35\textwidth]{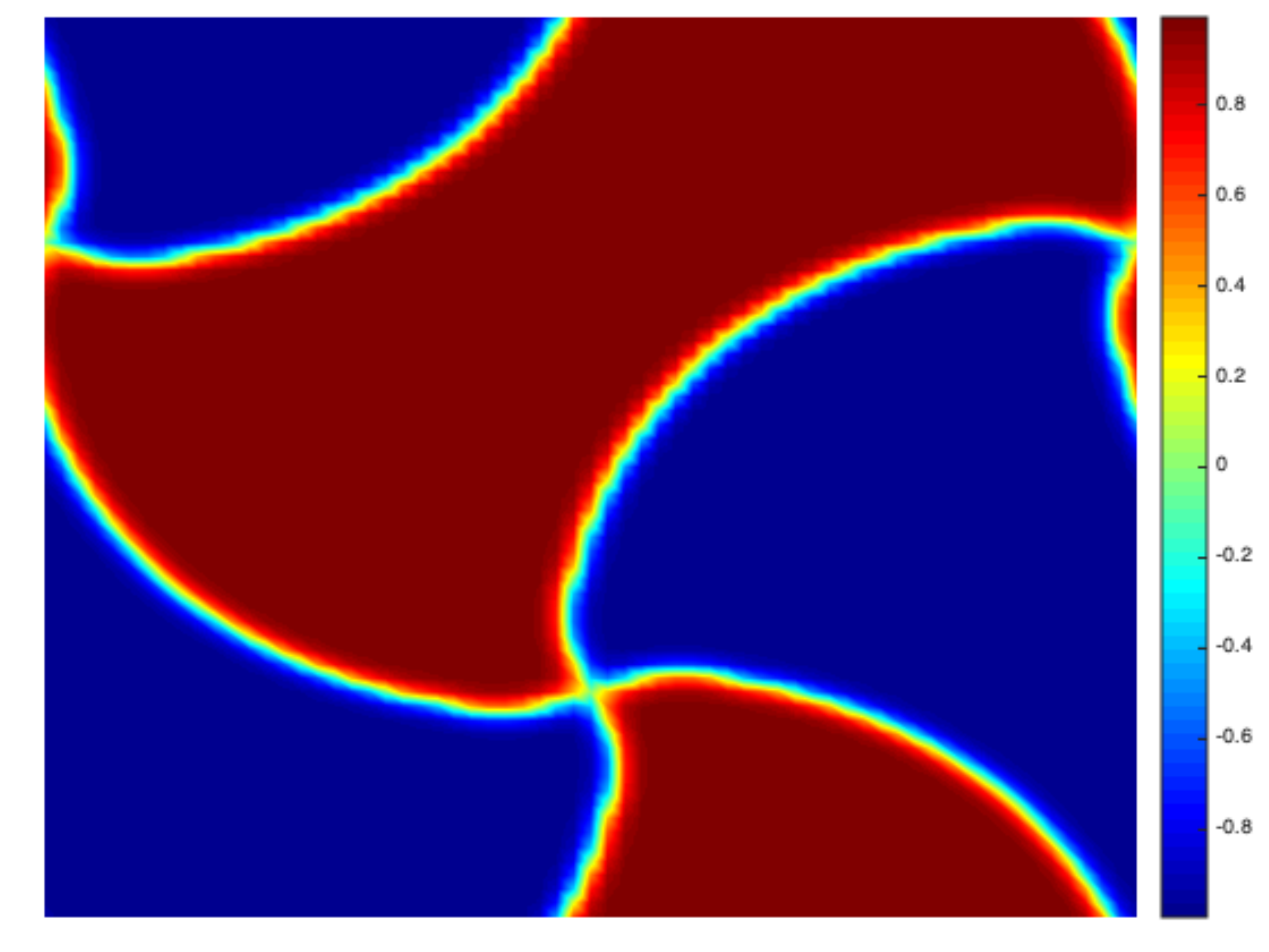}}
	{\includegraphics[width=0.35\textwidth]{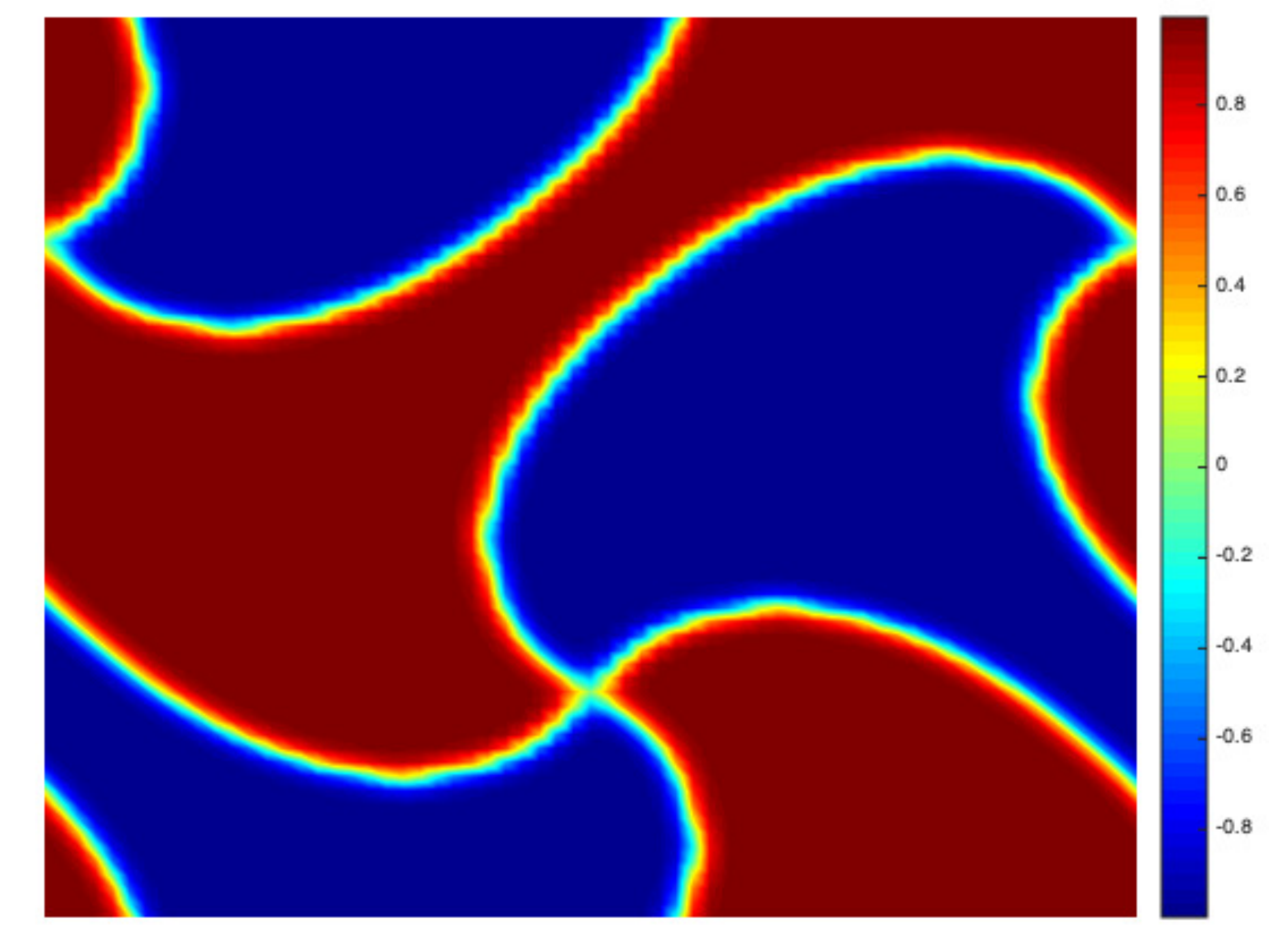}}
	\caption{The numerical solutions  at $t$ = 0.1, 1, 8, 50, respectively (top to bottom) for the 2D convective Allen-Cahn equation with the double-well potential. Left column: $\tau$=0.1, right column: $\tau$=0.01.}
	\label{fig-dw2D-phase}
\end{figure}

\begin{figure}[!ht]
	\centering
	{\includegraphics[width=0.45\textwidth]{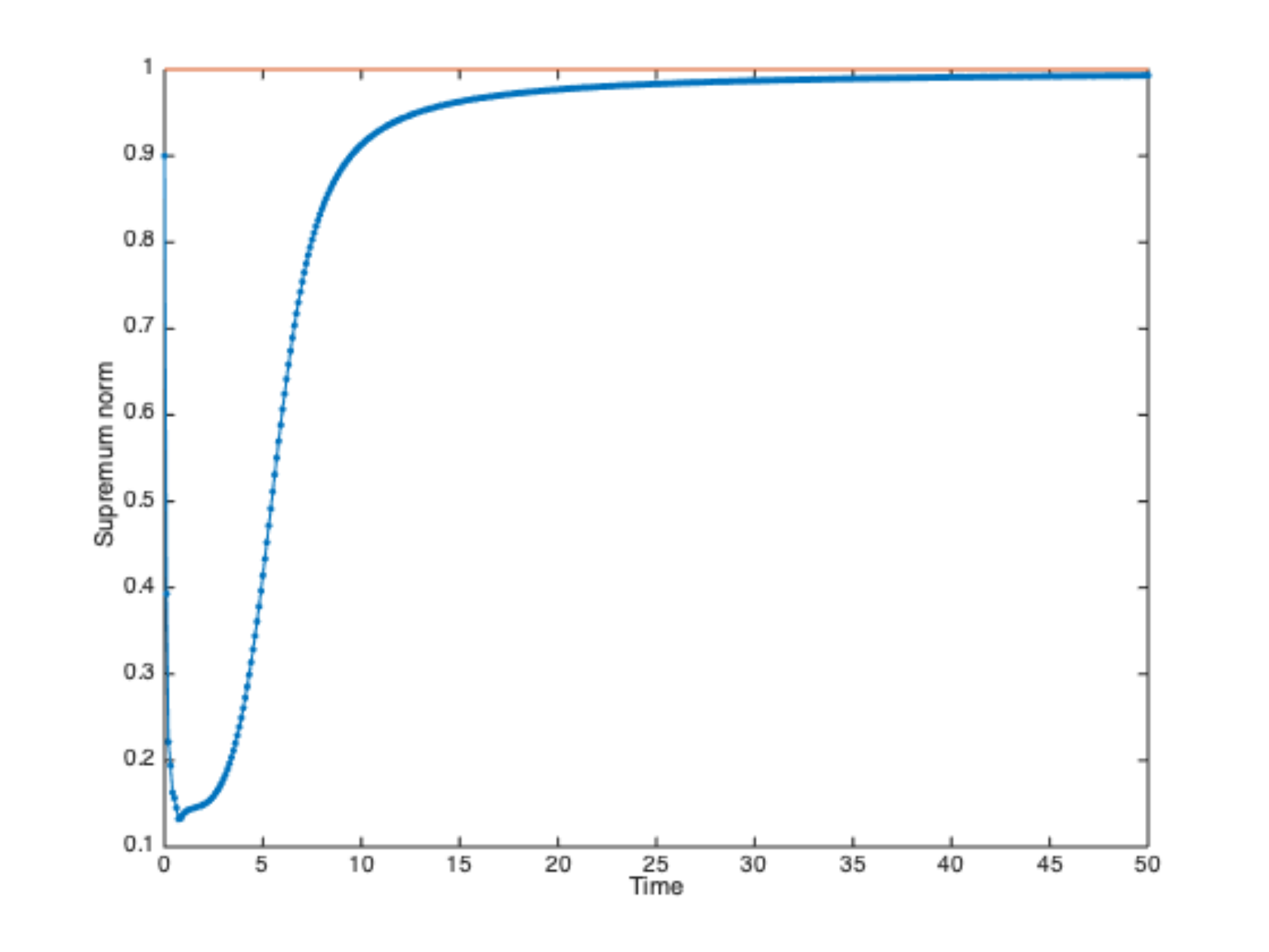}}\hspace{-0.3cm}
	{\includegraphics[width=0.45\textwidth]{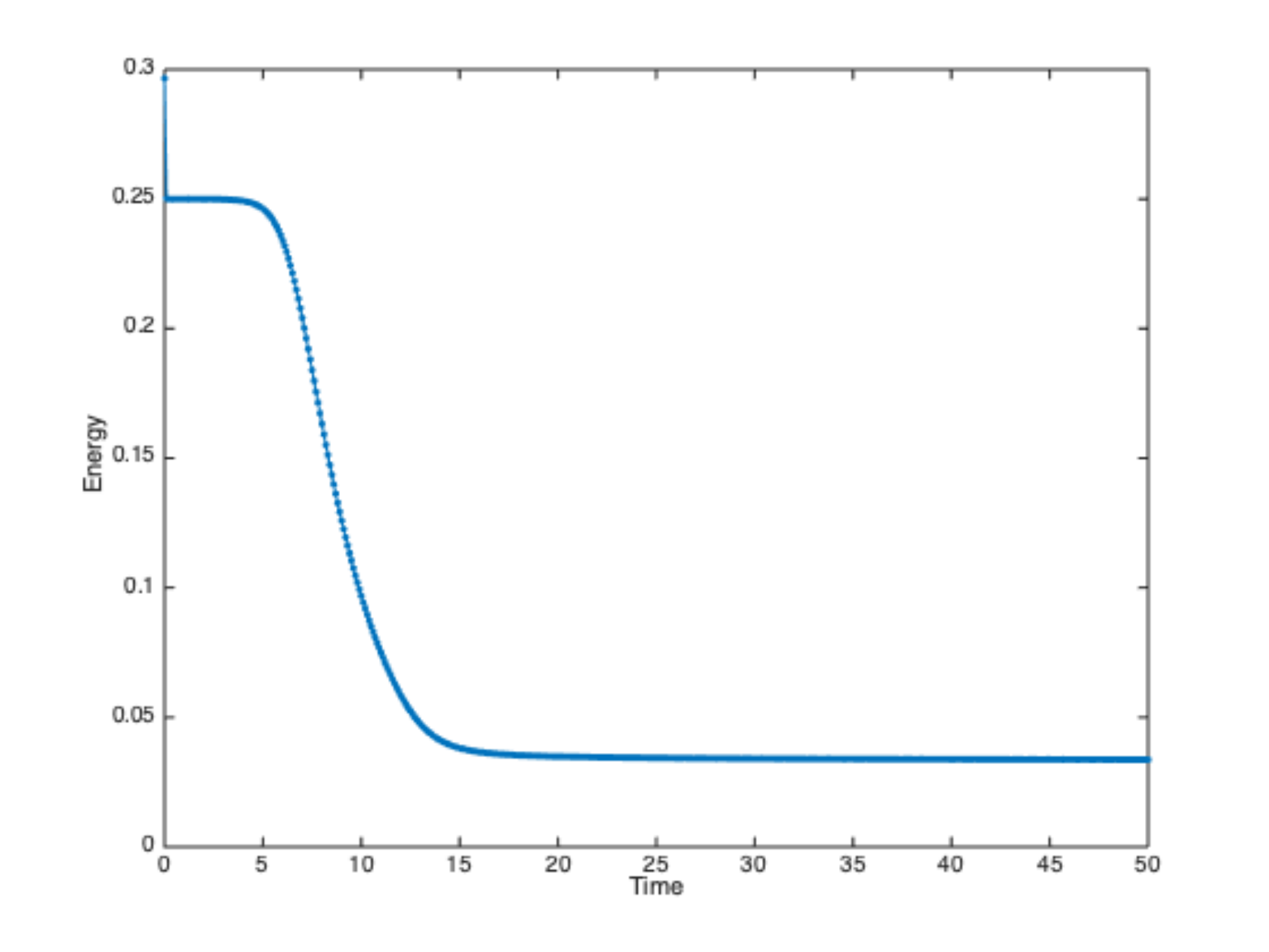}}
	{\includegraphics[width=0.45\textwidth]{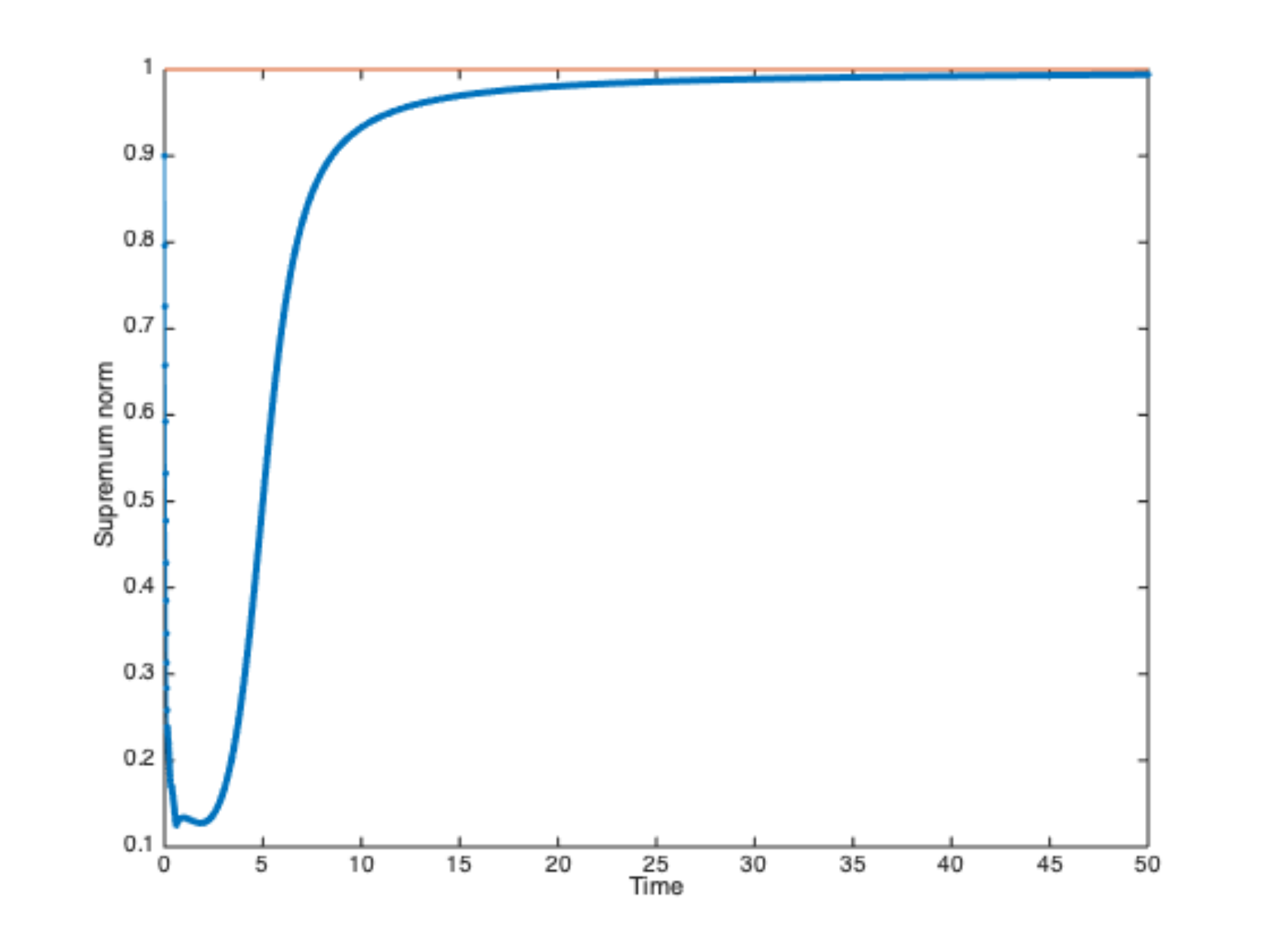}}\hspace{-0.3cm}
	{\includegraphics[width=0.45\textwidth]{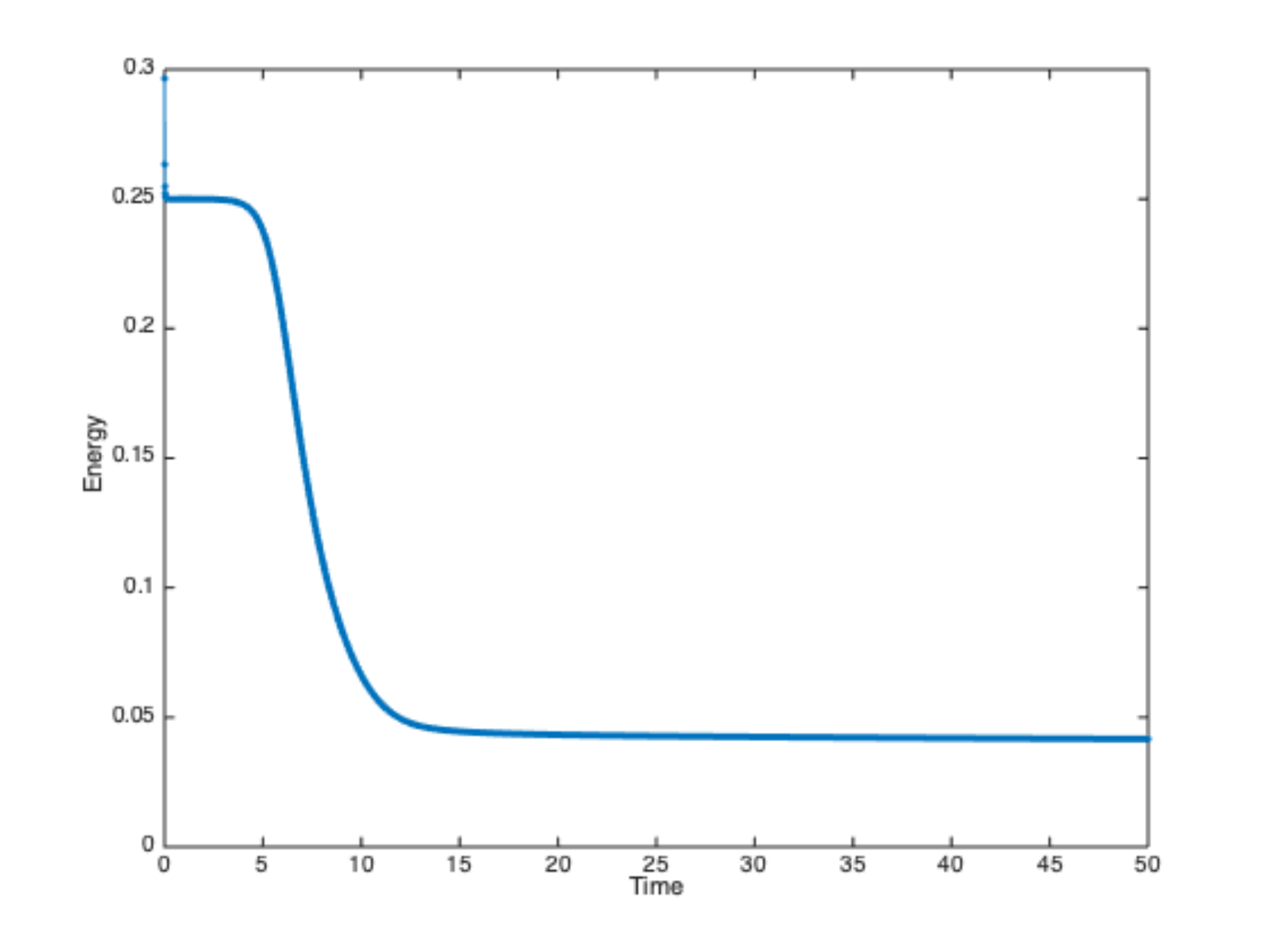}}
	\caption{Evolutions of the supremum norm (left) and the energy (right) of the numerical solutions for the 2D convective Allen-Cahn equation with the double-well potential. Top: $\tau$=0.1, bottom: $\tau$=0.01.}
	\label{fig-dw2D-phy}
\end{figure}

First,  the nonlinear reaction  $f=-F^\prime$ is chosen as the double-well potential case \eqref{dw}. The MBP bound constant is given by   $\beta=1$ and consequently and the stabilizing coefficient  $\kappa=1$ in view of the fact that $\max\limits_{|u|\leq1}\left|\widetilde{f}'(u)\right|=1$. Fig. \ref{fig-dw2D-phase} shows the snapshots of the numerical solutions at $t=0.1, 1, 8, 50$, respectively with different time step sizes $\tau=0.1$ and $\tau=0.01$. We can clearly see the ordering and coarsening phenomena, and the rotation effect caused by $\mathbf{v}(x,y,t)$ is well observed along the whole process. 
These two simulations produced by different time step sizes give us overall similar evolution processes.
The corresponding evolutions of the supremum norm and the classic free energy (defined in \eqref{energyorg}) of the numerical solutions are shown in Fig. \ref{fig-dw2D-phy}. We also observe that the  discrete MBP for the convective Allen-Cahn equation is  preserved perfectly along the time evolution. Moreover, the energy  also decays monotonically, although it is not guaranteed theoretically for the target problem. 

\begin{figure}[!ht]
	\centering
	{\includegraphics[width=0.35\textwidth]{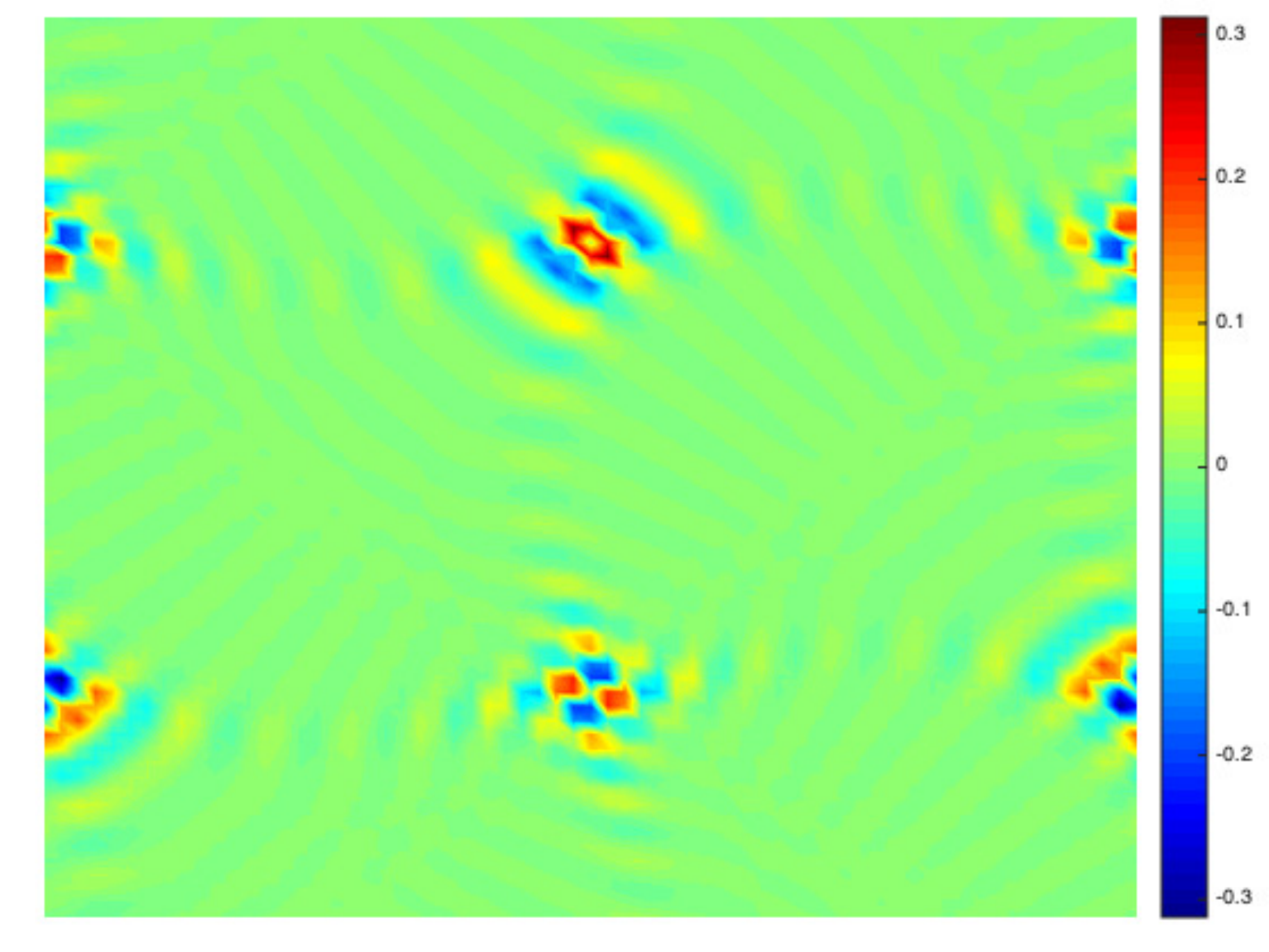}}
	{\includegraphics[width=0.35\textwidth]{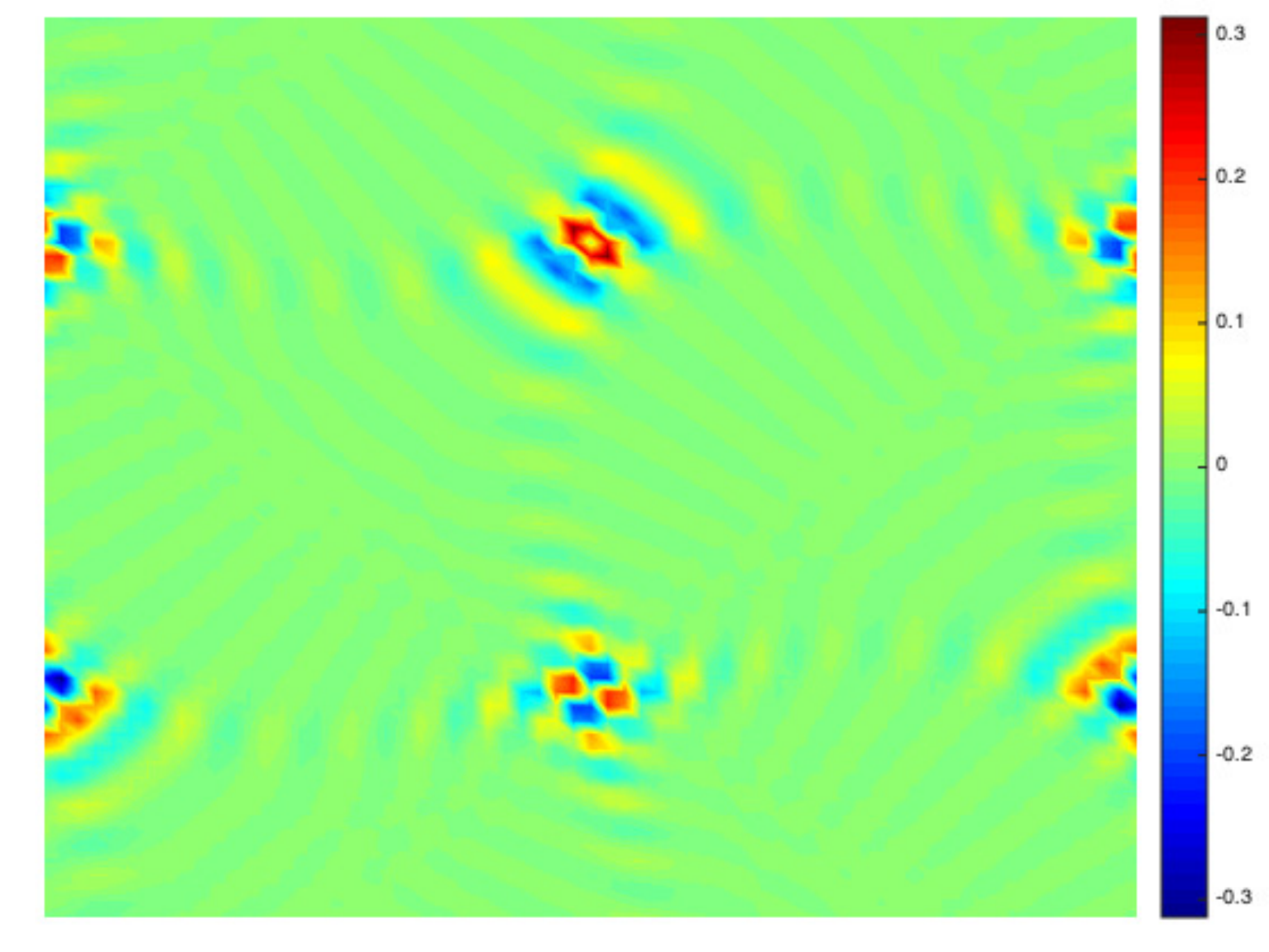}}
	{\includegraphics[width=0.35\textwidth]{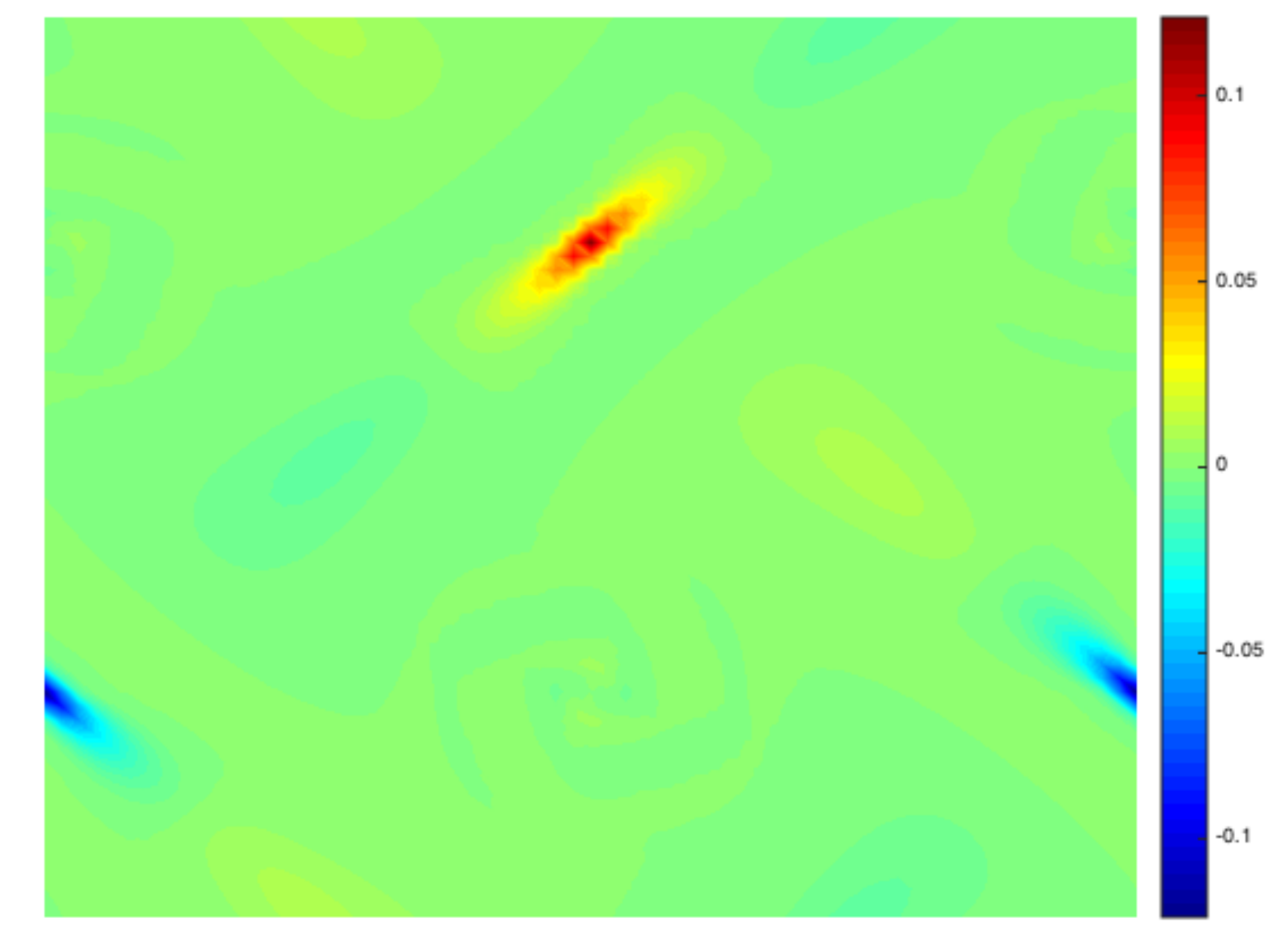}}
	{\includegraphics[width=0.35\textwidth]{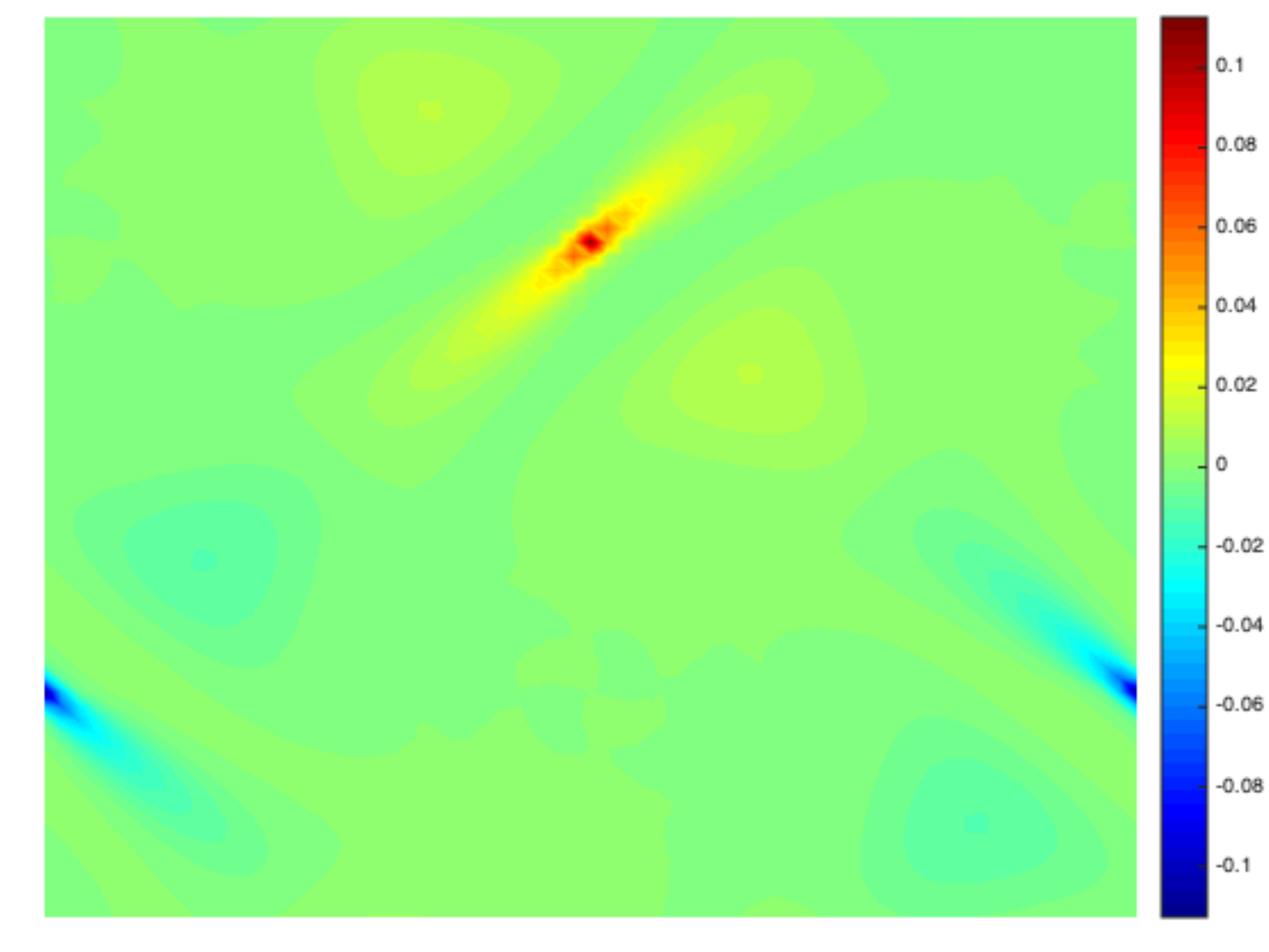}}
	{\includegraphics[width=0.35\textwidth]{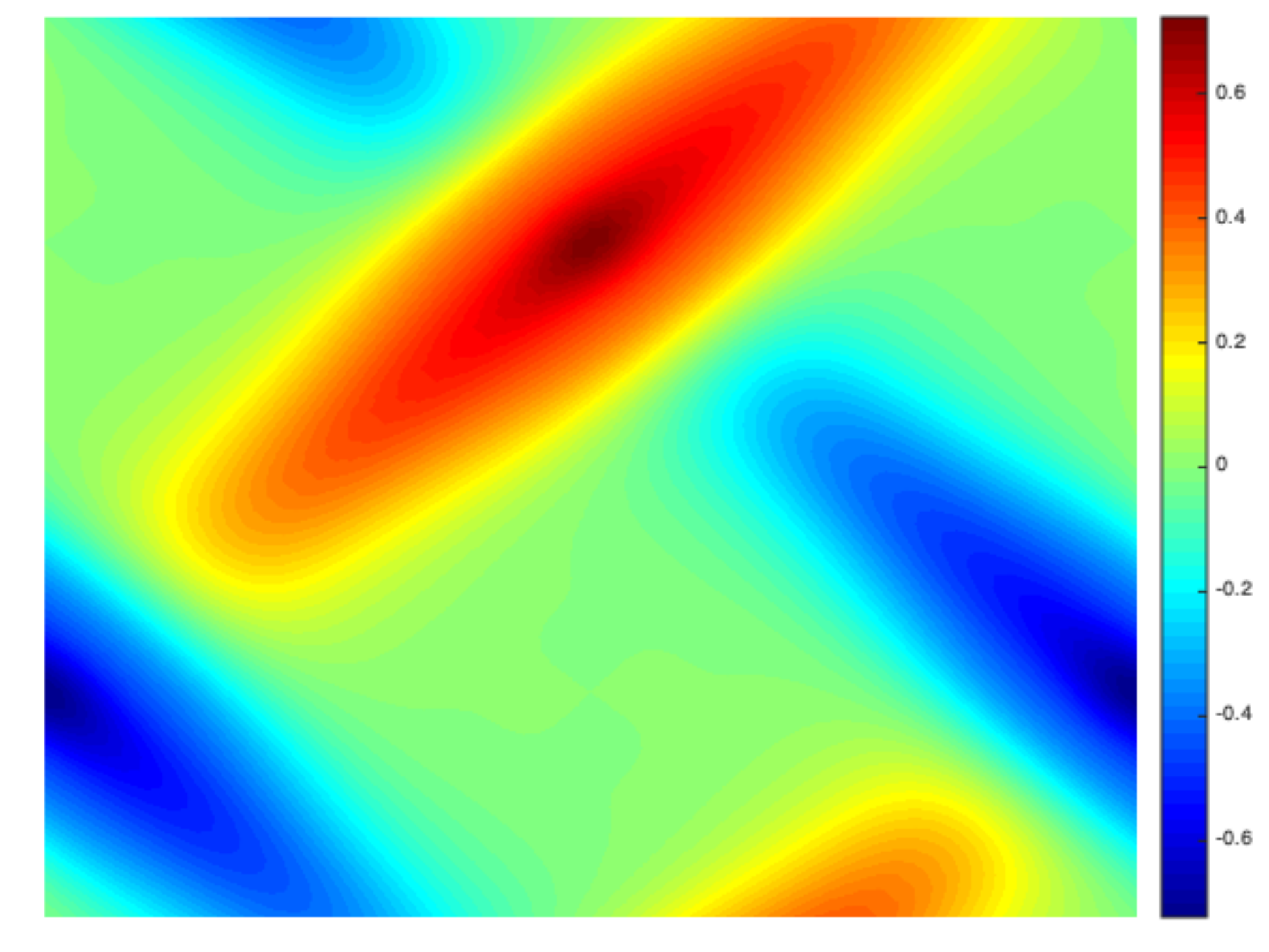}}
	{\includegraphics[width=0.35\textwidth]{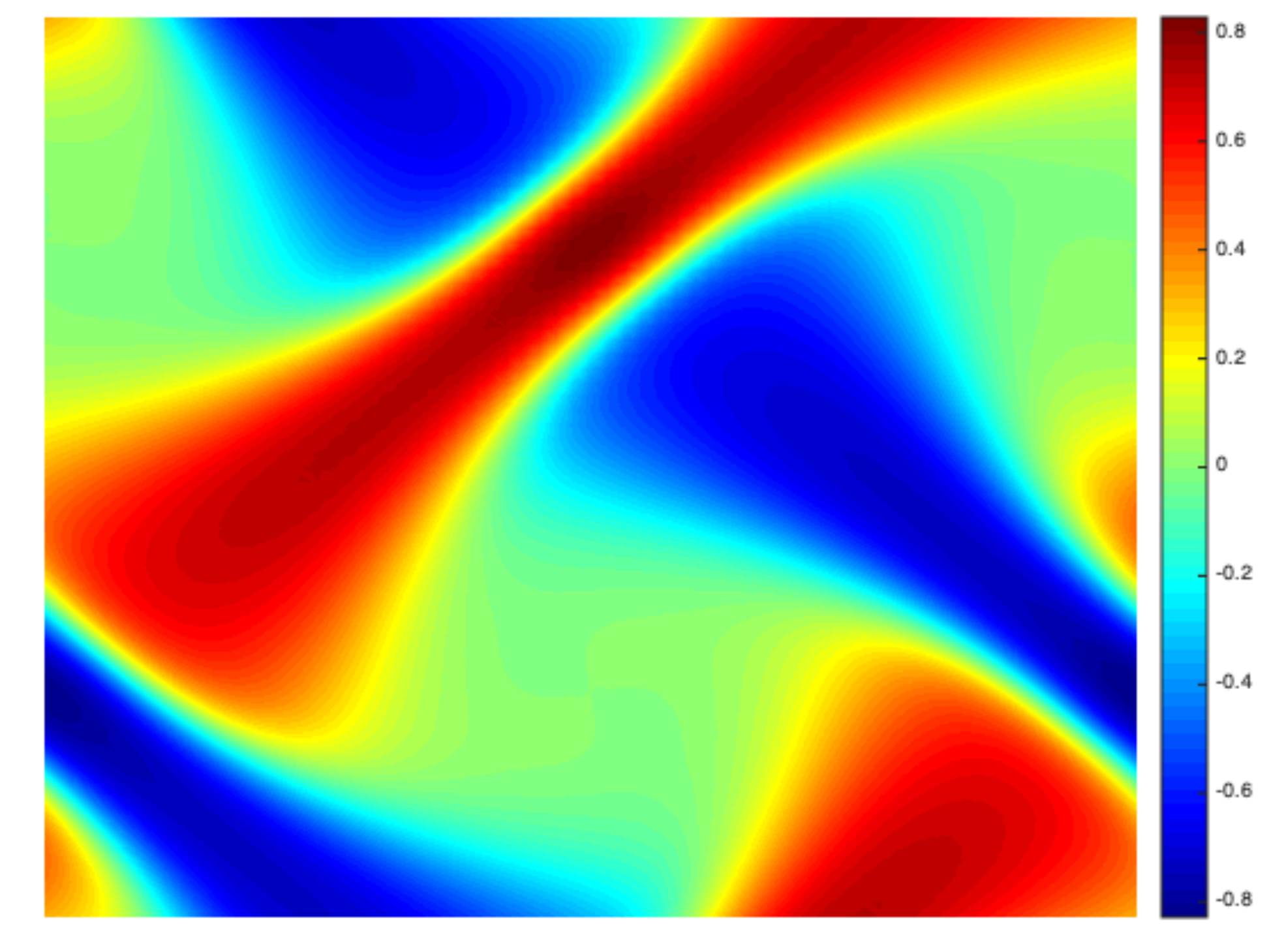}}
	{\includegraphics[width=0.35\textwidth]{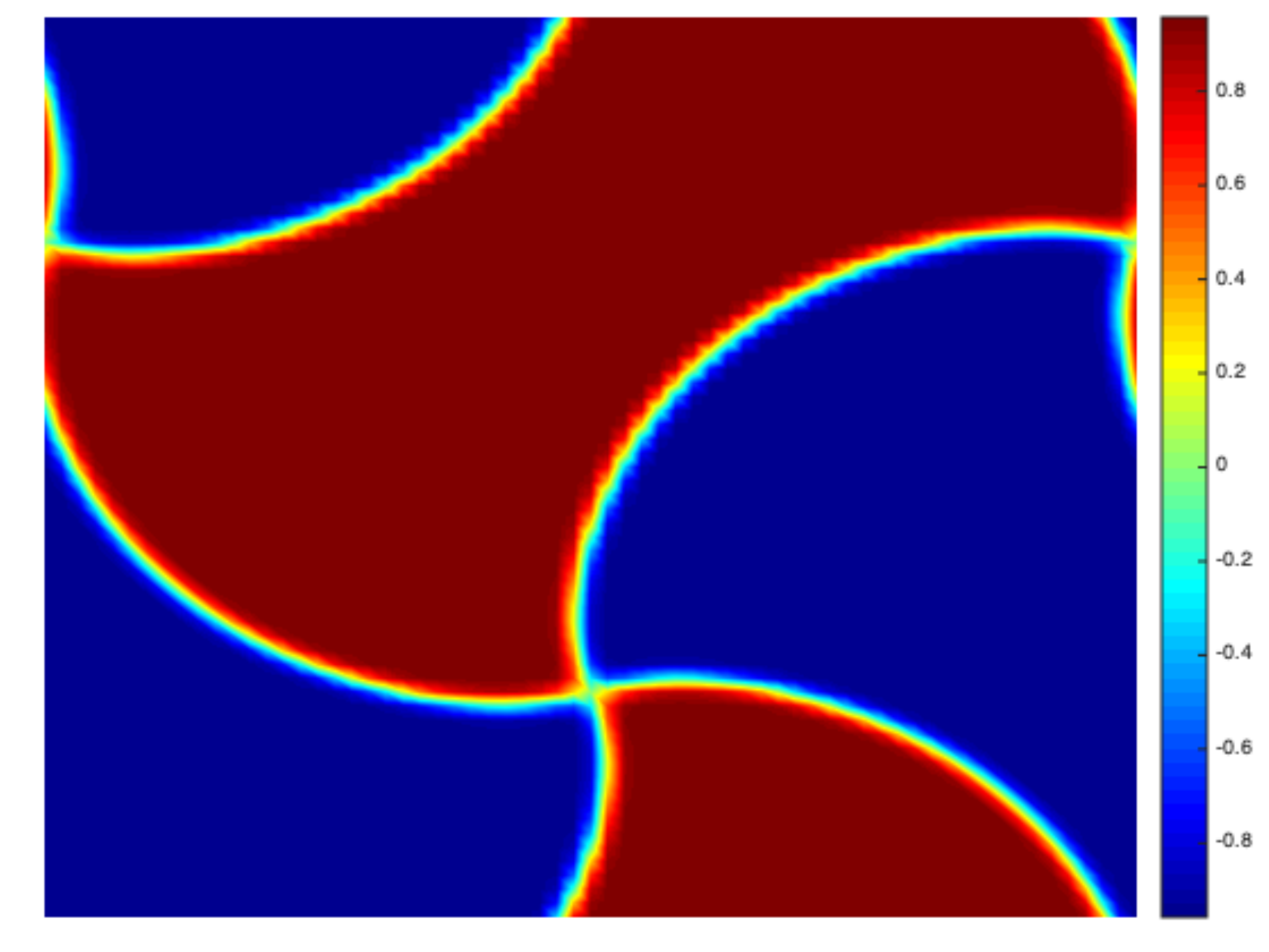}}
	{\includegraphics[width=0.35\textwidth]{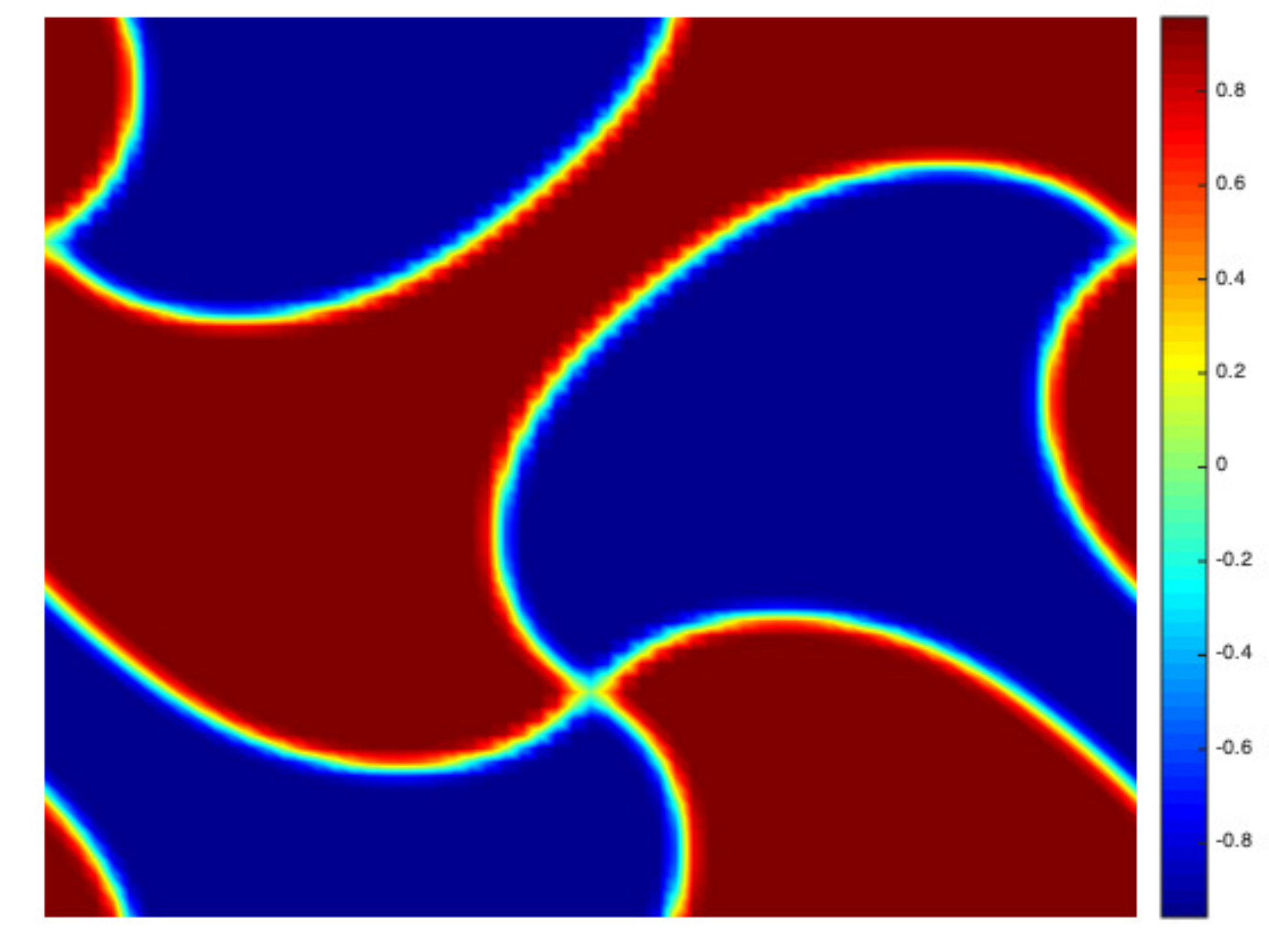}}
		\caption{The numerical solutions at $t$ = 0.1,1,8, 50, respectively (top to bottom) for the 2D convective Allen-Cahn equation with the Flory-Huggins potential.  Left column: $\tau$ = 0.1, right column: $\tau$ = 0.01.}
	\label{fig-fh2D-phase}
\end{figure}

\begin{figure}[!ht]
	\centering
	{\includegraphics[width=0.45\textwidth]{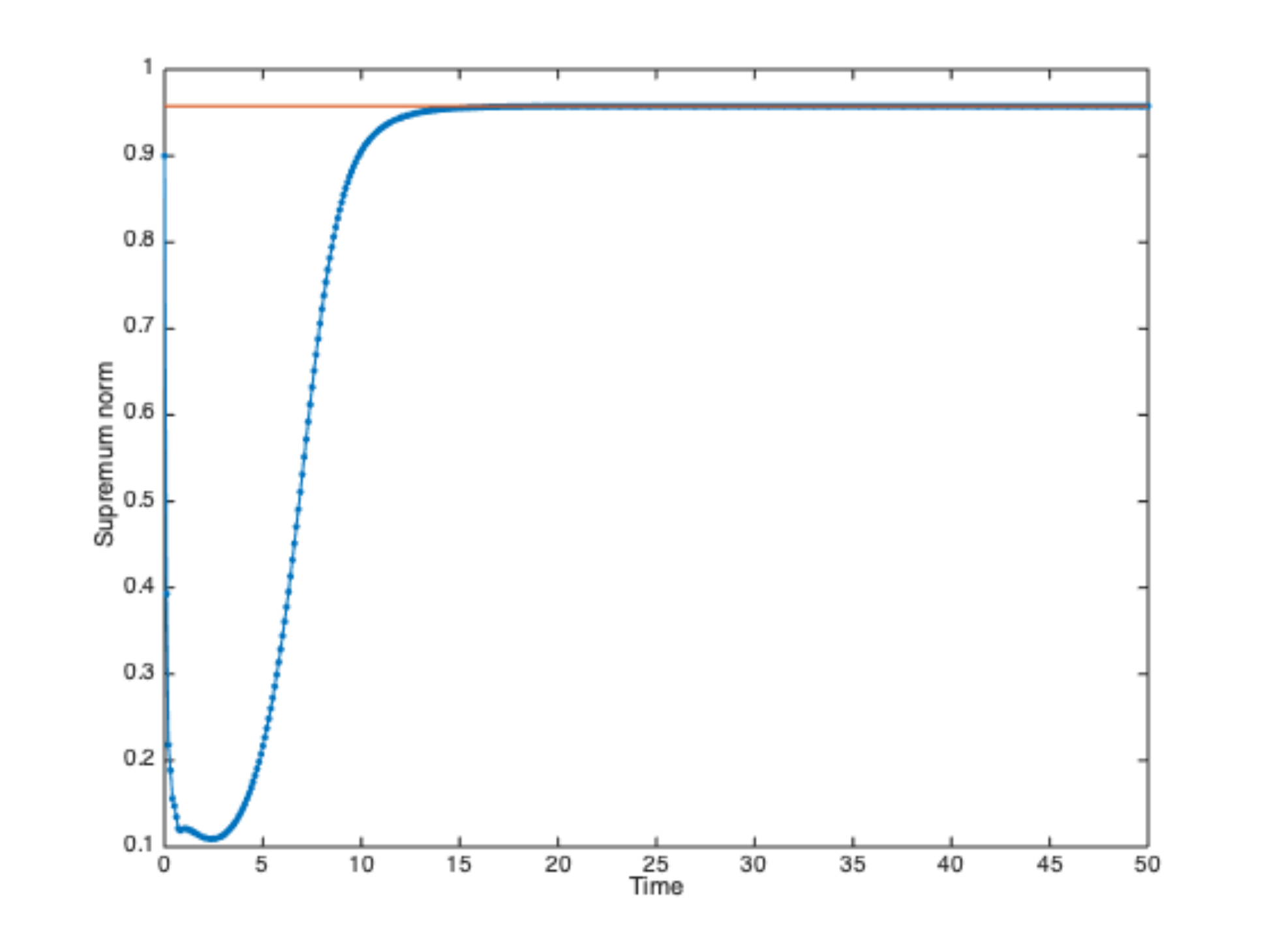}}\hspace{-0.3cm}
	{\includegraphics[width=0.45\textwidth]{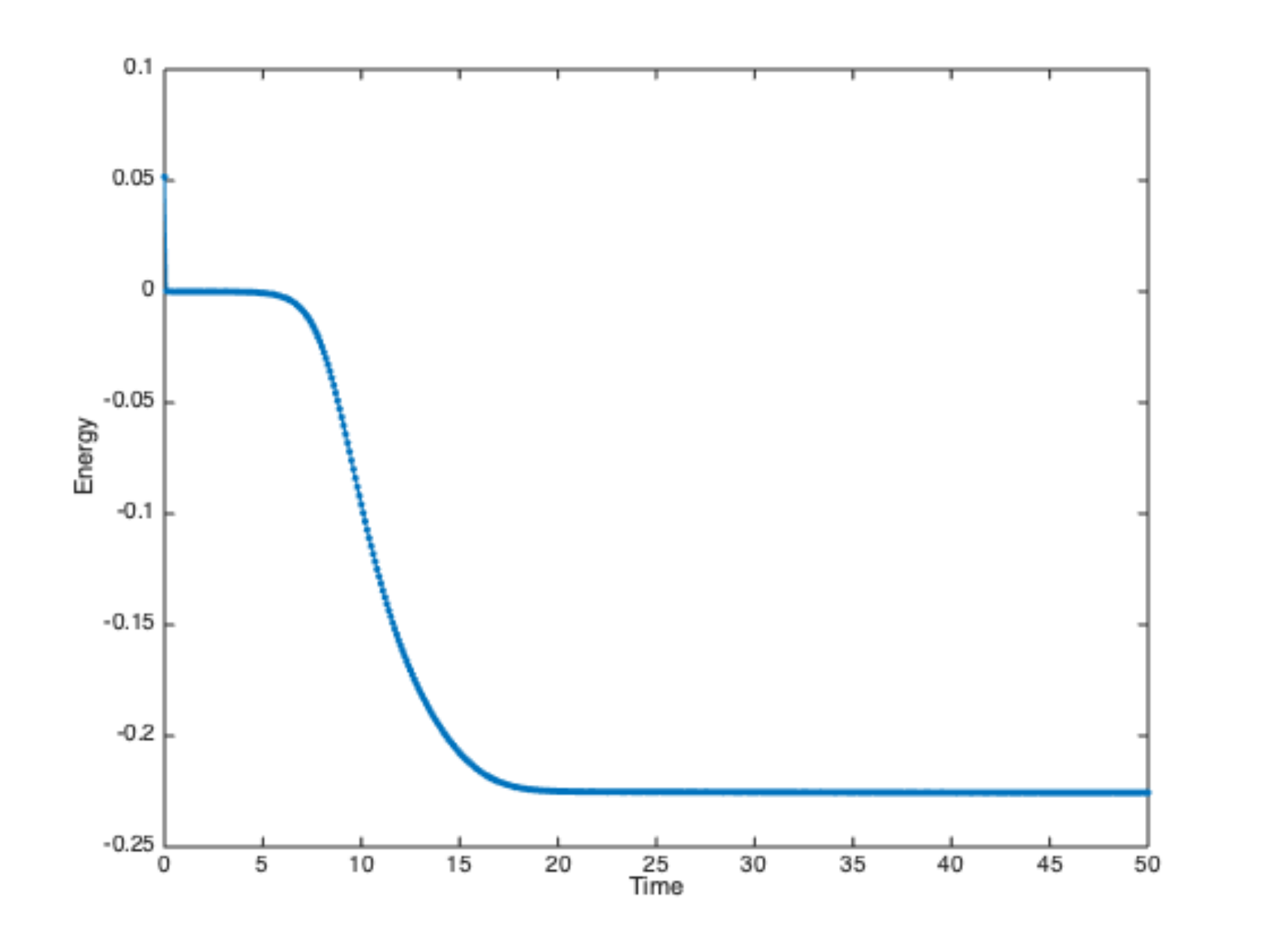}}
		{\includegraphics[width=0.45\textwidth]{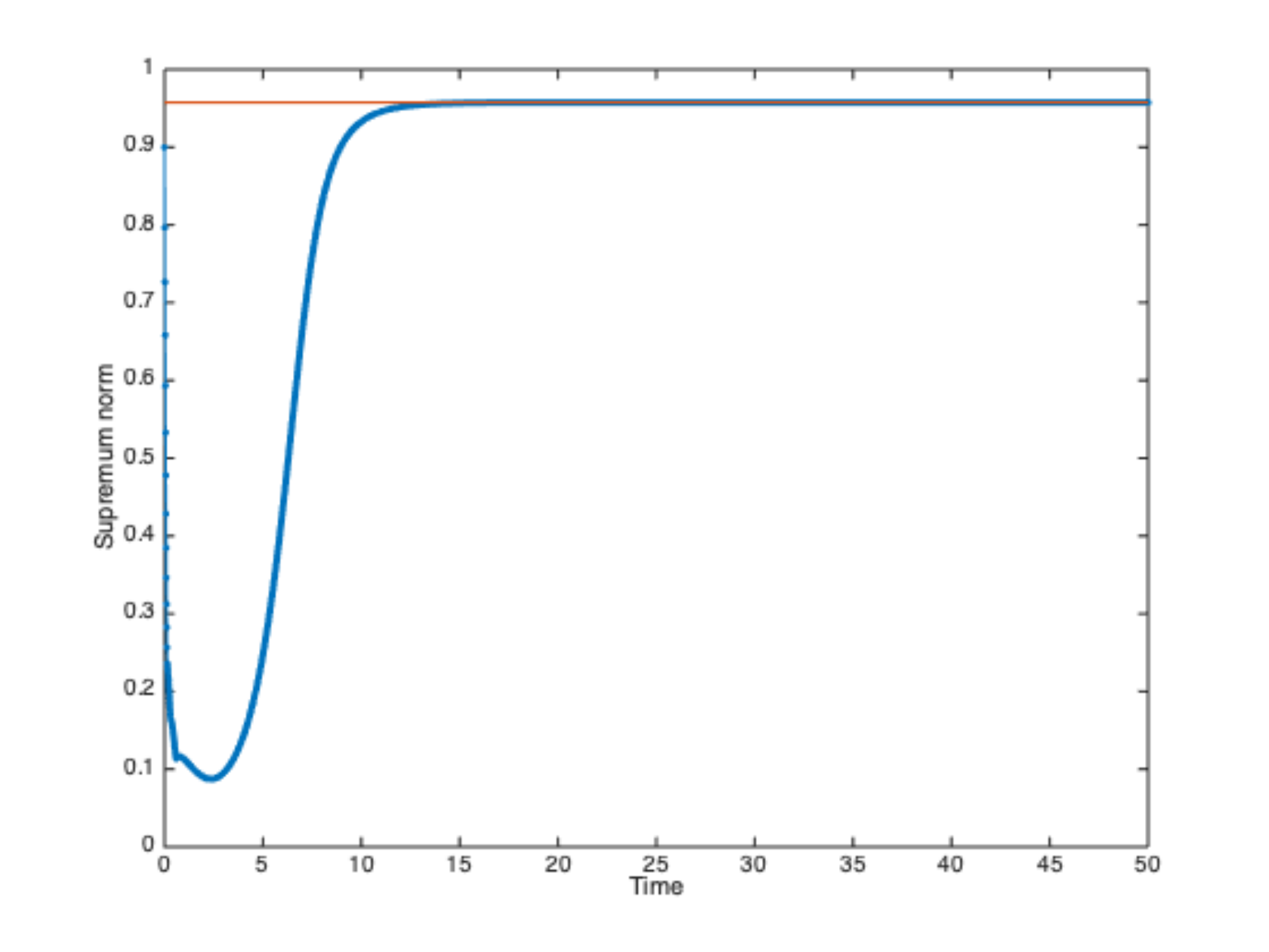}}\hspace{-0.3cm}
	{\includegraphics[width=0.45\textwidth]{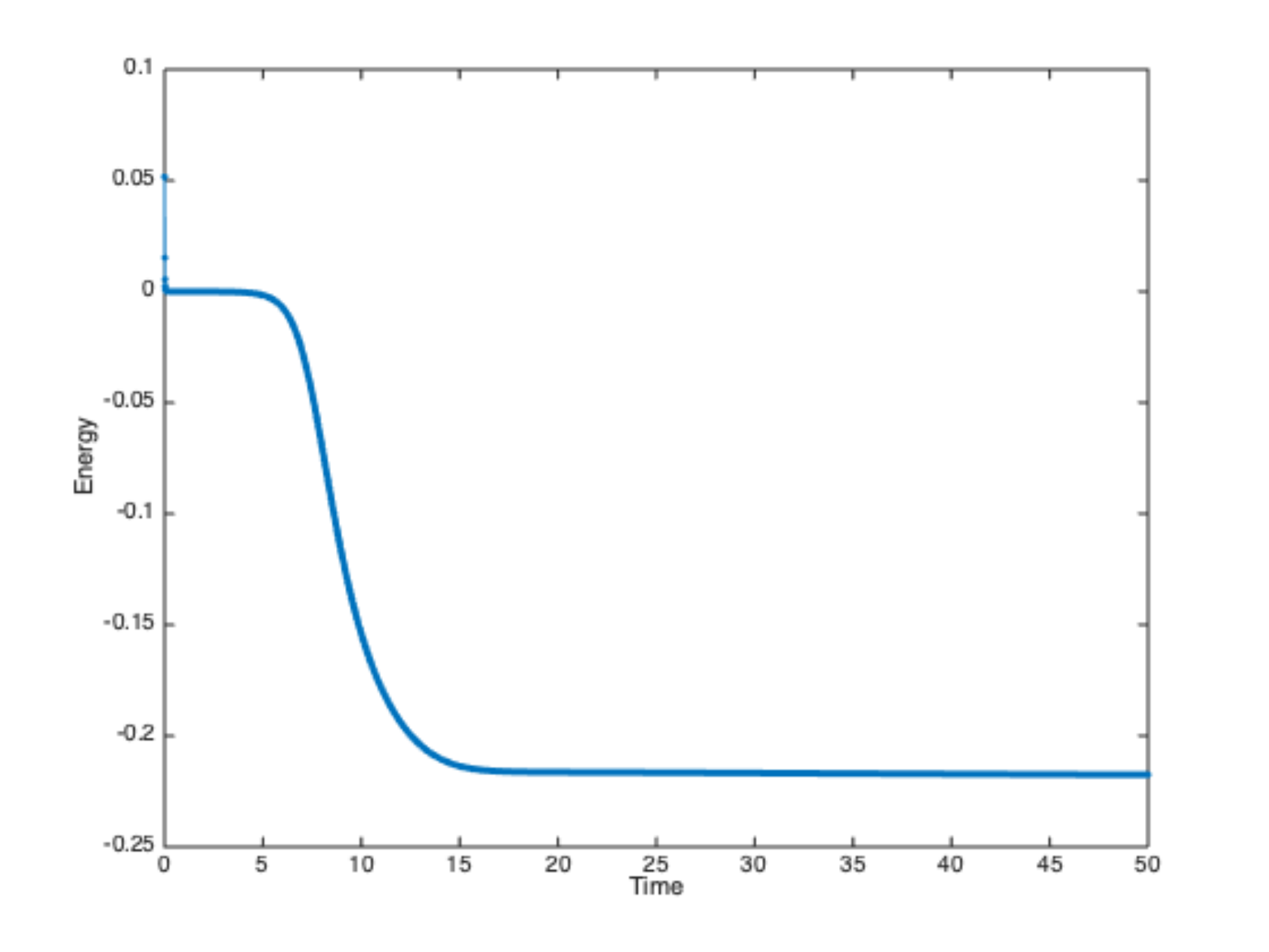}}
	\caption{Evolutions of the supremum norm (left) and the energy (right) of  the numerical solutions for the 2D convective Allen-Cahn equation with the Flory-Huggins potential. Top: $\tau$ = 0.1, bottom: $\tau$ = 0.01.}
	\label{fig-fh2D-phy}
\end{figure}

Next, $f=-F^\prime$ is chosen as the Flory-Huggins potential case \eqref{fh} with  $\theta=0.8$ and $\theta_c=1.6$.  The MBP bound constant is now given by   $\beta \approx 0.9575$, and consequently the stabilizing coefficient is set to be $\kappa=1$ since $ \max\limits_{|u|\leq\beta}\left|\widetilde{f}'(u)\right|\approx0.9801$ so that \eqref{coef} is satisfied in this case. Fig. \ref{fig-fh2D-phase} shows the snapshots of the numerical solutions at $t$=0.1, 1, 8, 50, respectively with $\tau=0.1$ and $\tau=0.01$, and  the corresponding time evolutions of the supremum norm and the energy are  presented in Fig. \ref{fig-fh2D-phy}. Again, the ordering and coarsening phenomena and the rotation effect caused by $\mathbf{v}(x,y,t)$ are clearly observed along the process. It is also clear that the energy decays monotonically and the discrete MBP for the convective Allen-Cahn equation is  well preserved numerically. The two simulations  by different time step sizes again produce similar evolution processes.

\subsection{Convective test}
\begin{figure}[!ht]
	\centering
	{\includegraphics[width=0.38\textwidth]{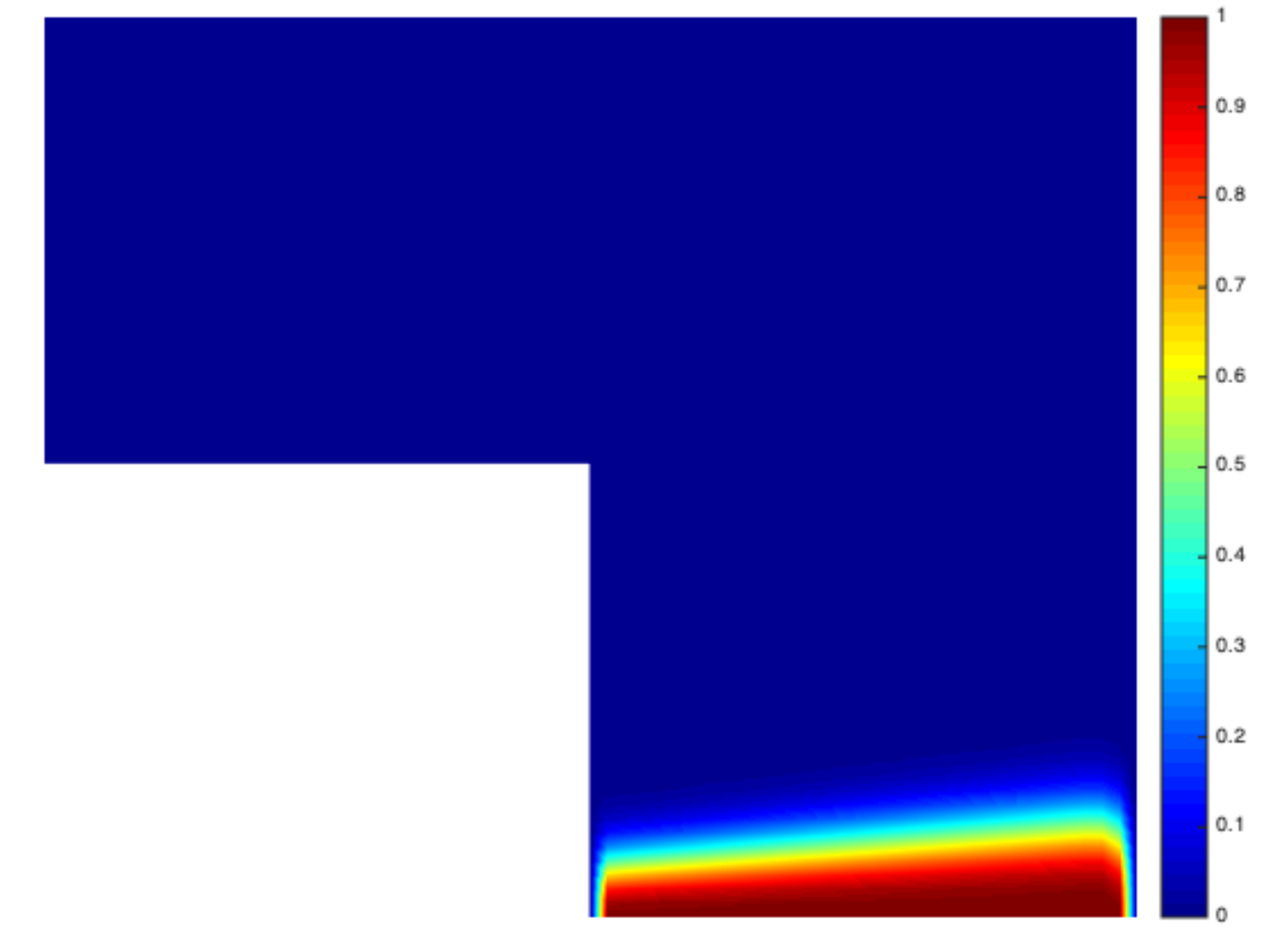}}
	{\includegraphics[width=0.38\textwidth]{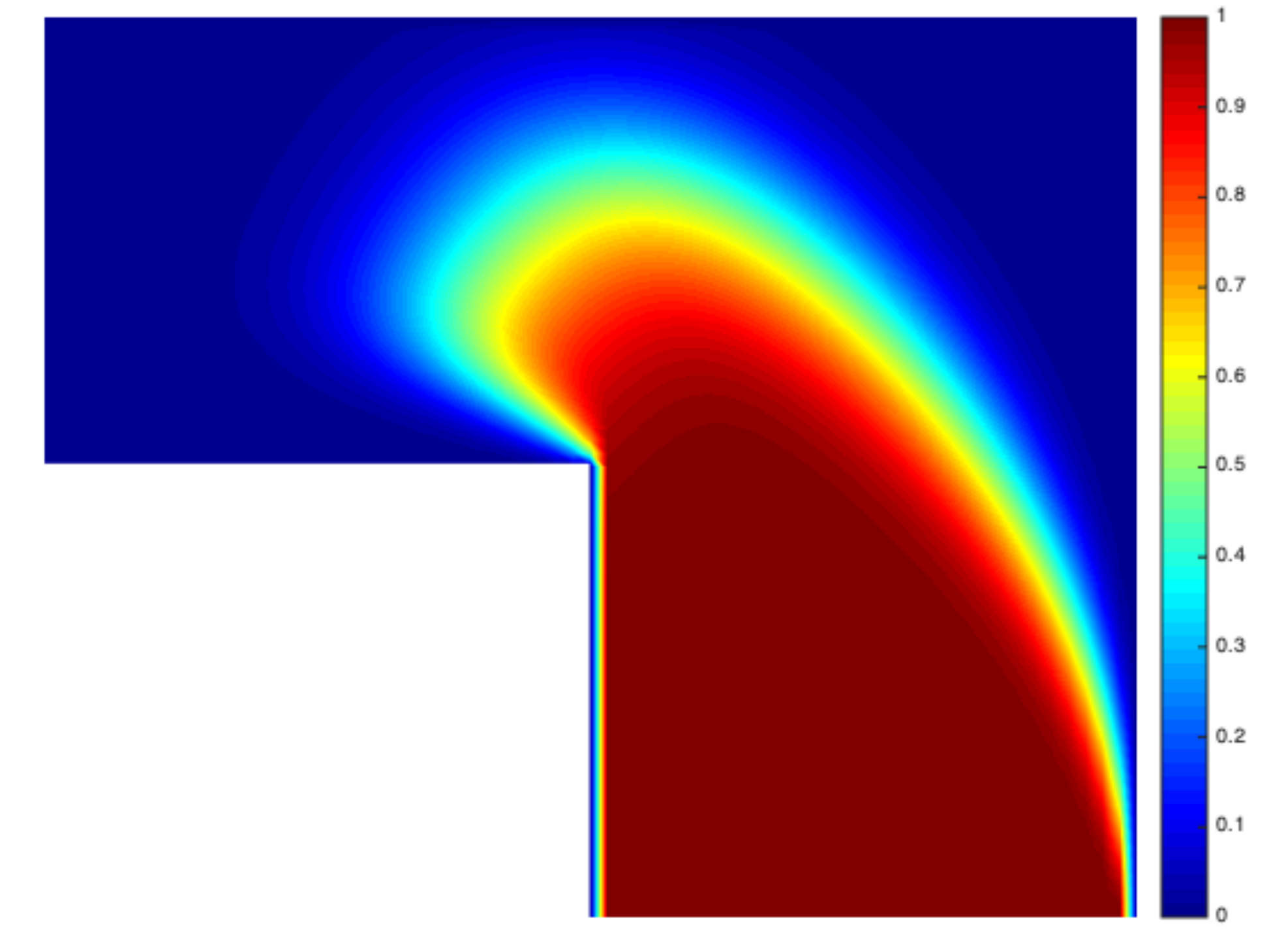}}
	{\includegraphics[width=0.38\textwidth]{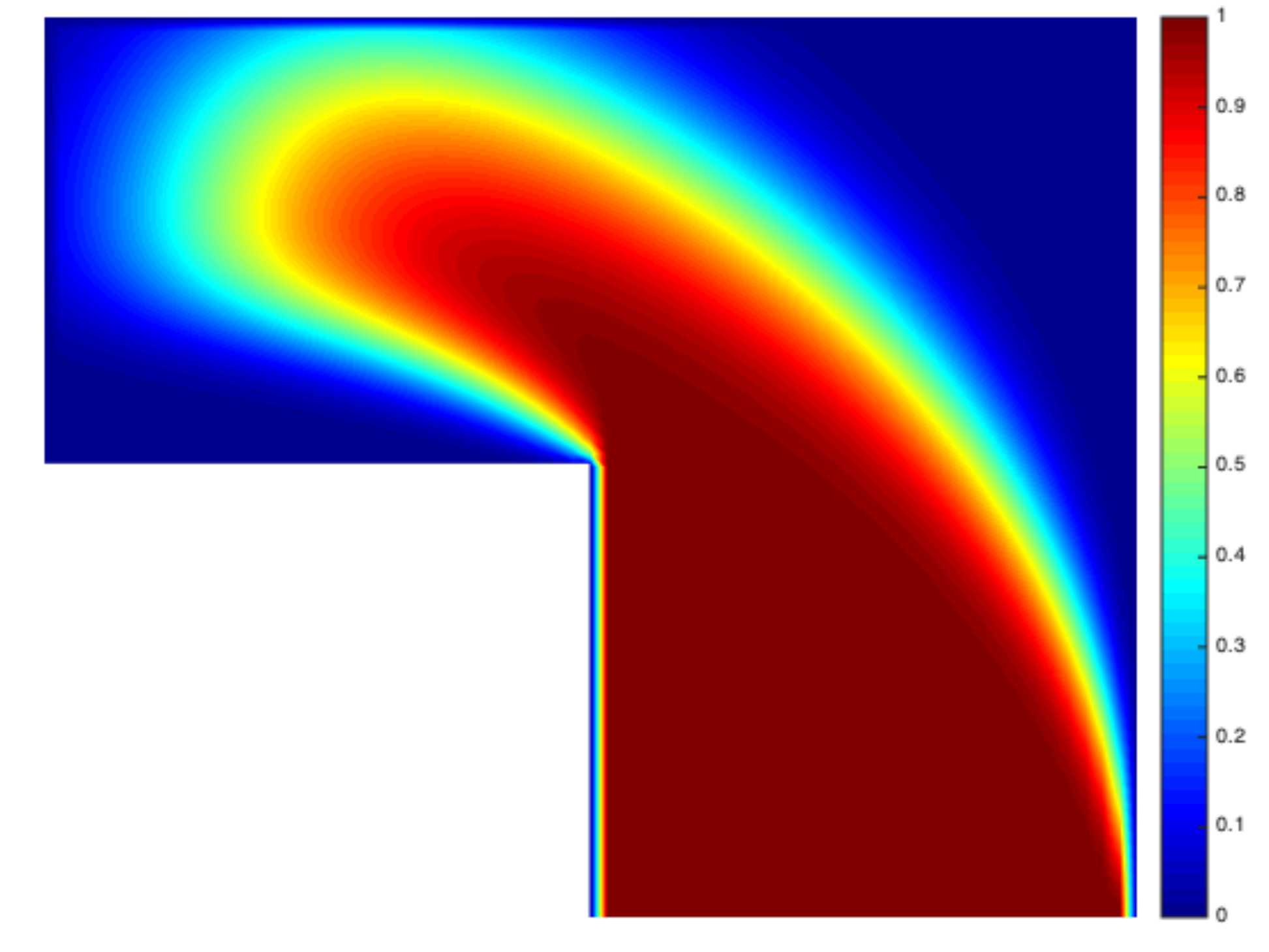}}
	{\includegraphics[width=0.38\textwidth]{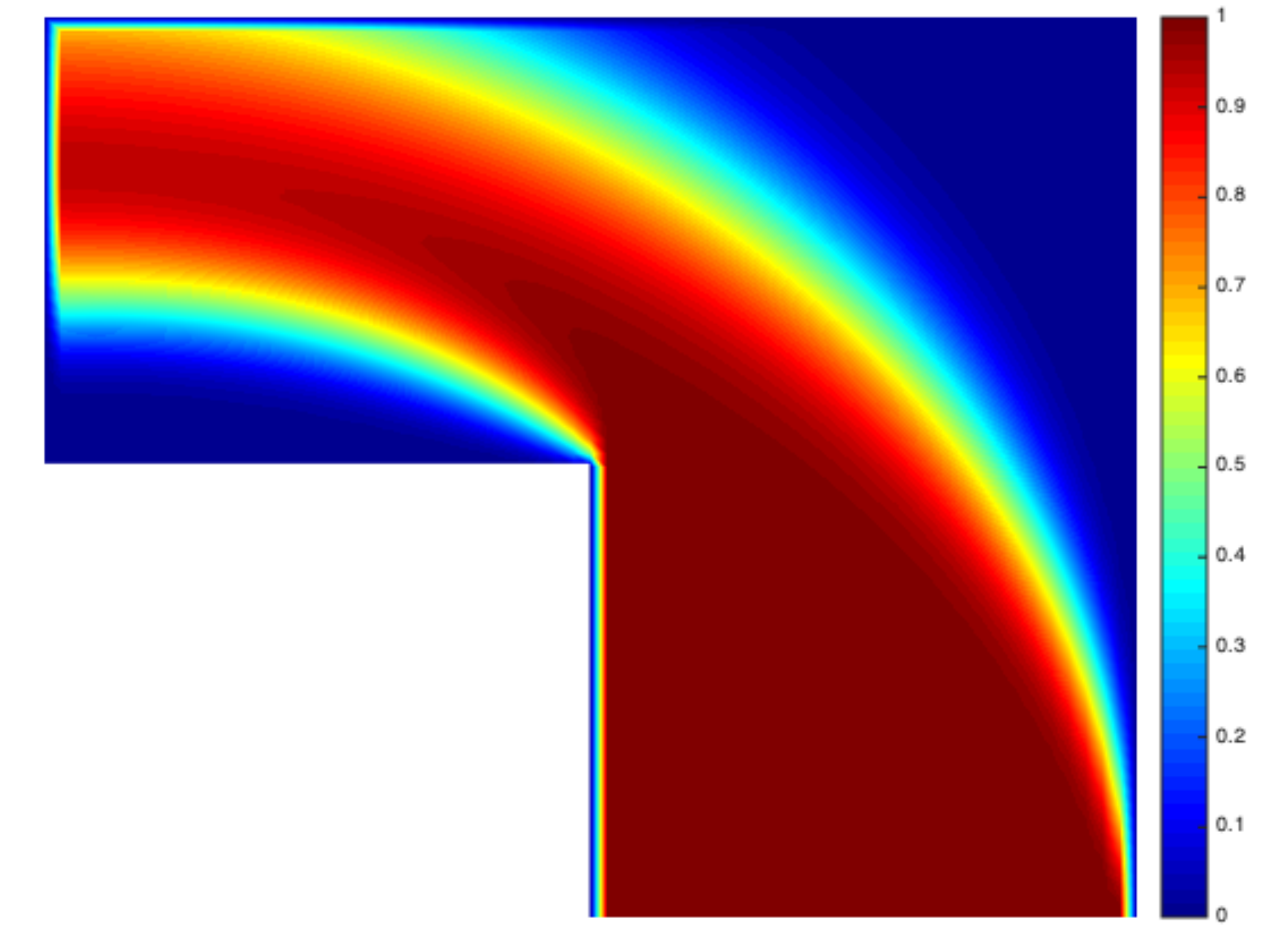}}
	\caption{The numerical solutions  at $t$=0.1,1,1.3,10, respectively (top to bottom and left to right) for the 2D convective Allen-Cahn equation with the double-well potential.}
	\label{fig-dw2D-LD}
\end{figure}

Now we consider the 2D Allen-Cahn equation \eqref{gac} with the  mobility $M(u)=1-u^2$, the velocity field  $\mathbf{v}(x,y,t)=[y,-x]^T$ and $\epsilon=0.01$ in the L-shape domain $(x,y)\in\Omega= (0,1)^2/[0,0.5]^2$ subjected to the  Dirichlet boundary condition specified on the boundary $\partial\Omega$ as
\begin{align*}
u(x,y,t)=\left\{
\begin{array}{ll}
1,&\quad y=0;\\
0,& \quad\text{otherwise}.
\end{array}
\right.
\end{align*}
The nonlinear reaction $f=-F^\prime$ is given by the double-well potential  \eqref{dw}. The initial data $u_0(x,y)=0$ except on the domain boundary $y=0$ where $u_0(x,0)=1$  (compatible with the boundary condition). For this special example \cite{QCRF14}, the solution $u(x,y,t)$ is always located in the interval $[0,1]$ according to the result $0\leq N(\xi)\leq\kappa$ for any $\xi\in[0,1]$ in Lemma \ref{lemma-stab-nonlinear}.
Fig. \ref{fig-dw2D-LD} shows the snapshots of the numerical solutions generated by the ETDRK2 scheme  with $\kappa=1$ at $t$=0.1,1,1.3,10, respectively. The ordering and coarsening phenomena as well as the counter-clockwise  rotation effect due to the convective term are clearly observed. Moreover, we see from Fig. \ref{fig-dw2D-LD} that the numerical solutions  remain between $0$ and $1$ very well at these times. 

\subsection{3D simulations}
In this subsection, some 3D simulations are performed for the convective Allen-Cahn equation \eqref{gac} with $\epsilon=0.01$ under the periodic boundary condition. The computational domain is set to be $\Omega=(-0.5,0.5)^3$. 
The initial configuration  is the quasi-uniform state 
$u_0(\cdot)=0.9\,\text{rand}(\cdot)$, where $\text{rand}(\cdot)$ generates the random numbers from uniform distribution on $[-1,1]$.  The mobility function is $M(u)\equiv1$ and the velocity field is $\mathbf{v}=[1,1,1]^T$. The time step is set to be $\tau=0.01$ and  the mesh size  is adopted with $h=1/128$. Again the ETDRK2 scheme is used.

\begin{figure}[!ht]
	\centering
	{\includegraphics[width=0.45\textwidth]{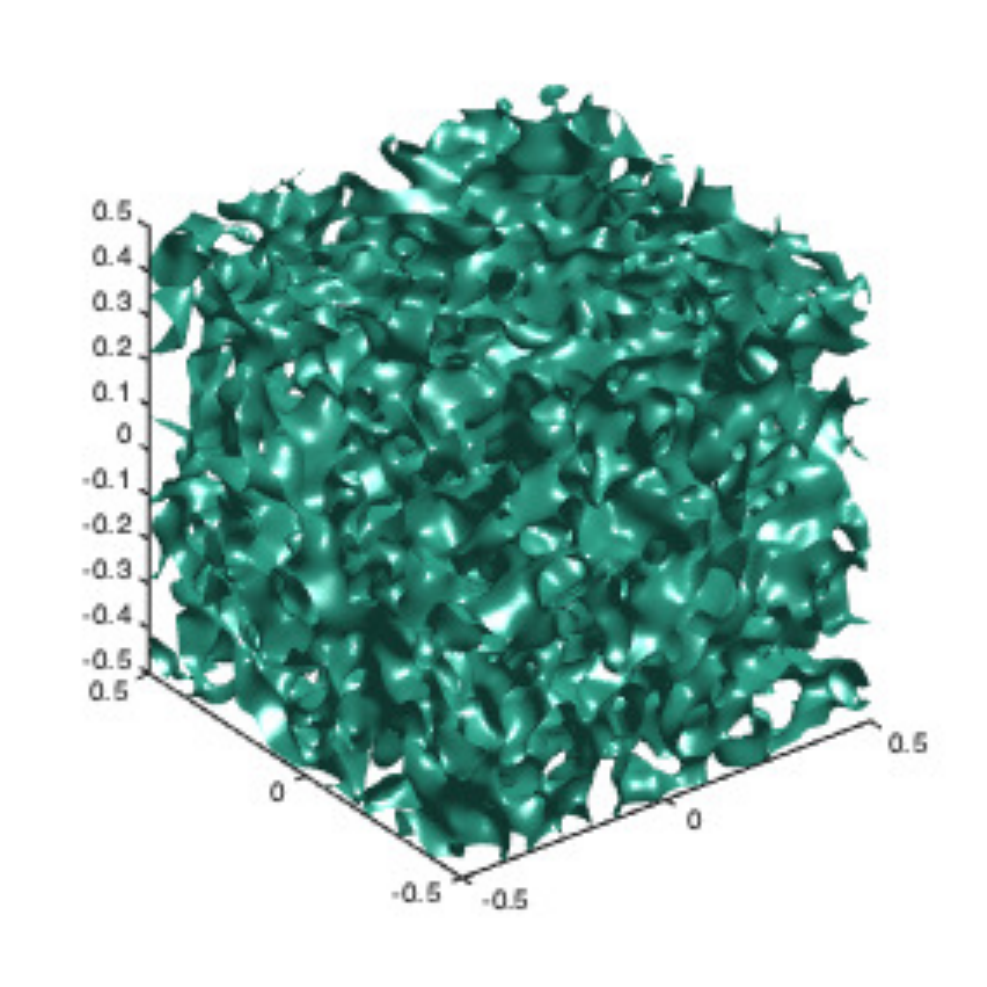}}\hspace{-0.1cm}
	{\includegraphics[width=0.45\textwidth]{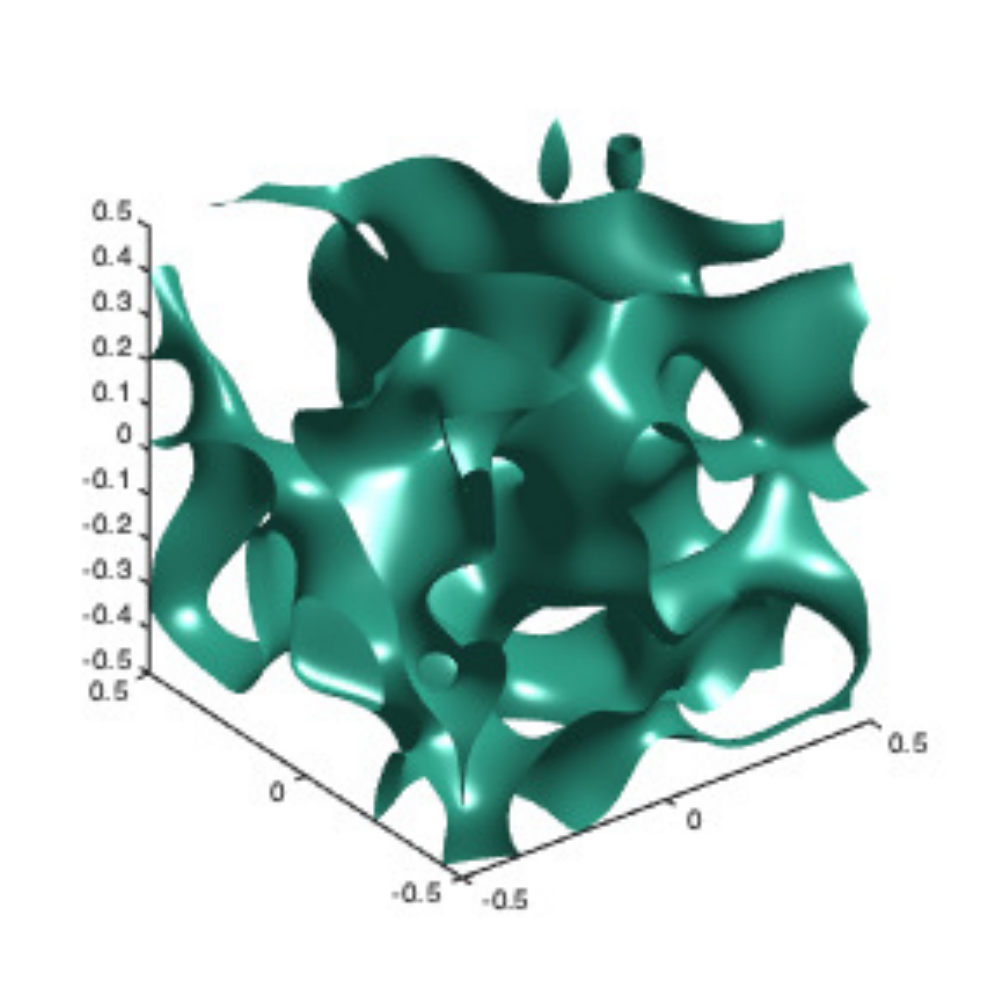}}\\ \vspace{-0.5cm}
	{\includegraphics[width=0.45\textwidth]{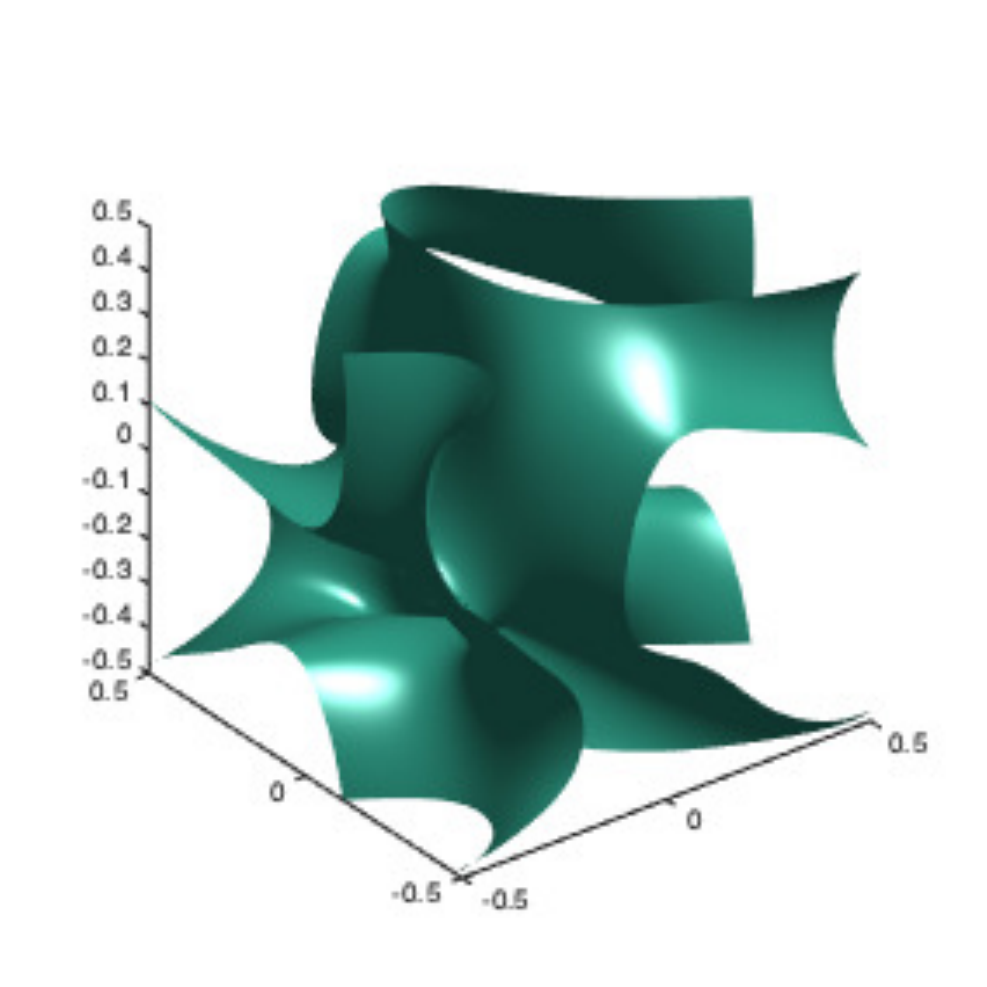}}\hspace{-0.1cm}
	{\includegraphics[width=0.45\textwidth]{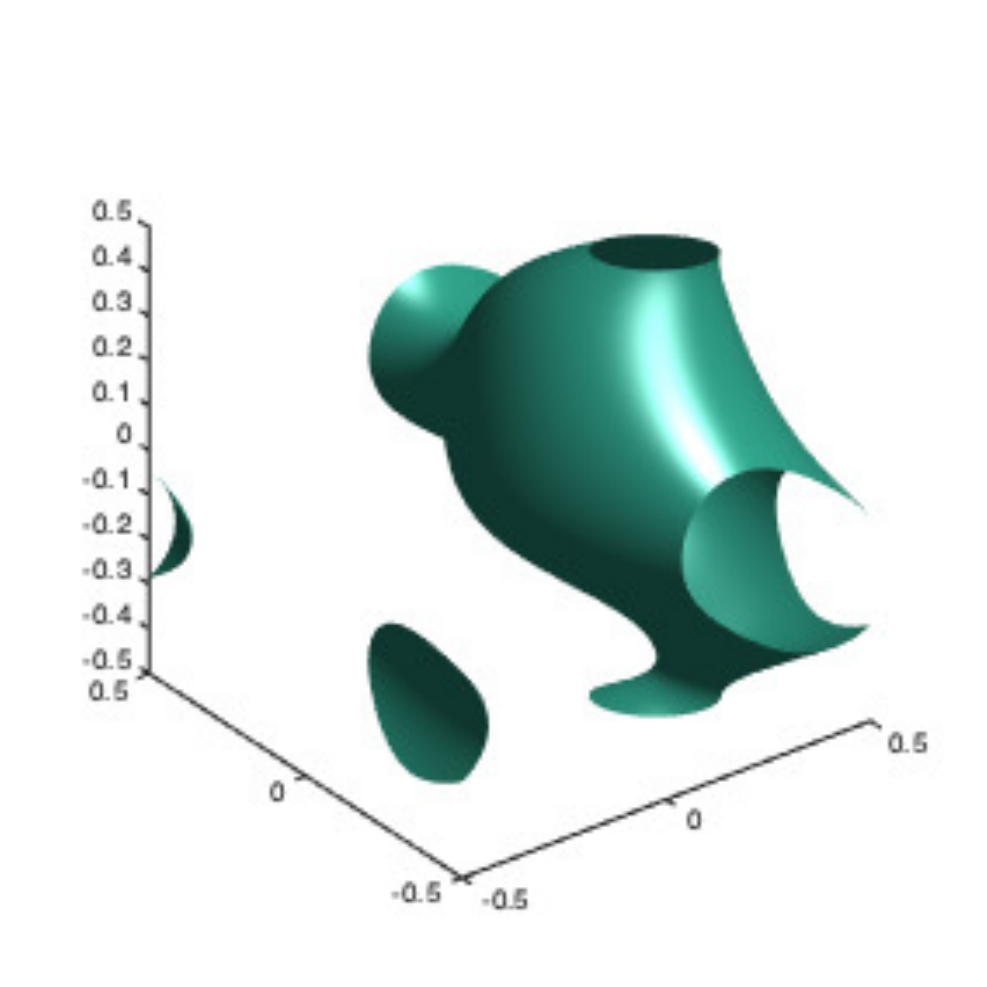}}
	\vspace{-0.2cm}
	\caption{The simulated phase structures at $t$ = 0.1,1,5,8, respectively (top to bottom and left to right) for the 3D convective Allen-Cahn equation with the double-well potential.}
	\label{fig-dw3D-phase}
\end{figure}

\begin{figure}[!ht]
	\centering
	{\includegraphics[width=0.45\textwidth]{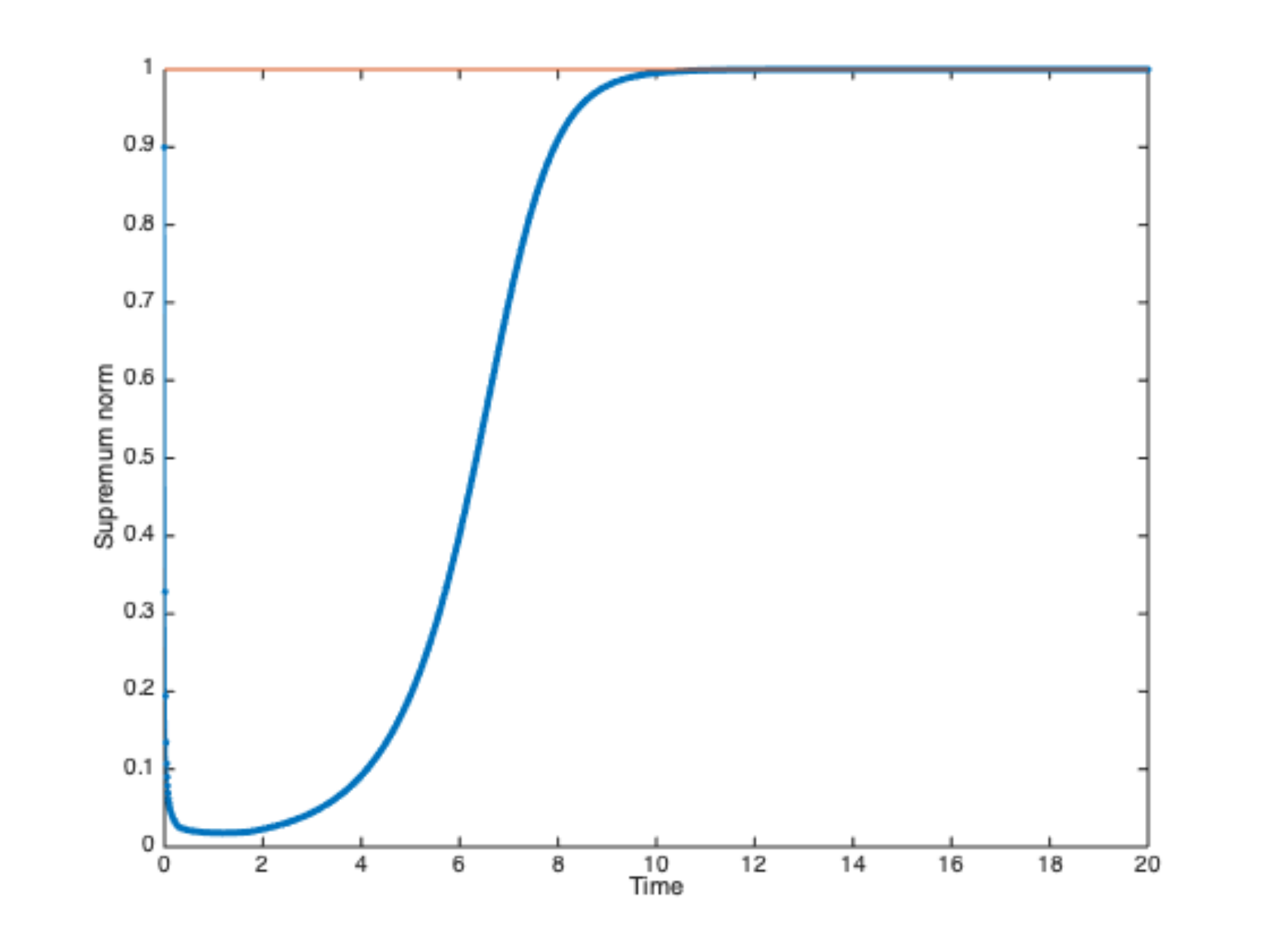}}\hspace{-0.3cm}
	{\includegraphics[width=0.45\textwidth]{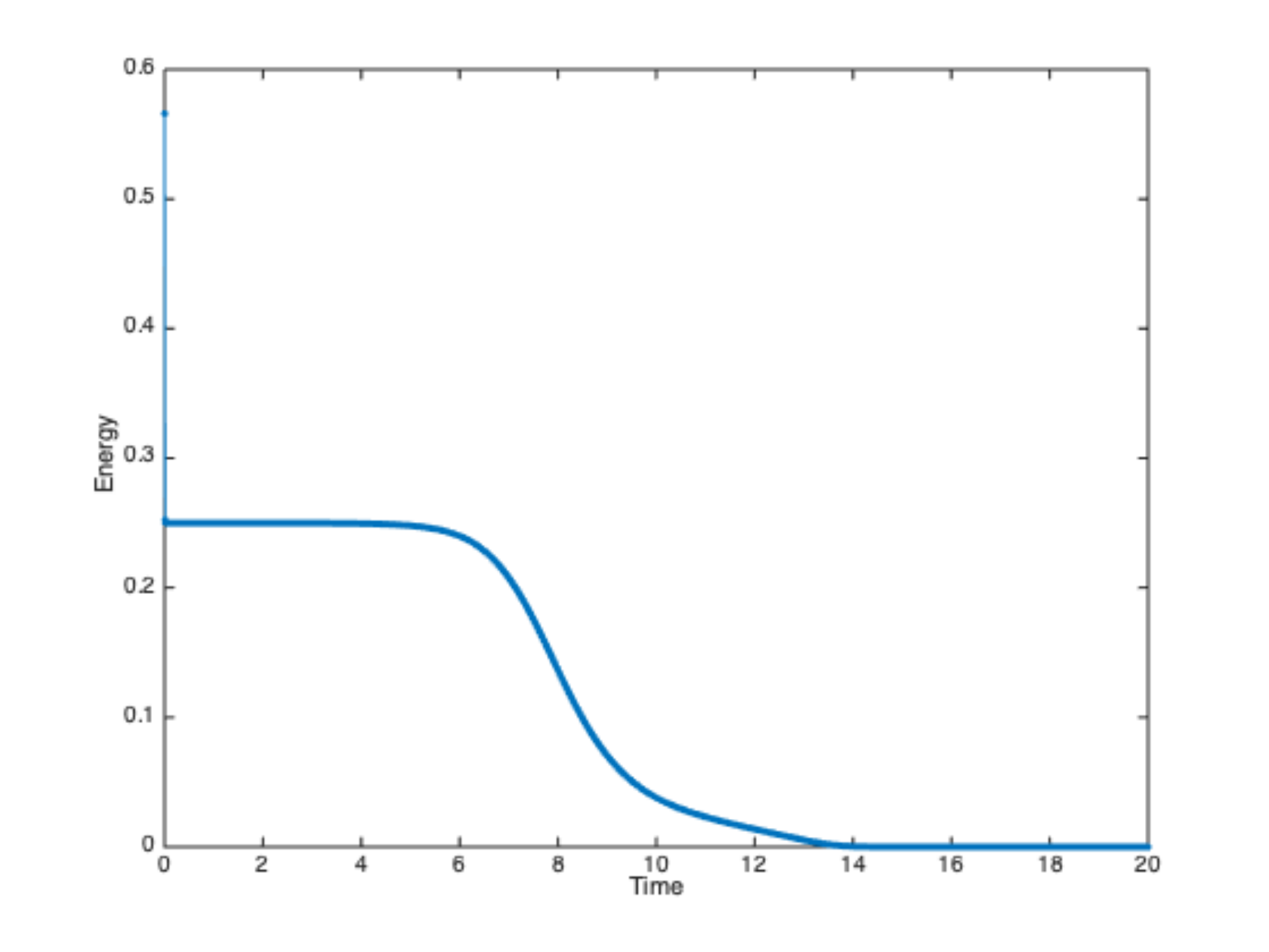}}
	\caption{The   evolutions of the supremum norm (left) and the energy (right) of numerical solutions for the 3D convective Allen-Cahn equation with the double-well potential.}
	\label{fig-dw3D-phy}
\end{figure}

First,  we choose the double-well potential case\eqref{dw} and the corresponding stabilizing coefficient is $\kappa=2$. Fig. \ref{fig-dw3D-phase} shows the phase structures  of the numerical solutions at $t$=0.1, 1, 5, 8, respectively. The time evolutions of the supremum norm and the energy of numerical solutions are presented in Fig. \ref{fig-dw3D-phy}. It is observed that the MBP for the convective Allen-Cahn equation is numerically preserved perfectly and the energy decays monotonically.

\begin{figure}[!ht]
	\centering
	{\includegraphics[width=0.45\textwidth]{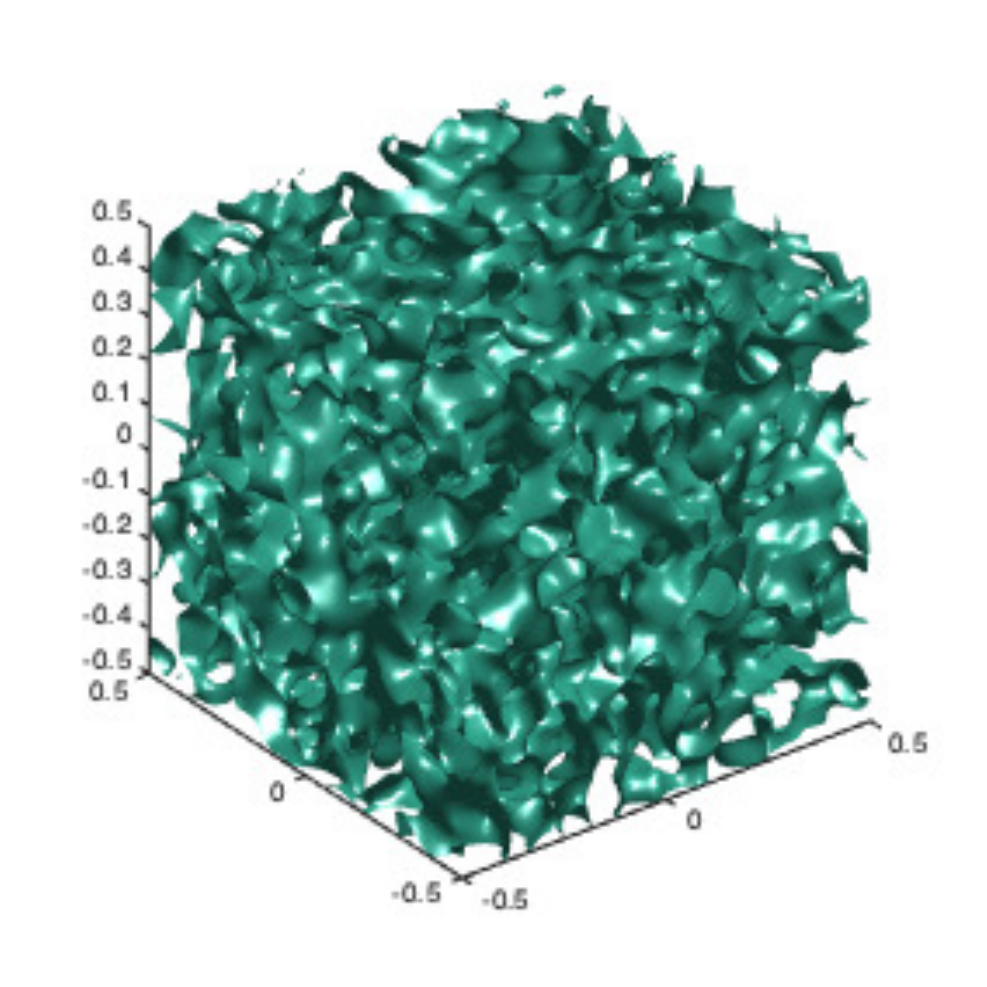}}\hspace{-0.1cm}
	{\includegraphics[width=0.45\textwidth]{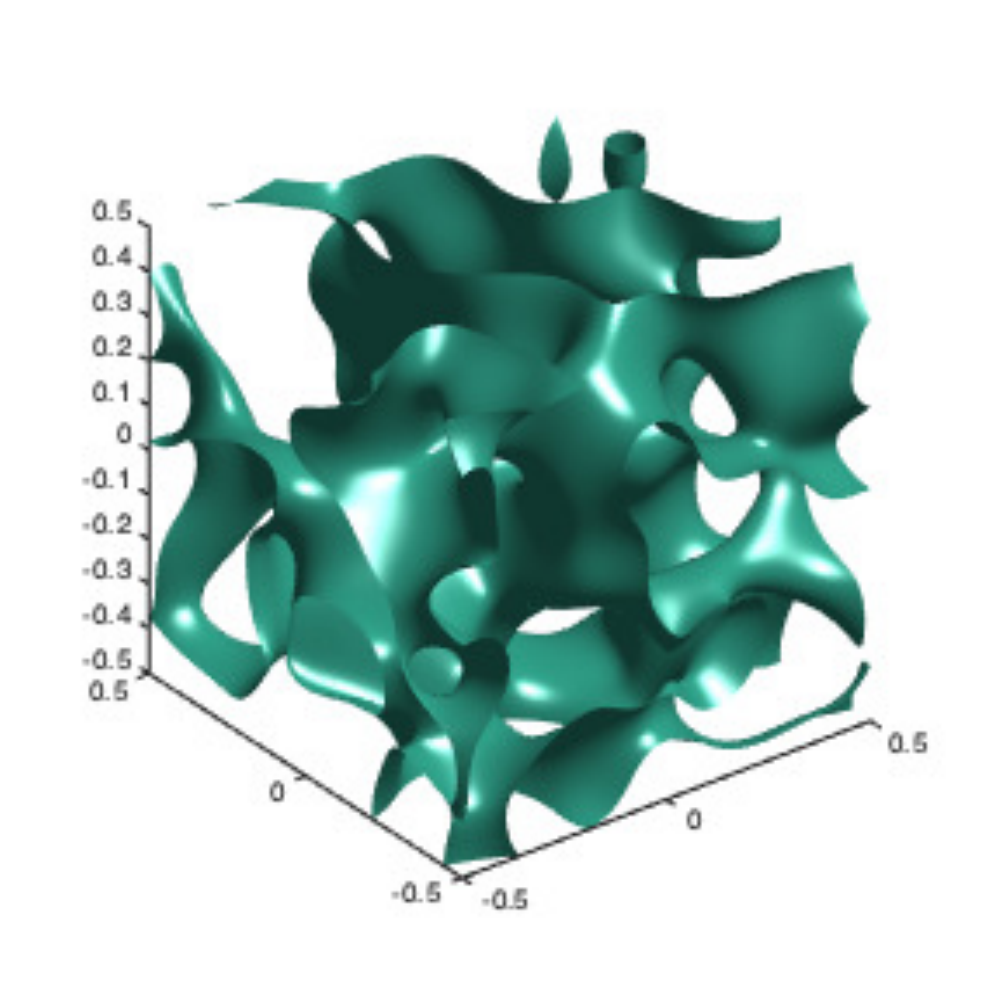}}\\ \vspace{-0.5cm}
	{\includegraphics[width=0.45\textwidth]{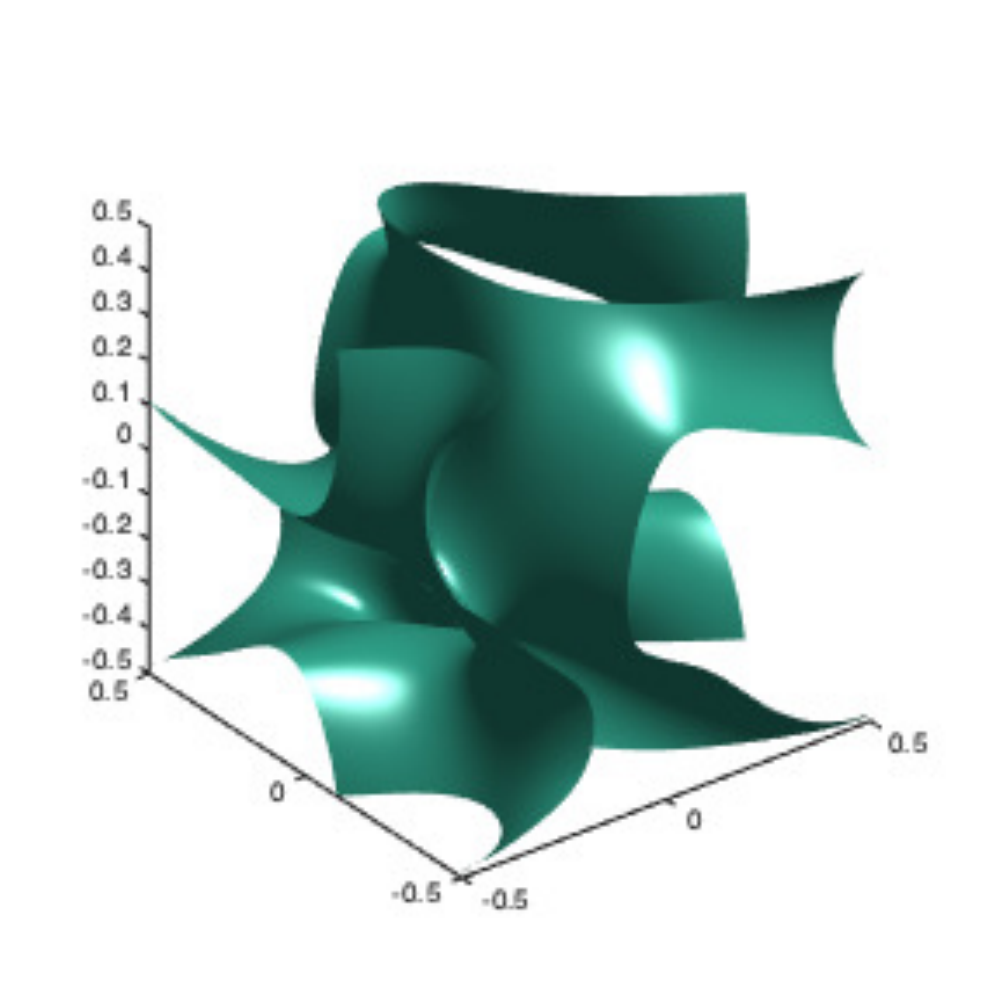}}\hspace{-0.1cm}
	{\includegraphics[width=0.45\textwidth]{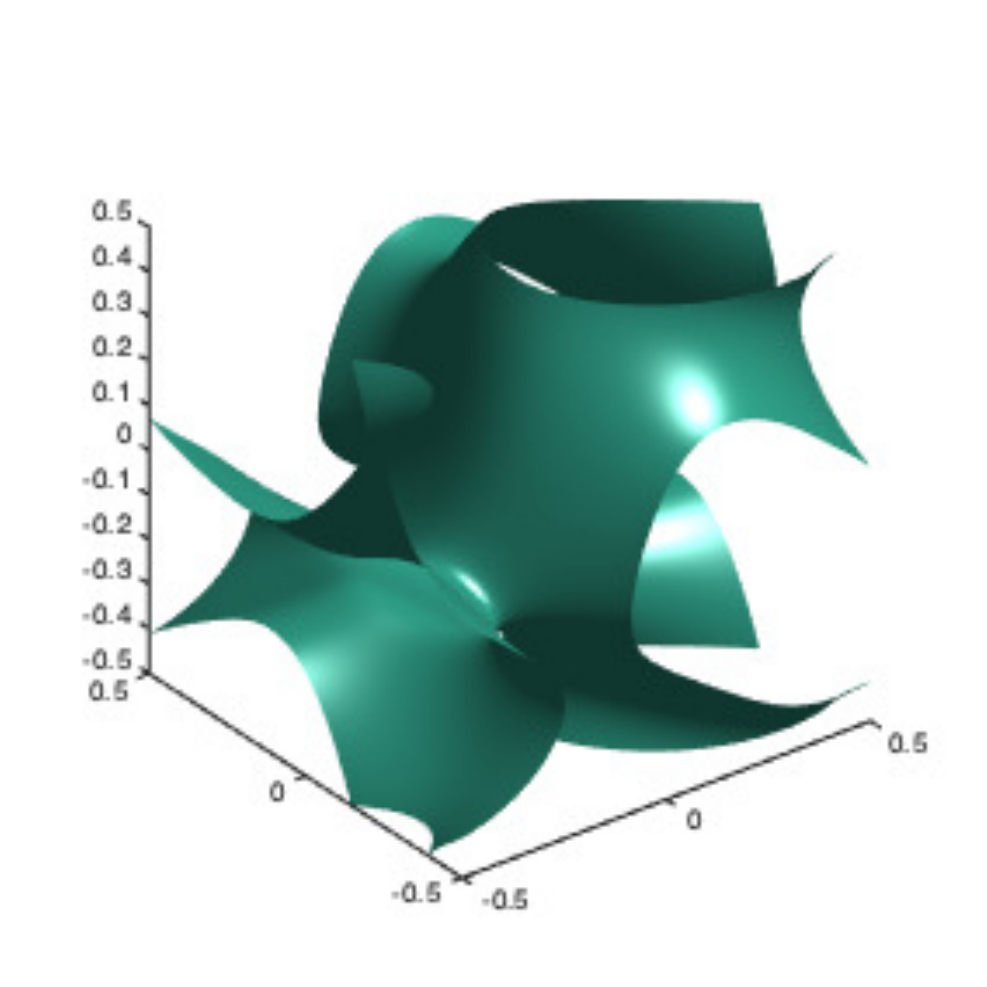}}\vspace{-0.2cm}
	\caption{The simulated phase structures at $t$ = 0.1,1,5,8, respectively (top to bottom and left to right) for the 3D convective Allen-Cahn equation with the Flory-Huggins potential.}
	\label{fig-fh3D-phase}
\end{figure}

\begin{figure}[!ht]
	\centering
	{\includegraphics[width=0.45\textwidth]{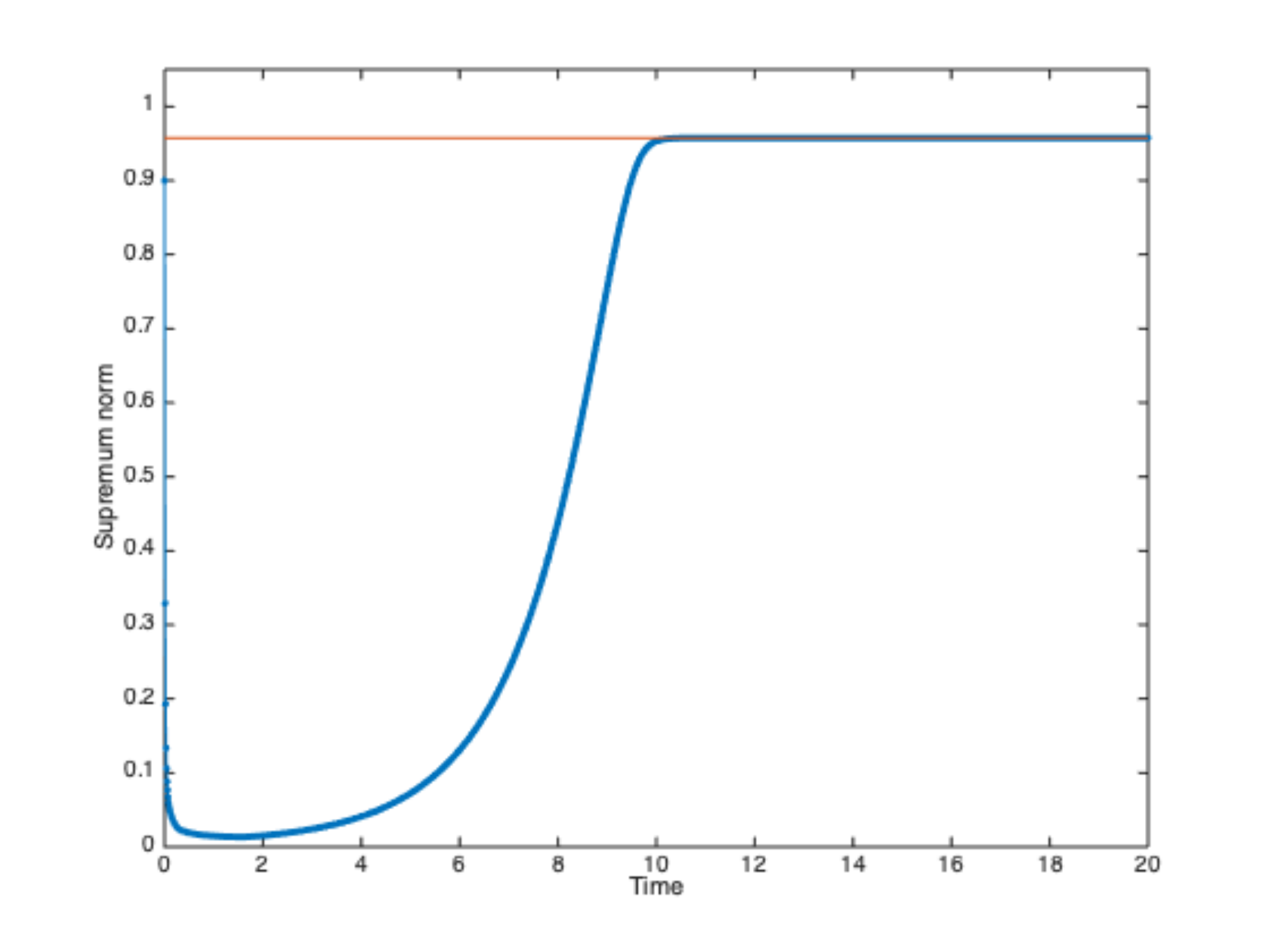}}\hspace{-0.3cm}
	{\includegraphics[width=0.45\textwidth]{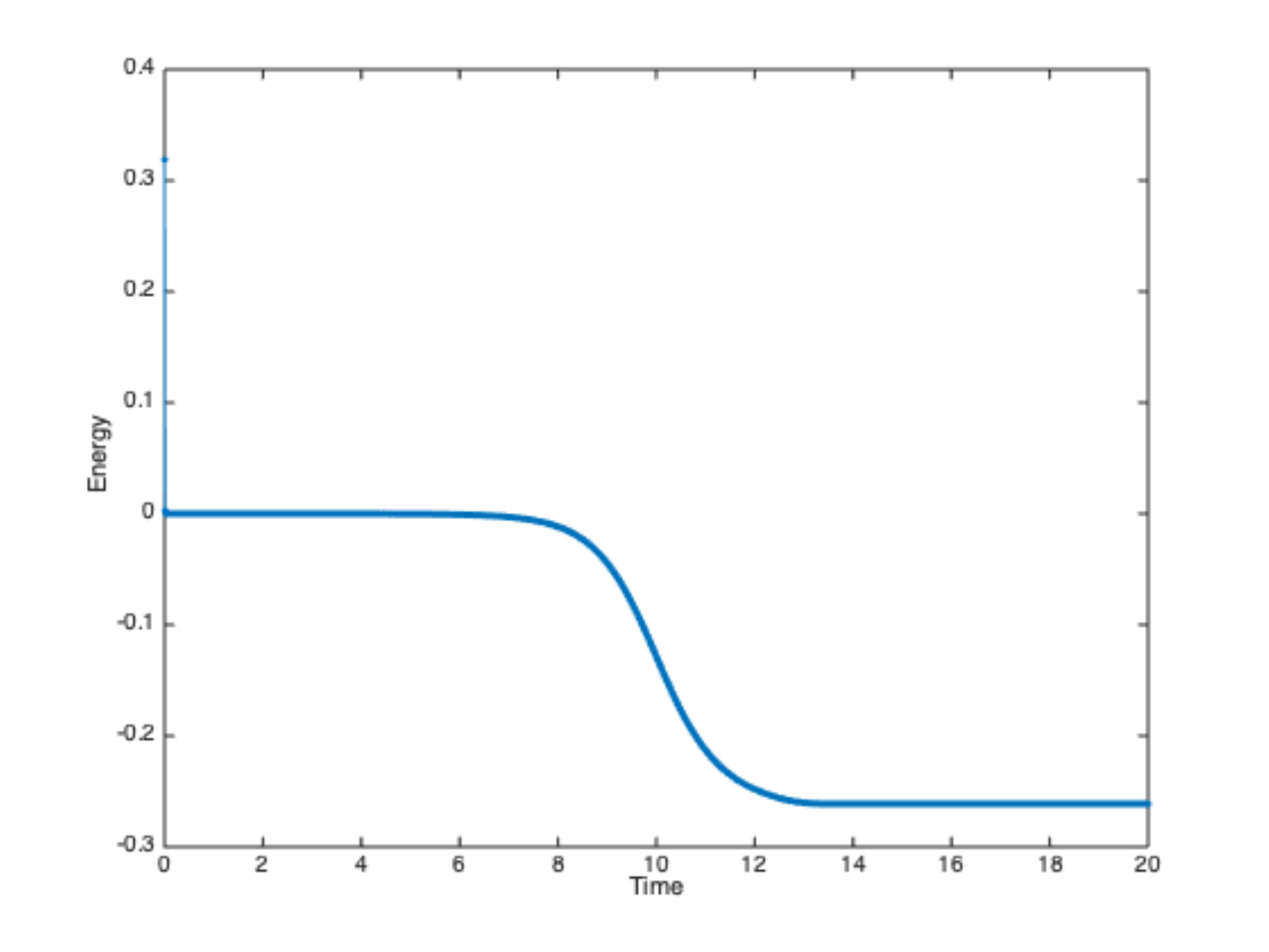}}
	\caption{The evolutions of the supremum norm (left) and the energy (right) of numerical solutions for the 3D convective Allen-Cahn equation with the Flory-Huggins potential.}
	\label{fig-fh3D-phy}
\end{figure}

Next, we test the Flory-Huggins potential case \eqref{fh} with $\theta=0.8$ and $\theta_c=1.6$ and the corresponding stabilizing coefficient is $\kappa=8.02$ \cite{DJLQ21}. Fig. \ref{fig-fh3D-phase} depicts the phase structures  of the numerical solutions at  $t$=0.1, 1, 5, 8, respectively. The time evolutions of the supremum norm and the energy of numerical solutions are shown in Fig. \ref{fig-fh3D-phy}. It is again observed that the MBP for the convective Allen-Cahn equation is numerically preserved perfectly and the energy decays monotonically as the double-well potential case.

\section{Conclusion}\label{con}
In this paper, we studied numerical solutions to the convective Allen-Cahn equation, including  nonlinear  mobility,  convective term, and nonlinear reaction. We proposed and analyzed stabilized ETD1 and ETDRK2  schemes which are linear and preserve the discrete MBP unconditionally. Various numerical examples are also  tested to verify our theoretical results on MBP and error estimates of the proposed schemes. 
It is worth noting that the upwind difference scheme used for discretizing the convective term in our method is only first-order accuracy in space. 
Many existing numerical results have indicated that the low-order spatial discrete schemes may not be able to capture the phase interface evolution well for the very small interface coefficient $\epsilon\ll1$. Thus, it is highly desirable to design high-order spatial approximation schemes  with the unconditional MBP preservation to simulate the convective Allen-Cahn equation, which remains as one of our future works.
In addition to the ETD schemes, the integrating factor method is also a widely-used time integration method. If the convective and nonlinear mobility terms are treated explicitly, the standard Runge-Kutta method could be applied. Therefore, designing accurate and efficient integrating factor Runge-Kutta (IFRK) schemes with  conditional or unconditional MBP preservation is also a very interesting topic for future study.



%
%

\end{document}